\documentclass[reqno]{amsart}
\usepackage{amsmath,amsfonts,amssymb,xcolor,enumitem,upgreek,graphicx}
\usepackage{caption}
\usepackage{subcaption}
\usepackage[english]{babel}
\usepackage[colorlinks=true,linkcolor=red!60!black,citecolor=green!60!black]{hyperref} 
\usepackage{mathrsfs}
\usepackage[colorlinks=true]{hyperref}

\newtheorem{theorem}{Theorem}[section]
\newtheorem{lemma}[theorem]{Lemma}
\newtheorem{proposition}[theorem]{Proposition}
\newtheorem{corollary}[theorem]{Corollary}

\theoremstyle{definition}
\newtheorem{definition}[theorem]{Definition}

\newtheorem{remark}[theorem]{Remark}
\newtheorem{remarks}[theorem]{Remarks}

\newcommand{\w}{\omega}

\numberwithin{equation}{section}

\begin{document}
\title[D-Concave Tipping]{Critical Transitions in D-Concave Nonautonomous \\ Scalar Ordinary Differential Equations \\ Appearing in Population Dynamics}
\author[J. Due\~{n}as]{Jes\'{u}s Due\~{n}as}
\author[C. N\'{u}\~{n}ez]{Carmen N\'{u}\~{n}ez}
\author[R. Obaya]{Rafael Obaya}
\address{Departamento de Matem\'{a}tica Aplicada, Universidad de Va\-lladolid, Paseo Prado de la Magdalena 3-5, 47011 Valladolid, Spain. The authors are members of IMUVA: Instituto de Investigaci\'{o}n en Matem\'{a}ticas, Universidad de Valladolid.}
\email[J.~Due\~{n}as]{jesus.duenas@uva.es}
\email[C.~N\'{u}\~{n}ez]{carmen.nunez@uva.es}
\email[R.~Obaya]{rafael.obaya@uva.es}
\thanks{All the authors were supported by Ministerio de Ciencia, Innovaci\'{o}n y Universidades (Spain) under project PID2021-125446NB-I00 and by Universidad de Valladolid under project PIP-TCESC-2020. J. Due\~{n}as was supported by Ministerio de Universidades (Spain) under programme FPU20/01627.}
\date{}
\begin{abstract}
A function with finite asymptotic limits gives rise to a transition equation between a ``past system" and a ``future system".
This question is analyzed in the case of nonautonomous coercive nonlinear scalar ordinary differential equations with concave derivative with respect to the state variable.
The fundamental hypotheses is the existence of three hyperbolic solutions for the limit systems,
in which case the upper and lower ones are attractive.
All the global dynamical possibilities are described in terms of the internal dynamics of the pullback attractor:
cases of tracking of the two hyperbolic attractive solutions or lack of it (tipping) may arise.
This analysis,
made in the language of processes and also in terms of the skewproduct formulation of the problem,
includes cases of rate-induced critical transitions,
as well as cases of phase-induced and size-induced tipping.
The conclusions are applied in models of mathematical biology and population dynamics.
Rate-induced tracking phenomena causing extinction of a native species or invasion of a non-native one are described,
as well as population models affected by a Holling type III functional response to predation where tipping due to the changes in the size of the transition may occur.
In all these cases,
the appearance of a critical transition can be understood as a consequence of the strength of Allee effect.
\end{abstract}
\renewcommand{\subjclassname}{\textup{2020} Mathematics Subject Classification}

\maketitle
\section{Introduction}
Tipping points or critical transitions are relevant nonlinear phenomena that can be described as large, sudden and often irreversible changes in the state of a given system in response to small and slow changes in the external conditions or external inputs of the physical phenomenon.
Motivated by the presence of these phenomena in climate \cite{aajq,awvc,lenton1,schellnhuber}, ecology \cite{altw,scheffer2,scheffer1,vanselow}, biology \cite{hill1,nene} or finance \cite{may1,yukalov}, among other scientific areas, this theory has been a focus of increasing interest in recent years.
A time-dependent dynamical transition connecting a past limiting system to a future one, both of them with known and similar internal dynamics,
is frequently analyzed in the mathematical formulation of this theory.
This formulation pays special attention to the local or global attractors,
which (partially or globally) concentrate the long-term dynamics.
Thus, the focus is put in understanding the evolution of the local or global attractor of the past system through the transition one.
In many cases,
the transition provides a full connection between the past and future dynamics:
the corresponding attractor globally connects with those of the past and future as time decreases and increases, respectively.
This behavior, which ensures the stability and robustness of the model,
is frequently called tracking.
In other cases, the set of trajectories which connect with the attractor of the past as time decreases is no longer the attractor of the transition system or it has a very different shape,
and it does not connect with the whole attractor of the future as time increases.
This behavior is usually called tipping.

Most of the deterministic papers on this subject consider past and future systems given by autonomous differential equations, connected by an asymptotically autonomous transition (see \cite{aspw}, \cite{kiersjones}, \cite{keeffe1}, \cite{sebastian1}
).
These works introduce several notions and mechanisms of tipping; among them, rate-induced tipping and phase-induced tipping.
A version of this theory where the past, future and transition systems are nonautonomous is considered in \cite{lno1} and \cite{lno2} for scalar quadratic ordinary differential equations,
for which the global dynamics is in general determined by the presence or absence of an attractor-repeller pair of bounded solutions.
Suitable one-parametric families of transition equations are considered in these papers in order to analyze the occurrence of rate-induced tipping, phase-induced tipping and tipping by size of transition.
In all these cases,
the presence of a tipping point is rigorously shown to be consequence of a saddle-node bifurcation of the nonautonomous transition differential equation as the parameter varies.

The more recent works \cite{dno1} and \cite{dno2} provide a bifurcation theory for nonautonomous scalar ordinary differential equations $x'=h(t,x)$, where $h\colon\mathbb{R}\times\mathbb{R}\rightarrow\mathbb{R}$ has concave derivative with respect to the state variable: d-concave equations, for short.
Conditions giving rise to saddle-node, transcritical and pitchfork nonautonomous bifurcation patterns for the one-parametric perturbations $x'=h(t,x)+ \lambda$, $x'=h(t,x)+ \lambda x$ and $x'=h(t,x)+ \lambda x^2$,
for $\lambda \in \mathbb{R}$, are rigorously described.
As in the case of \cite{lno1} and \cite{lno2}, some of the existing dynamical scenarios are exclusively nonautonomous,
in the sense that they cannot appear in the analogous autonomous cases:
time-dependent mathematical modeling can describe more real-world phenomena.
D-concave equations appear in mathematical models in biology and ecology (see e.g. \cite{courchamp1}, \cite{murray1}),
and hence the theory developed in \cite{dno1} and \cite{dno2} turns out to be interesting for applications in these fields.

This paper concerns the global dynamics induced by a nonautonomous transition equation with d-concave nonlinearities.
We discuss mechanisms providing rate-induced tipping/tracking, phase-induced tipping/tracking or size-induced tipping/tracking, showing in all the cases that they arise as a consequence of a nonautonomous saddle-node bifurcation of the equation.
The (quite usual) bistability of the dynamics induced by this type of d-concave differential equations makes it possible to incorporate the Allee effect to the mathematical models of some biological systems (see \cite{boukal1}, \cite{courchamp1}).
We will see how the appearance of a critical transition can also be understood as a consequence of the strength of the Allee effect in these population models.
This question was initially considered in \cite{remo1}.
In this work, we explore more types of critical transitions and of mechanisms causing them.

Let us describe the structure of the paper.
Section~\ref{sec:preliminaries} begins with some preliminary results for
nonautonomous scalar ordinary differential equations $x'=h(t,x)$ given by $C^2$-admissible functions $h$.
We recall that hyperbolic solutions are uniformly exponentially stable at $+\infty$ or $-\infty$, and prove a result on their persistence which is suitable for our setting.
We also recall the notions of locally pullback attractive or repulsive solutions and of global pullback attractor, and check that the existence of this one (which is determined by the globally bounded solutions) is ensured by a classical condition on coercivity.

Section~\ref{sec:2tipping} prepares the way to talk about critical transitions for some parametric problems.
Our starting point is a map $f\colon\mathbb{R}\times\mathbb{R}\times\mathbb{R} \rightarrow \mathbb{R}$ such that the nonautonomous equations $x'=f(t,x,\gamma_-)$ and $x'=f(t,x,\gamma_+)$ describe a past system and a future system.
The conditions assumed on $f$ include $C^2$-admissibility, d-concavity and coercivity,
as well as the fundamental hypotheses of existence of three uniformly separated bounded solutions for the asymptotic systems.
The transition equation is $x'=f(t,x,\Gamma(t))$ for a continuous function $\Gamma\colon\mathbb{R}\rightarrow\mathbb{R}$ such that $\lim_{t\to\pm\infty}\Gamma(t)=\gamma_\pm$.
Our main goal in this section is to classify the dynamics of these transition equations in three types which exhaust all the possibilities,
which is done in Section \ref{subsec:2transition}:
\textsc{Case A}, which occurs when the equation has three hyperbolic solutions which are the unique three uniformly separated;
\textsc{Case B}, which occurs when it has exactly two uniformly separated solutions, being one of them the unique hyperbolic one; and
\textsc{Case C}, which occurs when there are exactly two hyperbolic solutions but no uniformly separated bounded solutions.
The dynamics of the transition equations and the relation of all its trajectories with those of the past and future systems are exhaustively described.
In particular, \textsc{Case A} (which is persistent under small suitable variations) is the tracking situation,
since the pullback attractor of the transition equation connects that of the past system to that of the future system.
In \textsc{Cases B} (highly unstable) and \textsc{C} (again stable), the pullback attractor of the transition equation connects that of the past system with only a part of that of the future system:
the global connection is lost, and hence we are in tipping situations.

The hypotheses on admissibility on the coefficient function $f$ allow us to consider the skewproduct flows defined on the hulls for the compact-open topology of $f_\Delta$ for $\Delta=\gamma_-,\gamma_+,\Gamma$, where $f_\Delta(t,x)=f(t,x,\Delta(t))$.
In each case, the hull $\Omega_\Delta$ is a compact metric space admitting a continuous (time shift) flow,
and the skewproduct flow is defined on $\Omega_\Delta\times\mathbb R$ by considering simultaneously all the solutions of all the equations given by elements of the hull.
The condition on coercivity ensures the existence of a global (classical) attractor,
which is also determined by the globally bounded orbits and hence closely related to the initial global pullback attractor.
And the condition on d-concavity allows us to use the previous results of \cite{dno1} and \cite{dno2}.
This skewproduct formulation is a fundamental tool to prove part of our main results, and it has interest by its own:
in Section \ref{subsec:2skewproduct}, we show that \textsc{Case A} is characterized by the continuous variation on $\Omega_\Gamma$ of the upper and lower delimiters of the global attractor,
while in the \textsc{Cases B} and \textsc{C} this variation is only semicontinuous.

The arising of critical transitions is finally considered in Section~\ref{sec:mechanisms} under the additional assumption of monotonicity of $\gamma\to f(t,x,\gamma)$.
First, we consider one-parametric families of equations $x'=f(t,x,\Gamma(ct))$ for $c \in C=(0,\infty)$ and $x'=f(t,x,\Gamma(t+c))$ for $c \in C=\mathbb{R}$.
Let $\mathcal R_f$ be the set of values $\gamma$ of the parameter for which $x'=f(t,x,\gamma)$ is strictly d-concave and has three hyperbolic solutions.
We show that, if $\mathrm{cls}(\Gamma(\mathbb{R}))\subset\mathcal{R}_f$,
then the transition equations are in \textsc{Case A} for all the values of $c\in C$.
If $\mathrm{cls}(\Gamma(\mathbb{R}))$ is not contained in $\mathcal{R}_f$,
we establish additional conditions which allow us to determine the dynamical case of the equation by the sign of $\varphi(c)$ for a continuous function $\varphi\colon C \rightarrow\mathbb{R}$,
which can therefore be understood as a bifurcation map:
its simple zeroes determine a jump of \textsc{Case A} to \textsc{Case C} (through \textsc{Case B}).
The analysis of size-induced tipping which we made in Section~\ref{sec:mechanisms} reduces to the particular case of equations $x'=h(t,x-d\,\Gamma(t))$ for monotone $C^1$ maps $\Gamma$ with $\lim_{t\rightarrow\pm\infty}\Gamma'(t)=0$,
for which we show the existence of two tipping values of the parameter which delimit the interval of all the values of $d$ such that the equation is in \textsc{Case A} (and the two half-lines of all the values of $d$ such that the equation is in \textsc{Case C}).

In Section~\ref{section:singlespeciesAlleeeffect},
we apply part of the conclusions of the previous sections to mathematical population models given by time-dependent d-concave scalar differential equations with three uniformly separated bounded solutions.
These solutions are hyperbolic: attractive the upper and lower ones, repulsive the middle one.
We identify a strong or weak Allee effect with the attractive or repulsive character of the hyperbolic solution representing the lowest steady (nonnegative) state: a sparse or extinct population.
Since the model is nonautonomous, the strength of the Allee effect must be evaluated in average.
Suitable examples of populations modelled by d-concave cubic polynomials show that, when the system is subject to strong Allee effect, migration can cause two different types of critical transition, which may be the cause of extinction of the native species of the habitat, as well as of the invasion of a non-native species.
The parametric variation is on the rate, and in these cases we observe that tipping takes place at low transition rates while tracking appears for higher rates: rate-induced tracking occurs.
We also analyze analogous population models affected by a Holling type III functional response to predation.
New suitable examples show that, when the Allee effect is strong and the transition is large enough, size-induced tipping may cause the extinction of the species.
However, the upper steady population may persists as size increases if the Allee effect is weak.
The results of this section show that the d-concave equations are a suitable framework for the inclusion of the Allee effect in population dynamics models.
Finally, the use of the skewproduct formulation allows us to describe more indicators of the strength of Allee effect, of ergodic nature.

\section{Preliminaries}\label{sec:preliminaries}
Let the function $h\colon\mathbb{R}\times\mathbb{R}\rightarrow\mathbb{R}$, $(t,x)\mapsto h(t,x)$ be differentiable with respect to the state variable $x$, and such that $h$ and $h_x$ are jointly continuous. We consider the nonautonomous scalar differential equation
\begin{equation}\label{eq:1basicnonautonomous}
x'=h(t,x)\,,
\end{equation}
and represent by $t\mapsto x(t,s,x_0)$ the maximal solution of the initial value problem given by $x(s)=x_0$, defined on an interval $\mathcal{I}_{s,x_0}=(\alpha_{s,x_0},\beta_{s,x_0})$, with $-\infty\leq\alpha_{s,x_0}<s<\beta_{s,x_0}\leq\infty$. Hence, $x(s,s,x_0)=x_0$ and
the map $x\colon\bigcup_{(s,x_0)\in\mathbb R^2}(\mathcal{I}_{s,x_0}\times\{(s,x_0)\})\to\mathbb R$ satisfies the cocycle property $x(t_1,t_2,x(t_2,t_3,x_0))=x(t_1,t_3,x_0)$ if the left term is defined.
Given a set $A\subseteq\mathbb{R}$, we denote $x(t,s,A)=\{x(t,s,x_0)\colon\, x_0\in A\}$.
A set $B\subset\mathbb{R}\times\mathbb{R}$ is said to be \emph{invariant} (under the action of \eqref{eq:1basicnonautonomous}) if it is composed by graphs of globally defined functions.
Or, equivalently, if $x(t,s,B_s)=B_t$ for all $s,t\in\mathbb{R}$, where $B_s=\{x_0\colon\,(s,x_0)\in B\}$.

A function $h\colon\mathbb{R}\times\mathbb{R}\rightarrow\mathbb{R}$ is $C^0$-\emph{admissible} if it is bounded and uniformly continuous on $\mathbb{R}\times K$ for every compact subset $K\subset\mathbb{R}$. If, in addition, the derivative with respect to the state variable $h_x$ (resp. $h_x$ and $h_{xx}$) exists and is admissible, then $h$ is $C^1$-\emph{admissible} (resp. $C^2$-\emph{admissible}). We represent by $C^{0,n}$ the set of $C^n$-\emph{admisible} functions for $n=0,1,2$.

\subsection{Exponential dichotomies and hyperbolic solutions}\label{subsec:1dichoitomieshyperbolicity} In this subsection, we present the notions of exponential dichotomy and hyperbolic solution, which are fundamental for the theory of critical transitions: tipping points can be described as points of loss of hyperbolicity of some solutions of the transition equation.

Let $a\colon\mathbb{R}\rightarrow\mathbb{R}$ be a continuous function. The scalar linear differential equation $x'=a(t)\,x$ has \emph{exponential dichotomy} on $\mathbb{R}$ if there exist $k\geq1$ and $\beta>0$ such that either
\begin{equation}\label{eq:1hurwitzinfty}
\exp\int_s^t a(r)\; dr\leq ke^{-\beta (t-s)}\text{ whenever }t\geq s
\end{equation}
or
\begin{equation}\label{eq:1hurwitz-infty}
\exp\int_s^t a(r)\; dr\leq ke^{\beta (t-s)}\text{ whenever }t\leq s\,.
\end{equation}
The linear equation $x'=a(t)\, x$ is \emph{Hurwitz at $+\infty$} in case \eqref{eq:1hurwitzinfty} and \emph{at $-\infty$} in \eqref{eq:1hurwitz-infty}.

Given a bounded solution $\tilde x\colon\mathbb{R}\rightarrow\mathbb{R}$ of \eqref{eq:1basicnonautonomous}, we shall say that it is a \emph{hyperbolic solution} of \eqref{eq:1basicnonautonomous} if its variational equation $z'=h_x(t,\tilde x(t))\, z$ has exponential dichotomy on $\mathbb{R}$.
The (non unique) pair $(k,\beta)$ determining this property is called a \emph{dichotomy constant pair} for $\tilde x$. If \eqref{eq:1hurwitzinfty} (resp. \eqref{eq:1hurwitz-infty}) holds, then $\tilde x$ is said to be \emph{attractive} (resp. \emph{repulsive}). Proposition \ref{prop:2hiperbolicas} justifies these terms, and Theorem \ref{th:hypcontinuation} shows the persistence of this property under small variations on the coefficient function $h$. Both results are classical (see e.g. Lemma~3.3 of \cite{aloo} or Theorem~3.8 of \cite{potz} for Theorem \ref{th:hypcontinuation}), but we include their proofs in our particular setting.

\begin{proposition}\label{prop:2hiperbolicas}
Let $h\in C^{0,2}$, let $\tilde x$ be a hyperbolic solution
of \eqref{eq:1basicnonautonomous}, and let $(k,\beta)$
be a dichotomy constant pair for $\tilde x$.
\begin{itemize}
\item[\rm (i)] If $\tilde x$ is
attractive, then, given any $\bar\beta\in(0,\beta)$,
there exists $\rho>0$ such that, if $s\in\mathbb R$ and
$|\tilde x(s)-x_0|\le\rho$,
then $x(t,s,x_0)$ is defined for any $t\ge s$, and
\[
 |\tilde x(t)-x(t,s,x_0)|\le k\,e^{-\bar\beta\,(t-s)}|\tilde x(s)-x_0|
 \qquad\text{for $t\ge s$}\,.
\]
\item[\rm (ii)] If $\tilde x$ is repulsive, then, given any $\bar\beta\in(0,\beta)$,
there exists $\rho>0$ such that, if $s\in\mathbb R$ and
$|\tilde x(s)-x_0|\le\rho$,
then $x(t,s,x_0)$ is defined for any $t\le s$ and
\[
 |\tilde x(t)-x(t,s,x_0)|\le k\,e^{\bar\beta\,(t-s)}|\tilde x(s)-x_0|
 \qquad\text{for $t\le s$}\,.
\]
\end{itemize}
\end{proposition}
\begin{proof}
According to Taylor's Theorem, the change of variables $x=y+\tilde x(t)$ takes \eqref{eq:1basicnonautonomous} to
$y'=h_x(t,\tilde x(t))\,y+r(t,y)$, with $\lim_{y\to 0}r(t,y)/y=0$ uniformly in $t\in\mathbb{R}$.
Let $t\mapsto y(t,s,y_0)$ be the solution of this transformed equation satisfying
$y(s,s,y_0)=y_0$. According to the First Approximation Theorem
(see Theorem~III.2.4 of \cite{hale} and its proof), if
$\bar\beta\in(0,\beta)$, then there exists $\rho>0$ such that if $|y_0|\le\rho$, then
$y(t,s,y_0)$ is defined and satisfies $|y(t,s,y_0)|
\le k\,e^{-\bar\beta\,(t-s)}|y_0|$ for $t\ge s$. The inequality
in (i) follows from this, and the proof of (ii) is analogous.
\end{proof}

Given $\rho>0$ and $h\in {C}^{0,1}(\mathbb{R}\times\mathbb{R},\mathbb{R})$, we define the seminorm
\begin{equation*}
\|h\|_{1,\rho}=\sup_{(t,x)\in\mathbb{R}\times[-\rho,\rho]} |h(t,x)|+\sup_{(t,x)\in\mathbb{R}\times[-\rho,\rho]} |h_x(t,x)|\,,
\end{equation*}
and, given $x\in C(\mathbb{R},\mathbb{R})$ we define $\|x\|_\infty=\sup_{t\in\mathbb{R}}|x(t)|$.

\begin{theorem}\label{th:hypcontinuation} Let $h\in C^{0,1}$, let $\tilde x_h$ be an attractive (resp. repulsive) hyperbolic solution of \eqref{eq:1basicnonautonomous} with dichotomy constant pair $(k_0,\beta_0)$, and let $\rho>0$ satisfy $\|\tilde x_h\|_\infty<\rho$. Then, for every $\beta\in(0,\beta_0)$ and $\epsilon>0$, there exists $\delta_\epsilon>0$ such that, if $g\in C^{0,1}$ satisfies $\|h-g\|_{1,\rho}<\delta_\epsilon$, then there exists an attractive (resp. repulsive) hyperbolic solution $\tilde x_g$ of $x'=g(t,x)$ with dichotomy constant pair $(k_0,\beta)$ which satisfies $\|\tilde x_h-\tilde x_g\|_\infty<\epsilon$.
\end{theorem}
\begin{proof} Let us define $\delta_0=(\beta_0-\beta)/k_0$. Recall that Lecture~3 of \cite{coppel1} ensures that every continuous map $a\colon\mathbb{R}\rightarrow\mathbb{R}$ with $\|h_x({\cdot},\tilde x_h({\cdot}))-a({\cdot})\|_\infty<\delta_0$ determines a new Hurwitz equation, being $(k_0,\beta)$ a dichotomy constant pair for its hyperbolic solution 0. This fact will be used at the end of the proof.

We will work in the case in which $\tilde x_h$ is hyperbolic attractive. The change of variables $x=\tilde x_h(t)+y$ takes the equation
$x'=g(t,x)$ to
\begin{equation*}
y'=h_x(t,\tilde x_h(t))y+r_g(t,y)\,,
\end{equation*}
where $r_g(t,y)=g(t,\tilde x_h(t)+y)-h(t,\tilde x_h(t))-h_x(t,\tilde x_h(t))y$. Let $\epsilon_0\in(0,\rho]$ satisfy $\|\tilde x_h\|_\infty\leq \rho-\epsilon_0$. Since $h\in C^{0,1}$, there exists $\epsilon\in(0,\min\{1,\epsilon_0\})$ such that $|h_x(t,\tilde x_h(t)+y)-h_x(t,\tilde x_h(t))|\leq\delta_0/4$ for every $y\in [-\epsilon,\epsilon]$ and $t\in\mathbb{R}$.
Consequently, for any $g\in C^{0,1}$ such that $\|h-g\|_{1,\rho}\leq\delta_0/4$,
\begin{equation}\label{eq:2cita}
\begin{split}
&|g_x(t,\tilde x_h(t)+y)-h_x(t,\tilde x_h(t))| \\
&\qquad \leq |g_x(t,\tilde x_h(t)+y)-h_x(t,\tilde x_h(t)+y)|+|h_x(t,\tilde x_h(t)+y)-h_x(t,\tilde x_h(t))|\\
&\qquad \leq\|h-g\|_{1,\rho}+\frac{\delta_0}{4}\leq\frac{\delta_0}{2}
\end{split}
\end{equation}
for every $y\in[-\epsilon,\epsilon]$ and $t\in\mathbb{R}$. Notice that $r_g$ can be rewritten as
\begin{equation*}
r_g(t,y)=y\int_0^1 \big(g_x(t,\tilde x_h(t)+s y)-h_x(t,\tilde x_h(t))\big)\, ds+ g(t,\tilde x_h(t))-h(t,\tilde x_h(t))\,,
\end{equation*}
so \eqref{eq:2cita} ensures that $|r_g(t,y)|\leq \delta_0\epsilon/2+\|h-g\|_{1,\rho}\leq 3\delta_0\epsilon/4$ for every $y\in[-\epsilon,\epsilon]$, $t\in\mathbb{R}$ and $\|h-g\|_{1,\rho}\leq \delta_0\epsilon/4$. Since
\begin{equation*}
r_g(t,y_1)-r_g(t,y_2)=(y_1-y_2)\int_0^1 \big(g_x(t,\tilde x_h(t)+s y_1+(1-s)y_2)-h_x(t,\tilde x_h(t))\big)\, ds\,,
\end{equation*}
\eqref{eq:2cita} also ensures that $|r_g(t,y_1)-r_g(t,y_2)|\leq (\delta_0/2)\,|y_1-y_2|$ for every $y_1,y_2\in[-\epsilon,\epsilon]$, $t\in\mathbb{R}$ and $\|h-g\|_{1,\rho}\leq \delta_0/4$. Let us take $g\in C^{0,1}$ with $\|h-g\|_{1,\rho}\leq \delta_0\epsilon/4\le
\delta_0/4$. The results of Lecture~3 of \cite{coppel1} ensure that, for any $y_0\in C(\mathbb{R},\mathbb{R})$ with $\|y_0\|_\infty\leq \epsilon$, there is a unique bounded solution $Ty_0$ of $y'=h_x(t,\tilde x_h(t))y+r_g(t,y_0(t))$, given by
\begin{equation*}
Ty_0(t)=\int_{-\infty}^t u(t)u^{-1}(s)r_g(s,y_0(s))\, ds\,,
\end{equation*}
where $u(t)=\exp\int_0^t h_x(s,\tilde x_h(s))\, ds$. Therefore, $\|Ty_0\|_\infty\leq(3\epsilon/4)(\delta_0k_0/\beta)<3\epsilon/4$ and $\|Ty_1-Ty_2\|_\infty\leq(1/2)(\delta_0k_0/\beta)\, \|y_1-y_2\|_\infty<(1/2)\, \|y_1-y_2\|_\infty$. The map $T\colon C(\mathbb{R},[-\epsilon,\epsilon])\rightarrow C(\mathbb{R},[-\epsilon,\epsilon])$ is a contraction and thus it has a unique fixed point $y_g$.
It follows easily that $\tilde x_g=\tilde x_h+y_g$ is a bounded solution of $x'=g(t,x)$, and it satisfies $\|\tilde x_h-\tilde x_g\|_\infty\leq\epsilon$. The bound \eqref{eq:2cita} and the choice of $\delta_0$ at the beginning of the proof ensure that $\tilde x_g$ is an attractive hyperbolic solution of $x'=g(t,x)$ with dichotomy constant pair $(k_0,\beta)$. The proof is analogous in the repulsive case.
\end{proof}
\subsection{Pullback attraction and repulsion}\label{subsec:1pullback}
A solution $\bar x\colon(-\infty,\beta)\rightarrow\mathbb{R}$ of \eqref{eq:1basicnonautonomous} is said to be \emph{locally pullback attractive} if there exists  $s_0<\beta$ and $\delta>0$ such that, if $s\leq s_0$ and $|x_0-\bar x(s)|<\delta$, then $x(t,s,x_0)$ is defined for $t\in[s,s_0]$ and
\begin{equation*}
\lim_{s\rightarrow-\infty} \max_{x_0\in[\bar x(s)-\delta,\bar x(s)+\delta]}|\bar x(t)-x(t,s,x_0)|=0\qquad\text{for all }t\leq s_0\,.
\end{equation*}
In our scalar case, this definition is equivalent to the existence of  $s_0<\beta$ and $\delta>0$ such that, if $s\leq s_0$, then $x(t,s,\bar x(s)\pm\delta)$ is defined for $t\in[s,s_0]$, and
\begin{equation*}
\lim_{s\rightarrow-\infty}|\bar x(t)-x(t,s,\bar x(s)\pm\delta)|=0\qquad\text{for all }t\leq s_0\,.
\end{equation*}
Analogously, a solution $\bar x\colon(\alpha,\infty)\rightarrow\mathbb{R}$ of \eqref{eq:1basicnonautonomous} is said to be \emph{locally pullback repulsive} if and only if the solution $\bar y(t)=\bar x(-t)$ of $y'=-h(-t,y)$ is locally pullback attractive; or, equivalently, if
there exists  $s_0>\alpha$ and $\delta>0$ such that, if $s_0\leq s$ and $|x_0-\bar x(s)|<\delta$, then $x(t,s,x_0)$ is defined for $t\in[s_0,s]$ and
\begin{equation*}
\lim_{s\rightarrow \infty}|\bar x(t)-x(t,s,\bar x(s)\pm\delta)|=0\qquad\text{for all }t\geq s_0\,.
\end{equation*}

The next definition requires the notion of \emph{Hausdorff semidistance} between subsets $C_1$ and $C_2$ of $\mathbb{R}$,
\begin{equation*}
\text{dist}(C_1,C_2)=\sup_{x_1\in C_1}\left(\inf_{x_2\in C_2}|x_1-x_2|\right)\,.
\end{equation*}
\begin{definition}\label{def:1pullbackattractor} A family $\mathscr{A}=\{\mathscr{A}(t)\colon\, t\in\mathbb{R}\}$ is the \emph{pullback attractor} of \eqref{eq:1basicnonautonomous} if
\begin{enumerate}[label=\rm{(\roman*)}]
\item $\mathscr{A}(t)$ is a compact subset of $\mathbb{R}$ for each $t\in\mathbb{R}$;
\item $\mathscr{A}$ is invariant for \eqref{eq:1basicnonautonomous}, i.e., $\mathscr{A}(t)=x(t,s,\mathscr{A}(s))$ for all $s,t\in\mathbb R$;
\item $\mathscr{A}$ pullback attracts bounded subsets of $\mathbb{R}$, that is, for any bounded $D\subset\mathbb{R}$ and any $t\in\mathbb{R}$,
\begin{equation*}
\lim_{s\rightarrow-\infty} \mathrm{dist}(x(t,s,D),\mathscr{A}(t))=0\,;
\end{equation*}
\item $\mathscr{A}$ is the minimal family of closed sets with property $\mathrm{(iii)}$.
\end{enumerate}
\end{definition}
The pullback attractor $\mathscr{A}$ is said to be \emph{globally forward attractive} if, for every bounded subset $D\subset\mathbb{R}$ and every $s\in\mathbb{R}$,
\begin{equation}\label{eq:forwardattraction}
\lim_{t\rightarrow\infty} \mathrm{dist}(x(t,s,D),\mathscr{A}(t))=0\,,
\end{equation}
and it is said to be \emph{locally forward attractive} if there exists $\delta>0$ such that \eqref{eq:forwardattraction} holds for every $s\in\mathbb{R}$ and every bounded subset $D\subseteq \{x_1+x_2\colon$ exists $s\in\mathbb{R}$ such that $x_1\in\mathscr{A}(s)$ and $x_2\in[-\delta,\delta]\}$.
We recall that, in general, the pullback attraction property of Definition~\ref{def:1pullbackattractor} does not imply local forward attraction
(see e.g.~\cite{laos}). However, sometimes, our pullback attractors will also be locally or globally forward attractive.

The following proposition describes a dissipation hypothesis which ensures that all the solutions of the equation are globally forward defined and bounded.

\begin{proposition}\label{prop:1forwarddefinedcoercive} Let $h\in C^{0,1}$ satisfy the coercivity property $\lim_{|x|\rightarrow\infty} h(t,x)/x=-\infty$ uniformly in $t\in\mathbb{R}$. Then, all the maximal solutions of \eqref{eq:1basicnonautonomous} are globally forward defined and bounded, and the equation has a bounded pullback attractor $\mathscr{A}=\{\mathscr{A}(t)\colon\, t\in\mathbb{R}\}$, with $\mathscr{A}(t)$ composed by the values at $t$ of all the globally bounded solutions.
\end{proposition}
\begin{proof} Let us take $\rho>0$ with $h(t,x)/x\leq -1$ if $|x|\geq\rho$ and $t\in\mathbb{R}$. If $y(t)$ solves \eqref{eq:1basicnonautonomous}, then $y'(t)\leq-\rho<0$ if $y(t)\geq\rho$ and $y'(t)\geq\rho>0$ if $y(t)\leq-\rho$. Hence, first, $x(t,s,[-\rho,\rho])\subset[-\rho,\rho]$ for all $t\geq s$; second, if $x_0>\rho$, then $x(t,s,x_0)=x_0+\int_s^t h(l,x(l,s,x_0))\, dl\leq x_0-\rho(t-s)$ for those values of $t\geq s$ such that $x(t,s,x_0)\geq \rho$, which ensures the existence of $t_s\leq s+(x_0-\rho)/\rho$ such that $x(t,s,x_0)=\rho$; and third, analogously, if $x_0<-\rho$, then there exists $t_s\leq s-(x_0+\rho)/\rho$ such that $x(t,s,x_0)=-\rho$. This yields $\lim_{s\rightarrow-\infty}\mathrm{dist}(x(t,s,D),[-\rho,\rho])=0$ for all bounded set $D\subset\mathbb{R}$: $\mathcal{B}(t)=[-\rho,\rho]$ pullback attracts bounded sets in time $t$ (see Definition~1.11 of \cite{carvalho1}). Therefore, Theorem~2.12 of \cite{carvalho1} ensures the existence of the pullback attractor $\mathscr{A}=\{\mathscr{A}(t)\colon\, t\in\mathbb{R}\}$ with $\mathscr{A}(t)\subseteq[-\rho,\rho]$ for all $t\in\mathbb{R}$. The last assertion follows from Corollary~1.18 of \cite{carvalho1}, since the pullback attractor is bounded.
\end{proof}
\section{Tipping in d-concave equations}\label{sec:2tipping}
In this section, we will study scalar ordinary differential equations of the form
\begin{equation}\label{eq:2completeproblem}
x'=f_\Gamma(t,x)\,,
\end{equation}
for $f_\Gamma(t,x)=f(t,x,\Gamma(t))$, where $f\colon\mathbb{R}\times\mathbb{R}\times\mathbb{R}\rightarrow\mathbb{R}$ and $\Gamma\colon\mathbb{R}\rightarrow\mathbb{R}$ satisfy some of the following conditions (which we will specify at each statement):
\begin{enumerate}
\item[\hypertarget{h1}{\rm{\textbf{h1}}$_{\,~}$}] $\Gamma$ is continuous and has finite asymptotic limits $\gamma_\pm=\lim_{t\rightarrow\pm\infty}\Gamma(t)$.
\item[\hypertarget{h2}{\rm{\textbf{h2}}$_{\,~}$}] $f$ is continuous, its derivatives $f_x$ and $f_{xx}$ with respect to \emph{state variable} $x$ exist and are jointly continuous on $\mathbb{R}\times\mathbb{R}\times\mathbb{R}$, and the restrictions $f,f_x,f_{xx}\colon\mathbb{R}\times K_1\times K_2\rightarrow\mathbb{R}$ are bounded and uniformly continuous whenever $K_1$ and $K_2$ are compact subsets of $\mathbb{R}$.
\item[\hypertarget{h3}{\rm{\textbf{h3}}$_{\,~}$}] $f$ is \emph{coercive}, that is,
\begin{equation*}
\lim_{|x|\rightarrow\infty} \frac{f(t,x,\gamma)}{x}=-\infty
\end{equation*}
uniformly in $(t,\gamma)\in\mathbb{R}\times K$, for any compact set $K\subset\mathbb{R}$.
\item[\hypertarget{h4}{\rm{\textbf{h4}}$_{\,~}$}] $f$ is \emph{d-concave}, that is, its derivative $x\mapsto f_x(t,x,\gamma)$ is concave for all $(t,\gamma)\in\mathbb{R}\times \mathbb{R}$.
\item[\hypertarget{h5}{\rm{\textbf{h5}}$_\gamma$}] $f$ is $\gamma$-\emph{strictly concave}, that is,
\begin{equation}\label{eq:2strictdconcavityrfxx1}
\inf_{t\in\mathbb{R}}\big(f_{xx}(t,x_1,\gamma)-f_{xx}(t,x_2,\gamma)\big)>0
\end{equation}
whenever $x_1<x_2$.
\item[\hypertarget{h6}{\rm{\textbf{h6}}$_\gamma$}] $x'=f(t,x,\gamma)$ has three hyperbolic solutions.
\end{enumerate}
Hypothesis \hyperlink{h1}{\textbf{h1}} allows us to understand the equation \eqref{eq:2completeproblem} as a transition between
the \emph{past equation} and the \emph{future equation}, respectively given by
\begin{equation}\label{eq:2pastproblem}
x'=f_{\gamma_-}(t,x)\,,
\end{equation}
\begin{equation}\label{eq:2futureproblem}
x'=f_{\gamma_+}(t,x)\,,
\end{equation}
where $f_{\gamma_\pm}(t,x)=f(t,x,\gamma_\pm)$;
and to pose the question if also the dynamics of \eqref{eq:2completeproblem} is, to a certain extent, a transition between the dynamics of \eqref{eq:2pastproblem} and \eqref{eq:2futureproblem}.

The most simple example of function $f$ satisfying \hyperlink{h2}{\textbf{h2}}-\hyperlink{h4}{\textbf{h4}}
is a cubic polynomial $f(t,x,\gamma)=-a_3(t,\gamma)x^3+a_2(t,\gamma) x^2+a_1(t,\gamma)x+a_0(t,\gamma)$ such that
its coefficients are bounded and uniformly continuous on $\mathbb{R}\times K$ for any compact set $K\subset\mathbb{R}$,
and such that there exists $\delta>0$ satisfying $a_3(t,\gamma)\ge\delta$ for all $(t,\gamma)\in\mathbb R\times\mathbb{R}$.
We remark that, in this case, \hyperlink{h5}{\textbf{h5}$_\gamma$} holds for any $\gamma\in\mathbb{R}$, as $\inf_{t\in\mathbb R}(f_{xx}(t,x_1,\gamma)-f_{xx}(t,x_2,\gamma))\ge 6\delta(x_2-x_1)$ for every $x_1<x_2$.

\begin{remarks}\label{rm:twovariables}
\hypertarget{311}{1}.~All the pairs of maps $(f,\Gamma)$ giving rise to equations \eqref{eq:2completeproblem} will be assumed to satisfy
\hyperlink{h1}{\textbf{h1}} ($\Gamma$) and \hyperlink{h2}{\textbf{h2}} ($f$). For such a pair,
we shall say that \lq\lq $(f,\Gamma)$ satisfies \hyperlink{h5}{\textbf{h5}} (and \hyperlink{h6}{\textbf{h6}})" if \hyperlink{h5}{\textbf{h5}$_{\gamma_\pm}$} (and \hyperlink{h6}{\textbf{h6}$_{\gamma_\pm}$}) holds, where
$\gamma_\pm$ are the asymptotic limits of $\Gamma$ provided by \hyperlink{h1}{\textbf{h1}}.
When saying that $(f,\Gamma)$ satisfies \hyperlink{h3}{\textbf{h3}} or \hyperlink{h4}{\textbf{h4}},
we mean that $f$ does.
\par
\hypertarget{312}{2}.~Some of the results of this section and the next one
require to adapt the previous hypotheses to two-variable maps  $h\colon\mathbb{R}\times\mathbb{R}\rightarrow\mathbb{R}$. In this case, the list is:
\begin{enumerate}
\item[\hypertarget{h2t}{\rm{\textbf{h2$_*$}}\,}] $h$ is in $C^{0,2}$ (see Section~\ref{sec:preliminaries}),
\item[\hypertarget{h3t}{\rm{\textbf{h3$_*$}}\,}] $h$ is coercive, that is, $\lim_{|x|\to\infty}h(t,x)/x=-\infty$ uniformly in $t\in\mathbb R$,
\item[\hypertarget{h4t}{\rm{\textbf{h4$_*$}}\,}] $h$ is d-concave, that is, $x\mapsto h_x(t,x)$ is concave for all $t\in\mathbb{R}$,
\item[\hypertarget{h5t}{\rm{\textbf{h5$_*$}}\,}] $\inf_{t\in\mathbb R}(h_{xx}(t,x_1)-h_{xx}(t,x_2))>0$ for $x_1<x_2$,
\item[\hypertarget{h6t}{\rm{\textbf{h6$_*$}}\,}] $x'=h(t,x)$ has three hyperbolic solutions.
\end{enumerate}
Note that \hyperlink{h5t}{\textbf{h5$_*$}} is stronger than \hyperlink{h4t}{\textbf{h4$_*$}}.
\end{remarks}
It is easy to check that the map $f_\Gamma$ given by $f_\Gamma(t,x)=f(t,x,\Gamma(t))$ satisfies
\hyperlink{h2t}{\textbf{h2$_*$}}-\hyperlink{h4t}{\textbf{h4$_*$}}
if $\Gamma$ is continuous and
$f$ satisfies \hyperlink{h2}{\textbf{h2}}-\hyperlink{h4}{\textbf{h4}}.
If so, $f_\Gamma$ satisfies the hypotheses of
Proposition \ref{prop:1forwarddefinedcoercive}, which shows the boundedness as time increases
of all the solutions of \eqref{eq:2completeproblem} as well as the existence of
a bounded pullback attractor $\mathscr A_\Gamma=\{\mathscr A_\Gamma(s)\colon s\in\mathbb R\}$.
In particular, the set
\[
 \mathcal B_\Gamma=\{(s,x_0)\in\mathbb R\times\mathbb R\colon \,t\mapsto x_\Gamma(t,s,x_0) \text{ is
 bounded}\}
\]
of the initial data of bounded solutions is nonempty. Here,
$t\mapsto x_\Gamma(t,s,x_0)$ is the maximal solution of \eqref{eq:2completeproblem}
with $x_\Gamma(s,s,x_0)=x_0$. The boundedness of the pullback attractor ensures that
\begin{equation}\label{eq:3conjuntosolucionesacotadas}
 \mathcal B_\Gamma=\bigcup_{s\in\mathbb R}\big(\{s\}\times [l_\Gamma(s),u_\Gamma(s)]\big)\,,
\end{equation}
where $l_\Gamma$ and $u_\Gamma$ are the lower and upper bounded solutions of
\eqref{eq:2completeproblem}.
Note that $\mathcal{B}_\Gamma$ is invariant for \eqref{eq:2completeproblem}. Note also that $x_0\in [l_\Gamma(s),u_\Gamma(s)]$
if and only if $x_\Gamma(t,s,x_0)$ is globally defined and bounded for $t\le s$;
or, more precisely, that $x_\Gamma(t,s,x_0)$ is bounded from above on its domain if and only if
$x_0\le u_\Gamma(s)$ and from below if and only if $x_0\ge l_\Gamma(s)$. In addition,
\begin{remark}\label{remark:sectionsboundedsolutionsopullbackattractor}
Proposition~\ref{prop:1forwarddefinedcoercive} ensures that the section $(\mathcal B_\Gamma)_s=\{x_0\in\mathbb R
\colon\,(s,x_0)\in\mathcal B_\Gamma\}=[l_\Gamma(s),u_\Gamma(s)]$
coincides with $\mathscr A_\Gamma(s)$.
\end{remark}
We say that $n$ bounded solutions $x_1(t),x_2(t),\dots,x_n(t)$ of \eqref{eq:2completeproblem}, with $n\geq 2$,
are {\em uniformly separated\/} if $\inf_{t\in\mathbb R}|x_i(t)-x_j(t)|>0$ for any $1\leq i<j\leq n$. It is trivial
that $l_\Gamma$ and $u_\Gamma$ are uniformly separated in the case of
existence of (at least) two uniformly separated solutions, and that they are the upper and lower bounded solutions which may be uniformly separated from any other.
In particular, if \eqref{eq:2completeproblem} has three uniformly separated solutions and $m_\Gamma$ is the middle one, then $l_\Gamma<m_\Gamma<u_\Gamma$ are three uniformly separated solutions.

The next theorem establishes conditions ensuring the existence of at most three uniformly separated solutions,
in which case they are hyperbolic and determine
the global dynamics.
Its proof,
based on the results of \cite{dno1} and,
in turn,
in those of \cite{tineo1},
is postponed until Subsection~\ref{subsec:2skewproduct},
where the setting of \cite{dno1} is described:
the skewproduct formalism allowing us to define a flow from the solutions of the nonautonomous equation \eqref{eq:2completeproblem}.
\begin{theorem}\label{th:3three}
Let $(f,\Gamma)$ satisfy {\rm\hyperlink{h1}{\textbf{h1}}-\hyperlink{h5}{\textbf{h5}}}.
The following assertions are equivalent:
\begin{itemize}
\item[\rm (a)] \eqref{eq:2completeproblem} has three uniformly separated solutions
$l_\Gamma<m_\Gamma<u_\Gamma$.
\item[\rm (b)] \eqref{eq:2completeproblem} has three hyperbolic solutions
$\tilde l_\Gamma<\tilde m_\Gamma<\tilde u_\Gamma$.
\item[\rm (c)] \eqref{eq:2completeproblem} has three uniformly separated hyperbolic solutions
$\tilde l_\Gamma<\tilde m_\Gamma<\tilde u_\Gamma$.
\end{itemize}
In this case, $l_\Gamma=\tilde l_\Gamma$, $m_\Gamma=\tilde m_\Gamma$ and
$u_\Gamma=\tilde u_\Gamma$, so that $\mathcal B_\Gamma=\bigcup_{s\in\mathbb R}
(\{s\}\times [\tilde l_\Gamma(s),\tilde u_\Gamma(s)])$; and $\tilde l_\Gamma$ and
$\tilde u_\Gamma$ are attractive and $\tilde m_\Gamma$ is repulsive. In addition,
$\lim_{t\to\infty}(x_\Gamma(t,s,x_0)-\tilde u_\Gamma(t))=0$ if $x_0>\tilde m_\Gamma(s)$,
$\lim_{t\to\infty}(x_\Gamma(t,s,x_0)-\tilde l_\Gamma(t))=0$ if $x_0<\tilde m_\Gamma(s)$,
and $\lim_{t\to-\infty}(x_\Gamma(t,s,x_0)-\tilde m_\Gamma(t))=0$ if $x_0\in(\tilde l_\Gamma(s),
\tilde u_\Gamma(s))$.
In particular, \eqref{eq:2completeproblem} has at most three uniformly separated solutions.
\end{theorem}
Figure \ref{fig:CaseA} depicts this situation for a periodic nonautonomous equation.
\par
\begin{remark}\label{rm:3siempre}
It is important to point out that, as we will explain in Remark \ref{rm:3siempre-dos},
all the assertions of Theorem~\ref{th:3three} hold for the equation $x'=h(t,x)$ if $h\colon\mathbb{R}\times\mathbb{R}\rightarrow\mathbb{R}$ satisfies \hyperlink{h2t}{\textbf{h2}$_*$}-\hyperlink{h5t}{\textbf{h5}$_*$}.
In particular, if $(f,\Gamma)$ satisfies \rm{\hyperlink{h1}{\textbf{h1}}-\hyperlink{h5}{\textbf{h5}}}
(in the sense explained in Remark \ref{rm:twovariables}\hyperlink{311}{.1}), then Theorem \ref{th:3three} holds for the past and future equations \eqref{eq:2pastproblem} and
\eqref{eq:2futureproblem}, since the functions $(t,x)\mapsto f_{\gamma_\pm}(t,x)$ satisfy the required conditions.
\end{remark}
The proof of the following technical proposition, which will be used in the proof of Theorem~\ref{th:3general}, is also postponed until Section~\ref{subsec:2skewproduct}. This proposition explains that, if a coercive transition equation with three hyperbolic solutions has a strictly d-concave future equation (resp. past equation) with three hyperbolic solutions, then it inherits the essentials of the forward dynamics (resp. backward dynamics) of the future equation (resp. past equation). Recall that $f_{\gamma_\pm}(t,x)=f(t,x,\gamma_\pm)$.
\begin{proposition}\label{prop:3three}
Let $(f,\Gamma)$ satisfy {\rm\hyperlink{h1}{\textbf{h1}}-\hyperlink{h3}{\textbf{h3}}} and
assume that \eqref{eq:2completeproblem} has three uniformly separated hyperbolic solutions $\tilde l_\Gamma<\tilde m_\Gamma<\tilde u_\Gamma$.
Let $\tilde l_{\gamma_\pm}<\tilde m_{\gamma_\pm}<\tilde u_{\gamma_\pm}$ be the hyperbolic solutions of $x'=f_{\gamma_\pm}(t,x)$, if they exist.
\begin{enumerate}[label=\rm{(\roman*)}]
\item If $f_{\gamma_+}$ satisfies {\rm\hyperlink{h5t}{\textbf{h5}$_*$}}-{\rm\hyperlink{h6t}{\textbf{h6}$_*$}}, then $\lim_{t\to\infty}(x_\Gamma(t,s,x_0)-\tilde u_{\gamma_+}(t))=0$ for $x_0>\tilde m_\Gamma(s)$ and
$\lim_{t\to\infty}(x_\Gamma(t,s,x_0)-\tilde l_{\gamma_+}(t))=0$ for $x_0<\tilde m_\Gamma(s)$.
\item If $f_{\gamma_-}$ satisfies {\rm\hyperlink{h5t}{\textbf{h5}$_*$}}-{\rm\hyperlink{h6t}{\textbf{h6}$_*$}}, then
    $t\mapsto x_\Gamma(t,s,x_0)$ is bounded from above (resp. from below) as time decreases if and only if $x_0\leq \tilde u_\Gamma(s)$ (resp. $x_0\geq \tilde l_\Gamma(s)$);
    and $\lim_{t\to-\infty}(x_\Gamma(t,s,x_0)-\tilde m_{\gamma_-}(t))=0$ for $x_0\in(\tilde l_\Gamma(s), \tilde u_\Gamma(s))$.
\end{enumerate}
\end{proposition}
\subsection{The casuistic under a fundamental hypothesis: tracking or tipping}\label{subsec:2transition}
For the main results of this subsection,
we will assume that $(f,\Gamma)$ satisfies \hyperlink{h1}{\textbf{h1}}-\hyperlink{h6}{\textbf{h6}}
(see Remark \ref{rm:twovariables}\hyperlink{311}{.1} once again).
In this case, we represent the three hyperbolic solutions of \eqref{eq:2pastproblem} and \eqref{eq:2futureproblem} provided by \hyperlink{h6}{\textbf{h6}$_{\gamma_\pm}$}
by $\tilde l_{\gamma_\pm}<\tilde m_{\gamma_\pm}<
\tilde u_{\gamma_\pm}$, and observe that the dynamics induced by these equations is
that described by Theorem \ref{th:3three}:
see Remark \ref{rm:3siempre}.
In particular, $\tilde l_{\gamma_\pm}$ and $\tilde u_{\gamma_\pm}$
are attractive and $\tilde m_{\gamma_\pm}$ is repulsive.

\begin{definition}\label{def:3casesABC}
Let $(f,\Gamma)$ satisfy \hyperlink{h1}{\textbf{h1}}-\hyperlink{h6}{\textbf{h6}}.
We shall say that equation \eqref{eq:2completeproblem} is
\begin{itemize}
\item[-] in \textsc{Case} \hypertarget{A}{\textsc{A}} if it has three hyperbolic solutions, which are the unique three
uniformly separated solutions;
\item[-] in \textsc{Case} \hypertarget{B}{\textsc{B}} if it has exactly two uniformly separated solutions,
one of them being the only hyperbolic solution, of attractive type, and the other locally pullback attractive and repulsive;
\item[-] and in \textsc{Case} \hypertarget{C}{\textsc{C}} if it has no uniformly separated solutions and it has exactly two hyperbolic solutions,
which are attractive.
\end{itemize}
\end{definition}
Our goals in this subsection are to check that these cases exhaust the
possibilities when $(f,\Gamma)$ satisfies \hyperlink{h1}{\textbf{h1}}-\hyperlink{h6}{\textbf{h6}}, to describe the global
dynamics for each of them, and to relate these dynamics to those of
\eqref{eq:2pastproblem} and \eqref{eq:2futureproblem}. More precisely,
we will say that the graph of a solution $\bar x$ of \eqref{eq:2completeproblem} defined on a positive halfline (resp. negative halfline) {\em approaches that of a continuous map $c\colon\mathbb R\to\mathbb R$
as time increases} (resp. {\em as time decreases}) if $\lim_{t\to\infty}(\bar x(t)-c(t))=0$
(resp. $\lim_{t\to-\infty}(\bar x(t)-c(t))=0$). The next theorems describe the casuistic
for the global dynamics for \eqref{eq:2completeproblem} in terms of the approaching
behavior of the graphs of its solutions to the graphs of the hyperbolic
solutions of the limit equations \eqref{eq:2pastproblem} and \eqref{eq:2futureproblem}.
\begin{theorem}\label{th:3general}
Let $(f,\Gamma)$ satisfy {\rm\hyperlink{h1}{\textbf{h1}}-\hyperlink{h6}{\textbf{h6}}}. Then,
\begin{itemize}
\item[\rm(i)] the maximal solutions $u_\Gamma(t)$ and $l_\Gamma(t)$ of {\em \eqref{eq:2completeproblem}}
are the unique ones satisfying $\lim_{t\to-\infty}(u_\Gamma(t)-\tilde u_{\gamma_-}(t))=0$
and $\lim_{t\to-\infty}(l_\Gamma(t)-\tilde l_{\gamma_-}(t))=0$;
and $x_0\in(l_\Gamma(s),u_\Gamma(s))$ if and only if
$\lim_{t\to-\infty}(x_\Gamma(t,s,x_0)-\tilde m_{\gamma_-}(t))=0$.
\item[\rm(ii)] There exists a unique maximal solution
$m_\Gamma(t)$ satisfying $\lim_{t\to\infty}(m_\Gamma(t)-\tilde m_{\gamma_+}(t))=0$;
and, in addition, $\lim_{t\to\infty}(x_\Gamma(t,s,x_0)-\tilde u_{\gamma_+}(t))=0$
(resp. $\lim_{t\to\infty}(x_\Gamma(t,s,x_0)-\tilde l_{\gamma_+}(t))=0$) if and only if
$x_\Gamma(t,s,x_0)>m_\Gamma(t)$ (resp. $x_\Gamma(t,s,x_0)<m_\Gamma(t)$) for a value
of $t$ in the common interval of definition of both solutions, and
hence for all these values of $t$. In particular, given $(s,x_0)\in\mathbb R\times\mathbb R$,
one of these situations holds:
$\lim_{t\to\infty}(x_\Gamma(t,s,x_0)-\tilde u_{\gamma_+}(t))=0$,
$\lim_{t\to\infty}(x_\Gamma(t,s,x_0)-\tilde m_{\gamma_+}(t))=0$, or
$\lim_{t\to\infty}(x_\Gamma(t,s,x_0)-\tilde l_{\gamma_+}(t))=0$.
\item[\rm(iii)] $l_\Gamma$ and $u_\Gamma$ are locally pullback
attractive solutions of \eqref{eq:2completeproblem}, and $m_\Gamma$
is a locally pullback repulsive solution.
\end{itemize}
\end{theorem}
\begin{proof}
We take any $\epsilon>0$ with $\epsilon\le \epsilon_0=(1/3)\min\{\inf_{t\in\mathbb R}(\tilde u_{\gamma_{\pm}}(t)-\tilde m_{\gamma_{\pm}}(t)),$ $
\inf_{t\in\mathbb R}(\tilde m_{\gamma_{\pm}}(t)-\tilde l_{\gamma_{\pm}}(t))\}$, where $\tilde l_{\gamma_{\pm}}$,
$\tilde m_{\gamma_{\pm}}$ and $\tilde u_{\gamma_{\pm}}$ are the functions provided by \hyperlink{h6}{\textbf{h6}$_{\gamma_{\pm}}$}
(see Remark \ref{rm:twovariables}\hyperlink{311}{.1}).
We choose $\rho>0$ such that $f(t,x,\gamma)/x\leq -1$ if $(t,\gamma)\in\mathbb{R}\times\mathrm{cl}(\Gamma(\mathbb{R}))$ and $|x|\geq\rho$.
As seen in the proof of Proposition~\ref{prop:1forwarddefinedcoercive}, the sets
of bounded solutions of \eqref{eq:2completeproblem}, \eqref{eq:2pastproblem} and \eqref{eq:2futureproblem}
are contained in $\mathbb R\times[-\rho,\rho]$. Theorem~\ref{th:hypcontinuation}
provides $\delta_\epsilon>0$ such that, if $\|f_{\gamma_-}-g_-\|_{1,\rho}\le\delta_\epsilon$
and $\|f_{\gamma_+}-g_+\|_{1,\rho}\le\delta_\epsilon$,
then each one of the equations $x'=g_-(t,x)$ and $x'=g_+(t,x)$ has three
hyperbolic solutions which are respectively $\epsilon$-close to those of \eqref{eq:2pastproblem} and \eqref{eq:2futureproblem}, and
hence uniformly separated.
For each $r>0$, we define
\[
 \Gamma^-_r(t)=\left\{\begin{array}{ll}\Gamma(t)&\text{if $t\le-r$}\,,\\\Gamma(-r)&\text{if $t>-r$}\,,
 \end{array}\right.\qquad\text{and}\qquad
 \Gamma^+_r(t)=\left\{\begin{array}{ll}\Gamma(r)&\text{if $t< r$}\,,\\\Gamma(t)&\text{if $t\ge r$}\,.
 \end{array}\right.
\]
Condition \hyperlink{h1}{\textbf{h1}} yields $\lim_{r\rightarrow\infty}\|\gamma_\pm-\Gamma^\pm_r\|_\infty=0$.
The uniform continuity of $f$ and $f_x$ on $\mathbb R\times [-\rho,\rho]\times\mathrm{cl}(\Gamma(\mathbb R))$
provides $r_\epsilon\geq 1/\epsilon$ such that
$\|f_{\gamma_\pm}-f_{\Gamma^\pm_r}\|_{1,\rho}\le\delta_\epsilon$ for all $r\ge r_\epsilon$, where $f_{\Gamma_r^\pm}(t,x)=f(t,x,\Gamma_r^\pm(t))$.
For further purposes, we
observe that $f_{\Gamma_r^\pm}$ satisfy \hyperlink{h2}{\textbf{h2}$_*$}-\hyperlink{h4}{\textbf{h4}$_*$}.
An optimal choice ensures the nonincreasing character of $\epsilon\to r_\epsilon$.
Altogether, given $\epsilon\in(0,\epsilon_0]$, there exists $r_\epsilon>0$ (with $\lim_{\epsilon\to0}r_\epsilon=\infty$ and nonincreasing $\epsilon\mapsto r_\epsilon$) such that,
for $r\geq r_\epsilon$,
each one of the equations
\begin{equation}\label{eq:3Gammar}
 x'=f_{\Gamma^\pm_r}(t,x)
\end{equation}
has three uniformly separated hyperbolic
solutions $\tilde l_{\Gamma^\pm_r}$, $\tilde m_{\Gamma^\pm_r}$ and $\tilde u_{\Gamma^\pm_r}$,
satisfying $\|\tilde l_{\Gamma^\pm_r}-\tilde l_{\gamma_{\pm}}\|_\infty\le\epsilon$,
$\|\tilde m_{\Gamma^\pm_r}-\tilde m_{\gamma_{\pm}}\|_\infty\le\epsilon$ and
$\|\tilde u_{\Gamma^\pm_r}-\tilde u_{\gamma_{\pm}}\|_\infty\le\epsilon$.
It is easy to check that the pairs $(f,\Gamma_r^+)$ (resp. $(f,\Gamma_r^-)$) satisfy the hypotheses of Proposition~\ref{prop:3three}(i) (resp. (ii)), which provides fundamental information for the next steps.
\par
(i) Let us take $\epsilon\in(0,\epsilon_0]$ and the value of $r_\epsilon$ before described. We call $r=r_\epsilon$,
define $\mathfrak l^-_r$ and $\mathfrak u^-_r$ as the maximal
solutions of \eqref{eq:2completeproblem} satisfying $\mathfrak l^-_r(-r)=\tilde l_{\Gamma^-_r}(-r)$
and $\mathfrak u^-_r(-r)=\tilde u_{\Gamma^-_r}(-r)$, observe that
$\mathfrak l^-_r(t)=\tilde l_{\Gamma^-_r}(t)$
and $\mathfrak u^-_r(t)=\tilde u_{\Gamma^-_r}(t)$ for all $t\le -r$,
and deduce that $\mathfrak l^-_r$ and $\mathfrak u^-_r$ are bounded (see Proposition \ref{prop:1forwarddefinedcoercive}).
The characterizations ``$x_0\le \tilde u_{\Gamma^-_r}(-r)=\mathfrak u^-_r(-r)$
if and only if the solution $x_{\Gamma^-_r}(t,-r,x_0)$ is bounded from above as time decreases" (see Proposition~\ref{prop:3three}(ii))
and
``$x_0\le u_\Gamma(-r)$ if and only if the solution $x_\Gamma(t,-r,x_0)$ is bounded from above as time decreases",
combined with the equality $x_{\Gamma^-_r}(t,-r,x_0)=x_\Gamma(t,-r,x_0)$ for all $t\le -r$ in
the (common) interval of definition, ensure that $\mathfrak u^-_r(-r)=u_\Gamma(-r)$
and hence the global equality $\mathfrak u^-_r=u_\Gamma$. Therefore,
$|\tilde u_{\gamma_-}(t)-u_\Gamma(t)|=|\tilde u_{\gamma_-}(t)-\mathfrak u^-_r(t)|=
|\tilde u_{\gamma_-}(t)-\tilde u_{\Gamma^-_r}(t)|\le\epsilon$ for $t\le-r$. This shows that $\lim_{t\to -\infty}(u_\Gamma(t)-\tilde u_{\gamma_-}(t))=0$. Analogous arguments show that
$\mathfrak l^-_r=l_\Gamma$ and $\lim_{t\to -\infty}(l_\Gamma(t)-\tilde l_{\gamma_-}(t))=0$.

We take $x_0\in (l_\Gamma(s),u_\Gamma(s))$, note that $\bar x_0=x_\Gamma(-r,s,x_0)
\in(l_\Gamma(-r),u_\Gamma(-r))=(\tilde l_{\Gamma^-_r}(-r),\tilde u_{\Gamma^-_r}(-r))$,
and deduce from $x_\Gamma(t,s,x_0)=x_\Gamma(t,-r,\bar x_0)=x_{\Gamma^-_r}(t,-r,\bar x_0)$ for $t\le -r$ and from
Proposition~\ref{prop:3three}(ii) that $\lim_{t\to-\infty}(x_\Gamma(t,s,x_0)-\tilde m_{\gamma_-}(t))=0$.
This completes the proof of (i).
\par
(ii) Let us fix $\epsilon\in(0,\epsilon_0]$ and $r=r_\epsilon$ as above, and define
$m_\Gamma$ as the (perhaps locally defined) maximal
solution of \eqref{eq:2completeproblem} satisfying $m_\Gamma(r)=\tilde m_{\Gamma^+_r}(r)$.
We take $s$ in the domain of $m_\Gamma$ and $x_0>m_\Gamma(s)$,
and observe that, to prove that $\lim_{t\to \infty} (x_\Gamma(t,s,x_0)-\tilde u_{\gamma_+}(t))=0$,
there is no restriction in assuming that $s\ge r$. Then, $x_\Gamma(t,s,x_0)=x_{\Gamma_r^+}(t,s,x_0)$,
and the assertion follows from Proposition \ref{prop:3three}(i). The same argument
shows that $\lim_{t\to \infty} (x_\Gamma(t,s,x_0)-\tilde l_{\gamma_+}(t))=0$ if $x_0<m_\Gamma(s)$.
This duality shows that $m_\Gamma$ is unique and independent of the choice of $r$; i.e., of the choice
of $\epsilon$. In addition, $m_\Gamma(t)=\tilde m_{\Gamma^+_r}(t)$ for all $t\ge r$, which ensures that
$|m_\Gamma(t)-\tilde m_{\gamma_+}(t)|\le\epsilon$ for all $t\ge r$. The independence of
$\epsilon$ ensures that $\lim_{t\to\infty} (m_\Gamma(t)-\tilde m_{\gamma_+}(t))=0$.
The last assertion in (ii) is an immediate consequence of the previous ones.

(iii) We take $\epsilon\in(0,\epsilon_0]$ and $r=r_\epsilon$. The previously proved properties and
the repulsive hyperbolicity of $\tilde m_{\Gamma^+_r}$ provide $\delta\in(0,\epsilon)$, $k>0$ and $\beta>0$
such that $|m_\Gamma(t)-x_\Gamma(t,s,m_\Gamma(s)\pm\delta)|=
|\tilde m_{\Gamma^+_r}(t)-x_{\Gamma^+_r}(t,s,\tilde m_{\Gamma^+_r}(s)\pm\delta)|\le k\,\delta \,e^{\beta(t-s)}$ if $r\le s$ and $t\in [r,s]$,
which ensures that $m_\Gamma$ is locally
pullback repulsive. Analogous arguments prove the assertions for $l_\Gamma$ and $u_\Gamma$.
\end{proof}
The last assertion in point (ii) of the previous theorem
contributes to clarify the formulations of Theorems \ref{th:2massobrecasos} and \ref{th:2massobrecasos-2},
whose proofs use a previous
technical lemma. The (easy) proof of this lemma, which just requires \hyperlink{h1}{\textbf{h1}}-\hyperlink{h2}{\textbf{h2}},
is left to the reader, who can also find a similar one in Theorem~4.3.6 of \cite{tfmduenas}.
\begin{lemma}\label{lemm:3limit} Let $(f,\Gamma)$ satisfy {\rm\hyperlink{h1}{\textbf{h1}}} and {\rm\hyperlink{h2}{\textbf{h2}}}. Let $\bar x$ be a bounded solution of \eqref{eq:2completeproblem}.
\begin{itemize}
\item[\rm(i)] If the graph of $\bar x$ approaches that of an attractive (resp.~repulsive)
hyperbolic solution of \eqref{eq:2pastproblem} as time decreases, then there exist $s_0>0$, $k\ge 1$ and
$\beta>0$ such that $\exp\int_s^t (f_\Gamma)_x(l,\bar x(l))\,dl\le k e^{-\beta(t-s)}$ if $s\le t\le -s_0$
(resp.~$\exp\int_s^t (f_\Gamma)_x(l,\bar x(l))\,dl\le k e^{\beta(t-s)}$ if  $t\le s\le -s_0$).
\item[\rm(ii)] If the graph of $\bar x$ approaches that of an attractive (resp.~repulsive)
hyperbolic solution of \eqref{eq:2futureproblem} as time increases, then there exist $s_0>0$, $k\ge 1$ and
$\beta>0$ such that $\exp\int_s^t (f_\Gamma)_x(l,\bar x(l))\,dl\le k e^{-\beta(t-s)}$ if $s_0\le s\le t$
(resp.~$\exp\int_s^t (f_\Gamma)_x(l,\bar x(l))\,dl\le k e^{\beta(t-s)}$ if  $s_0\le t\le s$).
\end{itemize}
Consequently,
\begin{itemize}
\item[\rm(iii)] if its graph approaches that of an attractive (resp.~repulsive)
hyperbolic solution of \eqref{eq:2pastproblem} as time decreases and that of an attractive (resp.~repulsive)
hyperbolic solution of \eqref{eq:2futureproblem} as time increases, then $\bar x_\Gamma$ is an attractive (resp.~repulsive) hyperbolic
solution of \eqref{eq:2completeproblem}.
\item[\rm(iv)] If its graph approaches that of an attractive (resp.~repulsive)
hyperbolic solution of \eqref{eq:2pastproblem} as time decreases and that of a repulsive (resp.~attractive)
hyperbolic solution of \eqref{eq:2futureproblem} as time increases, then $\bar x_\Gamma$ is a nonhyperbolic
solution of \eqref{eq:2completeproblem}.
\end{itemize}
\end{lemma}
\begin{theorem}\label{th:2massobrecasos}
Let $(f,\Gamma)$ satisfy {\rm\hyperlink{h1}{\textbf{h1}}-\hyperlink{h6}{\textbf{h6}}} and
$\lim_{t\to\infty}(u_\Gamma(t)-\tilde u_{\gamma_+}(t))=0$. Then, $\tilde u_\Gamma=u_\Gamma$
is an attractive hyperbolic solution, and one of the following cases holds:
\begin{itemize}
\item[\rm (1)] $\lim_{t\to\infty}(l_\Gamma(t)-\tilde l_{\gamma_+}(t))=0$, in which case
$m_\Gamma$ is globally defined and satisfies $\lim_{t\to-\infty}(m_\Gamma(t)-\tilde m_{\gamma_-}(t))=0$;
if so, the equation \eqref{eq:2completeproblem} is in \textsc{Case \hyperlink{A}{A}}, being $\tilde l_\Gamma=l_\Gamma$, $\tilde m_\Gamma=m_\Gamma$ and
$\tilde u_\Gamma=u_\Gamma$ the three hyperbolic solutions.
\item[\rm (2)] $\lim_{t\to\infty}(l_\Gamma(t)-\tilde m_{\gamma_+}(t))=0$ or,
equivalently, $m_\Gamma=l_\Gamma$; if so, the equation \eqref{eq:2completeproblem} is in \textsc{Case \hyperlink{B}{B}},
being $\tilde u_\Gamma=u_\Gamma$ the unique hyperbolic solution;
\item[\rm (3)] $\lim_{t\to\infty}(l_\Gamma(t)-\tilde u_{\gamma_+}(t))=0$, in which case
$\lim_{t\to\infty}(x_\Gamma(t)-\tilde u_{\gamma_+}(t))=0$ for every bounded solution
$x_\Gamma(t)$ of \eqref{eq:2completeproblem}; if so, the equation \eqref{eq:2completeproblem} is in \textsc{Case \hyperlink{C}{C}},
being $\tilde u_\Gamma=u_\Gamma$ and $\tilde l_\Gamma=l_\Gamma$ the unique hyperbolic solutions.
\end{itemize}
In addition, if $m_\Gamma=l_\Gamma$ or if $\lim_{t\to\infty}(l_\Gamma(t)-\tilde u_{\gamma_+}(t))=0$,
then $\lim_{t\to\infty}(u_\Gamma(t)-\tilde u_{\gamma_+}(t))=0$.
\end{theorem}
\begin{proof}
Theorem \ref{th:3general}(ii) shows that (1), (2) and (3) exhaust the possibilities for the
limiting behavior of $l_\Gamma$, as well as the equivalence stated in (2). The hyperbolicity of
$u_\Gamma$ follows from the assumed condition $\lim_{t\to\infty}(u_\Gamma(t)-\tilde u_{\gamma_+}(t))=0$,
the fact that $\lim_{t\to-\infty}(u_\Gamma(t)-\tilde u_{\gamma_-}(t))=0$ (see Theorem \ref{th:3general}(i)),
and Lemma \ref{lemm:3limit}(iii).
\par
Let us analyze situation (1): $\lim_{t\to\infty}(l_\Gamma(t)-\tilde l_{\gamma_+}(t))=0$ and
$\lim_{t\to\infty}(u_\Gamma(t)-\tilde u_{\gamma_+}(t))=0$.
That is, the graphs of $l_\Gamma$, $m_\Gamma$ and $u_\Gamma$ respectively
approach those of $\tilde l_{\gamma_+}$, $\tilde m_{\gamma_+}$ and $\tilde u_{\gamma_+}$
as time increases (see Theorem \ref{th:3general}(ii)),
which ensures that $m_\Gamma(t)\in (l_\Gamma(t),u_\Gamma(t))$ for large enough $t$.
In particular, $m_\Gamma$ is globally defined. In these conditions, and according to
Theorem \ref{th:3general}(i), the graphs of $l_\Gamma$, $m_\Gamma$ and $u_\Gamma$
respectively approach those of $\tilde l_{\gamma_-}$, $\tilde m_{\gamma_-}$ and $\tilde u_{\gamma_-}$
as time decreases. Hence, $l_\Gamma<m_\Gamma<u_\Gamma$ are three uniformly
separated solutions, and Theorem \ref{th:3three} proves the remaining assertions.
\par
Assume that (2) holds (i.e., $l_\Gamma=m_\Gamma$: see Theorem \ref{th:3general}(ii)), and let $\bar x_\Gamma$ be a bounded solution
of \eqref{eq:2completeproblem} with $l_\Gamma<\bar x_\Gamma<u_\Gamma$. Theorem
\ref{th:3general}(i) and (ii) ensure that
$\lim_{t\to-\infty}(\bar x_\Gamma(t)-\tilde m_{\gamma_-}(t))=0$ and
$\lim_{t\to\infty}(\bar x_\Gamma(t)-\tilde u_{\gamma_+}(t))=0$.
Therefore, $l_\Gamma=m_\Gamma$ and $u_\Gamma$ are the unique uniformly separated
(bounded) solutions, with $\tilde u_\Gamma=u_\Gamma$ hyperbolic attractive and,
according to Theorem \ref{th:3general}(iii), $l_\Gamma$
locally pullback attractive and repulsive; and
Lemma \ref{lemm:3limit}(iv) ensures that $\bar x_\Gamma$ and $l_\Gamma$
are nonhyperbolic. We conclude that the dynamics fits in \textsc{Case \hyperlink{B}{B}}.
\par
Finally, we assume that (3) holds; i.e., that $\lim_{t\to\infty}(l_\Gamma(t)-\tilde u_{\gamma_+}(t))=0$.
In this case, any bounded solution $\bar x_\Gamma$ of \eqref{eq:2completeproblem}
with $l_\Gamma\le \bar x_\Gamma<u_\Gamma$ satisfies
$\lim_{t\to\infty}(\bar x_\Gamma(t)-\tilde u_{\gamma_+}(t))=0$, and if
$l_\Gamma<\bar x_\Gamma<u_\Gamma$, then $\lim_{t\to-\infty}(\bar x_\Gamma(t)-\tilde m_{\gamma_-}(t))=0$.
This precludes the existence of (bounded) uniformly separated solutions, as well
as the existence of hyperbolic solutions different from $u_\Gamma$ and $l_\Gamma$.
In addition, Lemma \ref{lemm:3limit}(iii) shows
the attractive hyperbolicity of $l_\Gamma$.
Altogether, we have check that the dynamics fits in \textsc{Case \hyperlink{C}{C}}.
\par
The last assertions of the theorem follow easily from Theorem \ref{th:3general}(ii).
\end{proof}

Analogous arguments lead to the proof of the next ``symmetric"~result:
\begin{theorem}\label{th:2massobrecasos-2}
Let $(f,\Gamma)$ satisfy {\rm\hyperlink{h1}{\textbf{h1}}-\hyperlink{h6}{\textbf{h6}}} and
$\lim_{t\to\infty}(l_\Gamma(t)-\tilde l_{\gamma_+}(t))=0$. Then, $\tilde l_\Gamma=l_\Gamma$
is an attractive hyperbolic solution, and one of the following cases holds:
\begin{itemize}
\item[\rm (1)] $\lim_{t\to\infty}(u_\Gamma(t)-\tilde u_{\gamma_+}(t))=0$, in which case
$m_\Gamma$ is globally defined and satisfies $\lim_{t\to-\infty}(m_\Gamma(t)-\tilde m_{\gamma_-}(t))=0$;
if so, the equation \eqref{eq:2completeproblem} is in \textsc{Case \hyperlink{A}{A}},
being $\tilde l_\Gamma=l_\Gamma$, $\tilde m_\Gamma=m_\Gamma$ and
$\tilde u_\Gamma=u_\Gamma$ the three hyperbolic solutions.
\item[\rm (2)] $\lim_{t\to\infty}(u_\Gamma(t)-\tilde m_{\gamma_+}(t))=0$ or,
equivalently, $m_\Gamma=u_\Gamma$; if so, the equation \eqref{eq:2completeproblem} is in \textsc{Case \hyperlink{B}{B}},
being $\tilde l_\Gamma=l_\Gamma$ the unique hyperbolic solution;
\item[\rm (3)] $\lim_{t\to\infty}(u_\Gamma(t)-\tilde l_{\gamma_+}(t))=0$, in which case
$\lim_{t\to\infty}(x_\Gamma(t)-\tilde l_{\gamma_+}(t))$ for every bounded solution
$x_\Gamma(t)$ of \eqref{eq:2completeproblem}; if so, the equation \eqref{eq:2completeproblem} is in \textsc{Case \hyperlink{C}{C}},
being $\tilde u_\Gamma=u_\Gamma$ and $\tilde l_\Gamma=l_\Gamma$ the unique hyperbolic solutions.
\end{itemize}
In addition, if $m_\Gamma=u_\Gamma$ or if $\lim_{t\to\infty}(u_\Gamma(t)-\tilde l_{\gamma_+}(t))=0$,
then $\lim_{t\to\infty}(l_\Gamma(t)-\tilde l_{\gamma_+}(t))=0$.
\end{theorem}
Reading carefully the statements of the last three theorems leads us to the following statement.
\begin{corollary}\label{coro:3exhaust}
Let $(f,\Gamma)$ satisfy {\rm\hyperlink{h1}{\textbf{h1}}-\hyperlink{h6}{\textbf{h6}}}. The situations described in Theorems~{\rm\ref{th:2massobrecasos}} and {\rm\ref{th:2massobrecasos-2}} exhaust the dynamical possibilities for \eqref{eq:2completeproblem}.
Consequently,
\textsc{Cases \hyperlink{A}{A}}, \hyperlink{B}{\textsc{B}} and \hyperlink{C}{\textsc{C}} exhaust the dynamical possibilities for \eqref{eq:2completeproblem}.

Furthermore, if we also assume $\lim_{t\to\infty}(u_\Gamma(t)-\tilde u_{\gamma_+}(t))=0$ (resp. $\lim_{t\to\infty}(l_\Gamma(t)-\tilde l_{\gamma_+}(t))=0$), then
\textsc{Case \hyperlink{A}{A}} holds if and only if there exists
$t_\Gamma\in\mathbb R$ such that $l_\Gamma(t_\Gamma)<m_\Gamma(t_\Gamma)$ (resp. $m_\Gamma(t_\Gamma)<u_\Gamma(t_\Gamma)$),
\textsc{Case \hyperlink{B}{B}} holds if and only if there exists
$t_\Gamma\in\mathbb R$ such that $l_\Gamma(t_\Gamma)=m_\Gamma(t_\Gamma)$ (resp. $m_\Gamma(t_\Gamma)=u_\Gamma(t_\Gamma)$),
and \textsc{Case \hyperlink{C}{C}} holds if and only if there exists
$t_\Gamma\in\mathbb R$ such that $l_\Gamma(t_\Gamma)>m_\Gamma(t_\Gamma)$ (resp. $m_\Gamma(t_\Gamma)>u_\Gamma(t_\Gamma)$).
\end{corollary}
Hereafter, \textsc{Cases \hyperlink{B}{B}} and \hyperlink{C}{\textsc{C}} will be called \hypertarget{B1}{\textsc{Cases B1}} and \hypertarget{C1}{\textsc{C1}} under the hypotheses of cases (2) and (3) of Theorem~\ref{th:2massobrecasos}, and \hypertarget{B2}{\textsc{Cases B2}} and \hypertarget{C2}{\textsc{C2}} under the hypotheses of cases (2) and (3) Theorem~\ref{th:2massobrecasos-2}.

The results seen so far show that, if $(f,\Gamma)$ satisfies {\rm\hyperlink{h1}{\textbf{h1}}-\hyperlink{h6}{\textbf{h6}}}, then the upper and lower
bounded solutions of \eqref{eq:2completeproblem} always ``connect"~with those of the past:
$\lim_{t\to-\infty}(l_\Gamma(t)-\tilde l_{\gamma_-}(t))=
\lim_{t\to-\infty}(u_\Gamma(t)-\tilde u_{\gamma_-}(t))=0$.
In addition,
\par
\begin{figure}[h]
\centering
\includegraphics[width=0.48\textwidth]{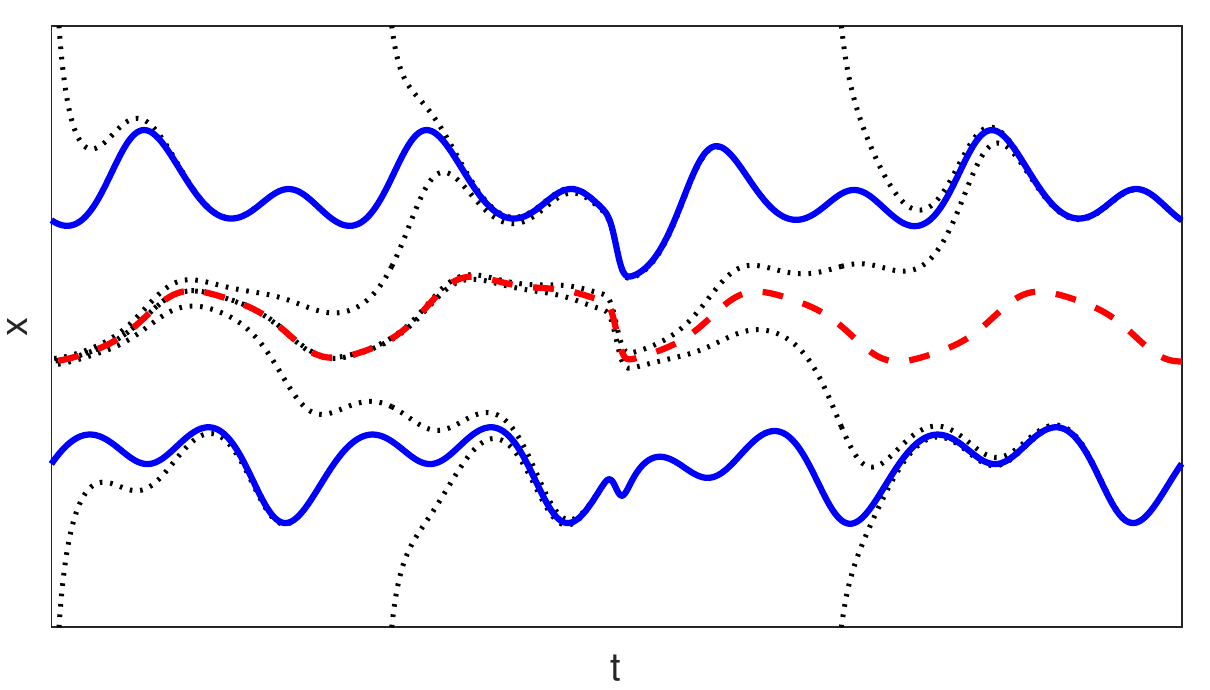}\par
\caption{Example of \textsc{Case \protect\hyperlink{A}{A}}. In blue solid line the attractive hyperbolic solutions $l_\Gamma$ and $u_\Gamma$, and in red dashed line the repulsive hyperbolic solution $m_\Gamma$. In black dotted lines, as in Figures \ref{fig:CaseB} and \ref{fig:CaseC}, typical solutions different from the significant ones.}\label{fig:CaseA}
\end{figure}

- \textsc{Case \hyperlink{A}{A}} holds when the transition equation has three hyperbolic solutions and they
connect those of the past equation to those of the future equation. More precisely, when
$\lim_{t\to \infty}(l_\Gamma(t)-\tilde l_{\gamma_+}(t))=
\lim_{t\to \infty}(u_\Gamma(t)-\tilde u_{\gamma_+}(t))=0$.
In this case, $m_\Gamma$ exists globally, and $\lim_{t\to-\infty}(m_\Gamma(t)-\tilde m_{\gamma_-}(t))=\lim_{t\to\infty}(m_\Gamma(t)-\tilde m_{\gamma_+}(t))=0$. The dynamics of the transition equation is described in Theorem \ref{th:3three} and
depicted in Figure \ref{fig:CaseA}. Note that \textsc{Case \hyperlink{A}{A}} is the unique one in which a
repulsive hyperbolic solution exists.
In addition, Proposition \ref{prop:3forward} below shows that \textsc{Case \hyperlink{A}{A}}
holds if and only if its pullback attractor has the property of forward attraction.
In addition, \textsc{Case \hyperlink{A}{A}} is a {\em robust\/} property: Proposition~\ref{prop:3CaseCrobust} below proves that it is
persistent under $\|{\cdot}\|_\infty$-small perturbations.
\par
\begin{figure}[h]
\centering
\begin{subfigure}[b]{0.48\textwidth}
         \centering
         \includegraphics[width=\textwidth]{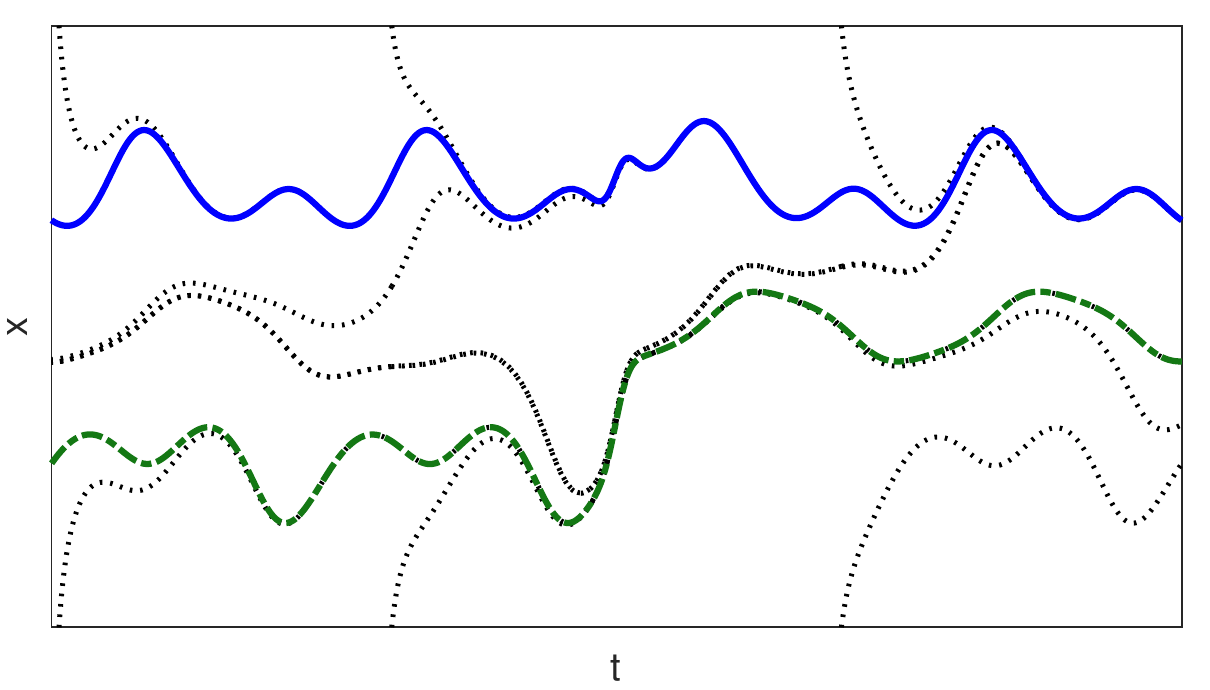}
         \caption{\textsc{Case \protect\hyperlink{B1}{B1}}. In blue solid line the attractive hyperbolic solution $u_\Gamma$, in green dashed-dotted line $l_\Gamma=m_\Gamma$.}
         \label{fig:CASOB1}
     \end{subfigure}
\hfill
\begin{subfigure}[b]{0.48\textwidth}
         \centering
         \includegraphics[width=\textwidth]{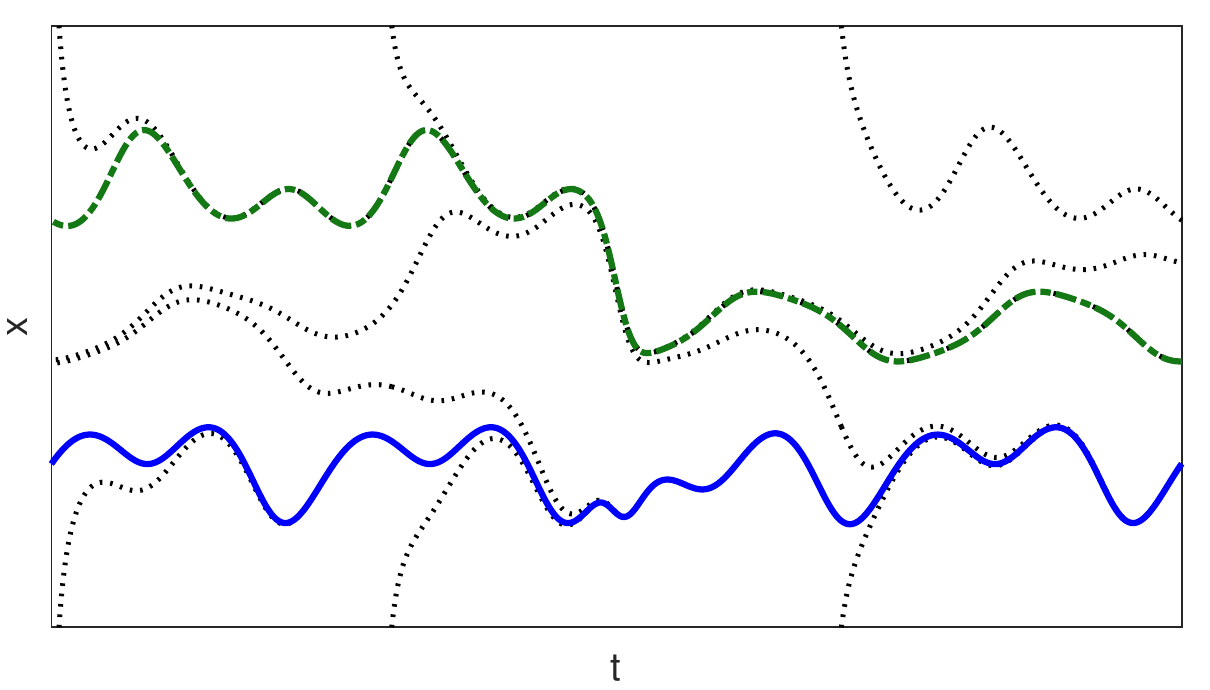}
         \caption{\textsc{Case \protect\hyperlink{B2}{B2}}. In blue solid line the attractive hyperbolic solution $l_\Gamma$, in green dashed-dotted line $u_\Gamma=m_\Gamma$.}
         \label{fig:CASOB2}
     \end{subfigure}
\caption{Examples of \textsc{Case \protect\hyperlink{B}{B}}.}\label{fig:CaseB}
\end{figure}

- \textsc{Case \hyperlink{B}{B}} holds when the transition equation has exactly one attractive hyperbolic
solution, given by the upper or lower bounded solution, which connects the
corresponding upper or lower bounded solution of the past to that of the future.
In this case, the other limiting bounded solution connects that of the past with the
repulsive hyperbolic solution of the future.
More precisely, \textsc{Case \hyperlink{B1}{B1}} holds when $\lim_{t\to \infty}(l_\Gamma(t)-\tilde m_{\gamma_+}(t))=
\lim_{t\to \infty}(u_\Gamma(t)-\tilde u_{\gamma_+}(t))=0$ and \textsc{Case \hyperlink{B2}{B2}} holds when $\lim_{t\to \infty}(l_\Gamma(t)-\tilde l_{\gamma_+}(t))=
\lim_{t\to \infty}(u_\Gamma(t)-\tilde m_{\gamma_+}(t))=0$. The dynamics of
the transition equation in \textsc{Cases \hyperlink{B1}{B1}} and \hyperlink{B2}{\textsc{B2}}, described in Theorem
\ref{th:3general}, is depicted in Figures \ref{fig:CASOB1} and \ref{fig:CASOB2}, respectively.
There, we observe that the functions $\tilde l_{\gamma_+}$ and
$\tilde u_{\gamma_+}$ also play a fundamental role in \textsc{Cases \hyperlink{B1}{B1}} and \hyperlink{B2}{\textsc{B2}}, respectively:
they determine the limiting behaviour as time increases of the solutions which, as time
decreases, are unbounded from below (in \textsc{Case \hyperlink{B1}{B1}}) or from above (in \textsc{Case \hyperlink{B2}{B2}}).
Numerical examples show that this type of dynamics is, in general, highly unstable, far away
from persisting under small perturbations of $\Gamma$. This lack of robustness fits with the
intuitive lack of robustness of the conditions $l_\Gamma(t_\Gamma)=m_\Gamma(t_\Gamma)$ or $u_\Gamma(t_\Gamma)=m_\Gamma(t_\Gamma)$ of
Corollary \ref{coro:3exhaust}.
\par
\begin{figure}[h]
\centering
\centering
\begin{subfigure}[b]{0.48\textwidth}
         \centering
         \includegraphics[width=\textwidth]{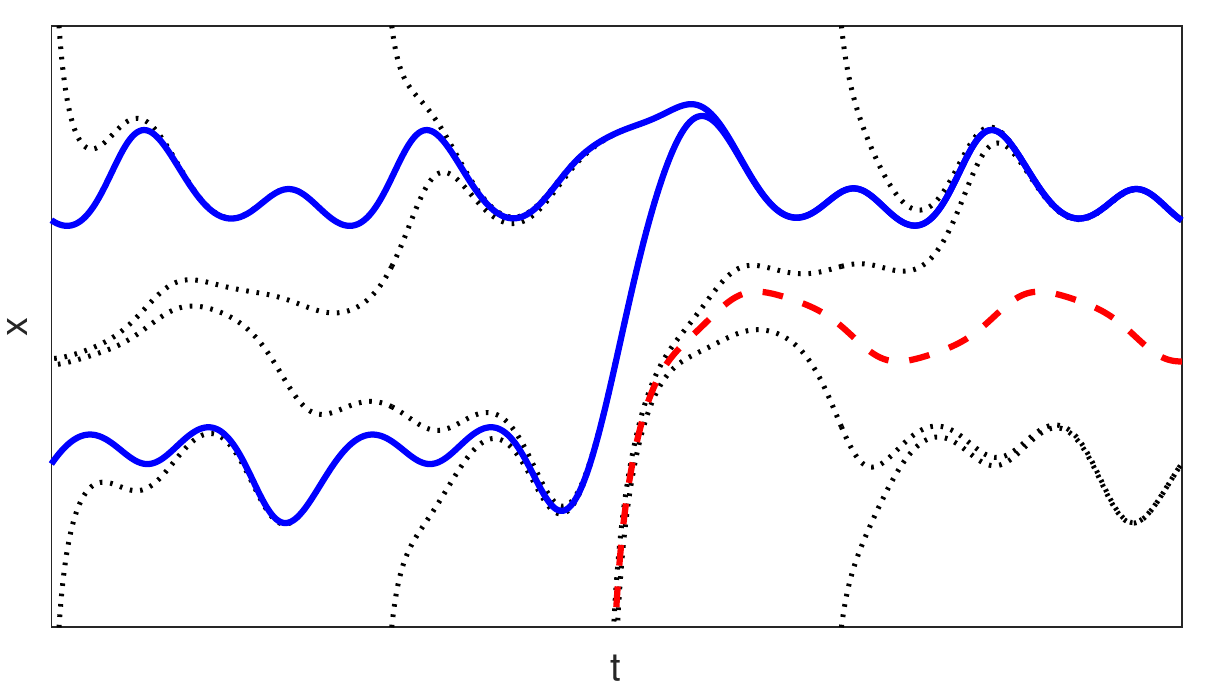}
         \caption{\textsc{Case \protect\hyperlink{C1}{C1}}. In blue solid line the attractive hyperbolic solutions $l_\Gamma$ and $u_\Gamma$, in red dashed line $m_\Gamma$.}
         \label{fig:CASOC1}
     \end{subfigure}
\hfill
\begin{subfigure}[b]{0.48\textwidth}
         \centering
         \includegraphics[width=\textwidth]{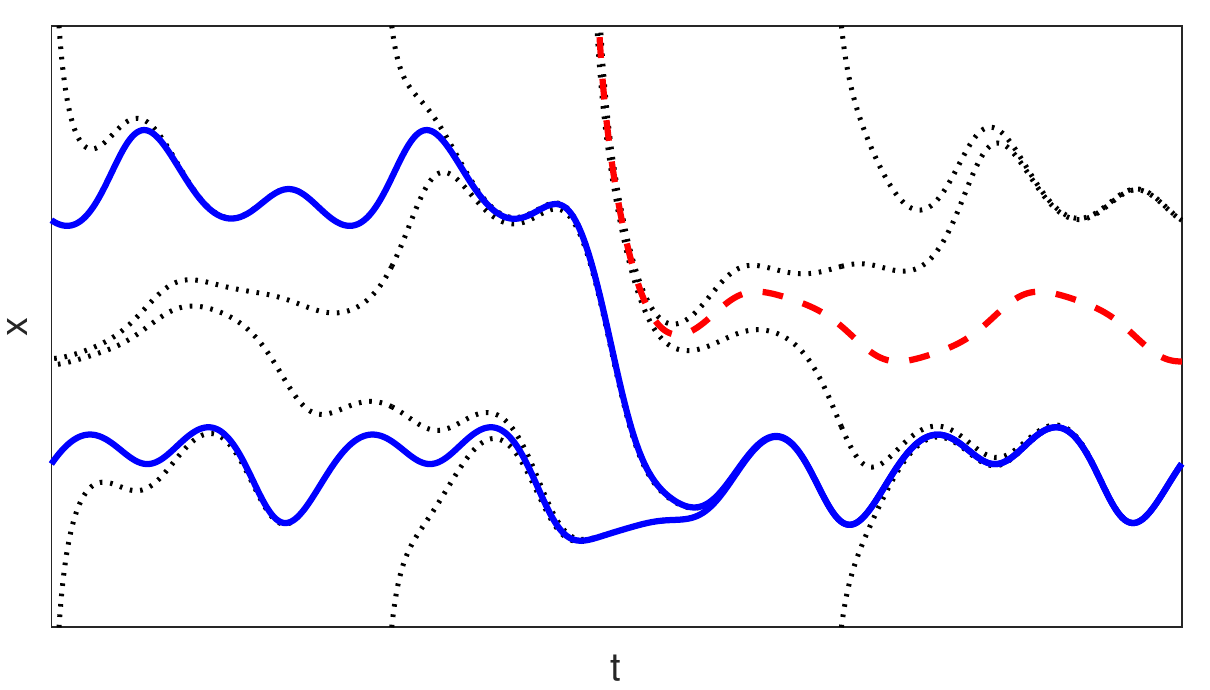}
         \caption{\textsc{Case \protect\hyperlink{C2}{C2}}. In blue solid line the attractive hyperbolic solutions $l_\Gamma$ and $u_\Gamma$, in red dashed line $m_\Gamma$.}
         \label{fig:CASOC2}
     \end{subfigure}
\caption{Examples of \textsc{Case \protect\hyperlink{C}{C}}.}\label{fig:CaseC}
\end{figure}

- \textsc{Case \hyperlink{C}{C}} holds when the transition equation has exactly two hyperbolic
solutions, which are attractive and given by the (non uniformly separated)
upper and lower bounded solutions, and which connect both the
upper and lower bounded solutions of the past to either the upper one
or the lower one of the future.
More precisely, \textsc{Case \hyperlink{C1}{C1}} holds when
$\lim_{t\to \infty}(l_\Gamma(t)-\tilde u_{\gamma_+}(t))=\lim_{t\to \infty}(u_\Gamma(t)-\tilde u_{\gamma_+}(t))=0$
and \textsc{Case \hyperlink{C2}{C2}} holds when $\lim_{t\to \infty}(l_\Gamma(t)-\tilde l_{\gamma_+}(t))=
\lim_{t\to \infty}(u_\Gamma(t)-\tilde l_{\gamma_+}(t))=0$. The dynamics of
the transition equation in \textsc{Cases \hyperlink{C1}{C1}} and \hyperlink{C2}{\textsc{C2}}, described in Theorem
\ref{th:3general}, are depicted in Figures \ref{fig:CASOC1} and \ref{fig:CASOC2},
where we observe the role played by all the hyperbolic solutions of
the future equation: they determine the limiting behaviour as time increases of the solutions
which, as time decreases, are unbounded from below (in \textsc{Case \hyperlink{C1}{C1}}) or from above (in \textsc{Case \hyperlink{C2}{C2}}).
Proposition \ref{prop:3CaseCrobust} below also shows the robustness of these dynamical situations under small perturbations of $\Gamma$.
\par
There is a fundamental difference between \textsc{Case \hyperlink{A}{A}} and \textsc{Cases \hyperlink{B}{B}} or \hyperlink{C}{\textsc{C}}, which can be formulated
in terms of the limiting behavior as time increases of the set of bounded solutions (closely
related to the pullback attractor, as explained in Subsection \ref{subsec:1pullback}):
\textsc{Case \hyperlink{A}{A}} holds if and only if the sections of this set ``connect" those of the past with those of the future,
in the sense before explained. In \textsc{Cases \hyperlink{B}{B}} and \hyperlink{C}{\textsc{C}}, the set of bounded solutions
of the transition equation connects that of the past
with just a part of that of the future. For this reason,
by analogy with similar situations (see \cite{aspw}, \cite{lno1,lno2}),
we refer to this \textsc{Case \hyperlink{A}{A}} as {\em tracking}, and talk about {\em tipping\/}
in \textsc{Cases \hyperlink{B}{B}} and \hyperlink{C}{\textsc{C}}.
\textsc{Cases \hyperlink{B}{B}} and \hyperlink{C}{\textsc{C}} can be distinguished by the boundedness or unboundedness of $m_\Gamma$, but
one of the most important differences between \textsc{Cases \hyperlink{B}{B}} and \hyperlink{C}{\textsc{C}} is the
robustness of \textsc{Case \hyperlink{C}{C}}.

We complete this subsection by proving the two results mentioned in the previous description. As said in Subsection \ref{subsec:1pullback}, in general, pullback and forward attraction are unrelated (see e.g.~\cite{laos}). The next result studies the local or global forward attraction properties of the pullback attractor $\mathscr{A}_\Gamma=\{\mathscr{A}_\Gamma(t)\colon\, t\in\mathbb R\}$ of \eqref{eq:2completeproblem}, whose existence is guaranteed by Proposition \ref{prop:1forwarddefinedcoercive} under our hypotheses on $(f,\Gamma)$.
\begin{proposition}\label{prop:3forward}
Let $(f,\Gamma)$ satisfy {\rm\hyperlink{h1}{\textbf{h1}}-\hyperlink{h6}{\textbf{h6}}}.
Then, the set $\mathcal B_\Gamma$ given by \eqref{eq:3conjuntosolucionesacotadas} is uniformly exponentially stable
if and only if \eqref{eq:2completeproblem} is in \textsc{Cases \hyperlink{A}{A}} or \textsc{\hyperlink{C}{C}}. In addition,
\begin{enumerate}[label=\rm{(\roman*)}]
\item $\mathscr{A}_\Gamma$ is globally forward attractive if and only if \eqref{eq:2completeproblem} is in \textsc{Case \hyperlink{A}{A}}.
\item $\mathscr{A}_\Gamma$ is not locally forward attractive if and only if
    \eqref{eq:2completeproblem} is in \textsc{Case \hyperlink{B}{B}}.
\item $\mathscr{A}_\Gamma$ is locally but not globally forward attractive if and only if \eqref{eq:2completeproblem} is in \textsc{Case \hyperlink{C}{C}}.
\end{enumerate}
\end{proposition}
\begin{proof}
Recall that $\mathscr{A}_\Gamma(t)=[l_\Gamma(t),u_\Gamma(t)]$: see Remark \ref{remark:sectionsboundedsolutionsopullbackattractor}.
In \textsc{Cases \hyperlink{A}{A}} and \textsc{\hyperlink{C}{C}}, the attractive hyperbolic properties of $u_\Gamma$ and $l_\Gamma$ established in Theorems~\ref{th:2massobrecasos} and \ref{th:2massobrecasos-2} ensure that $\mathcal B_\Gamma$ is uniformly exponentially stable, and it is easy to deduce that $\mathscr A_\Gamma$ is locally forward attractive. Hence, to prove (ii), it suffices to check that $\mathscr A_\Gamma$ is not locally forward attractive in \textsc{Cases \hyperlink{B}{B}}. In \textsc{Case \hyperlink{B1}{B1}}, we take $s\in\mathbb R$ for which there exists $m_\Gamma(s)$, and any $x_0<l_\Gamma(s)=m_\Gamma(s)$. Theorem~\ref{th:3general}(ii) ensures that $\lim_{t\to\infty}(x_\Gamma(t,s,x_0)-\tilde l_{\gamma_+}(t))=0$, which combined with $\lim_{t\to\infty}(l_\Gamma(t)-\tilde m_{\gamma_+}(t))=0$ (guaranteed by Theorem \ref{th:2massobrecasos}(2)) and with the uniform separation of $\tilde l_{\gamma_+}$ and $\tilde m_{\gamma_+}$ precludes $\lim_{t\to\infty}(x_\Gamma(t,s,x_0)-l_\Gamma(t))=0$. In \textsc{Case \hyperlink{B2}{B2}}, we take any $x_0>u_\Gamma(s)=m_\Gamma(s)$ to preclude $\lim_{t\to\infty}(x_\Gamma(t,s,x_0)-u_\Gamma(t))=0$.
Any of these two situations contradicts the local forward attraction of $\mathscr A_\Gamma$.

Observe that $\mathscr A_\Gamma$ is globally forward attractive if and only if
$\lim_{t\rightarrow\infty}(x_\Gamma(t,s,x_0)-l_\Gamma(t))=0$ whenever $x_0<l_\Gamma(s)$ and
$\lim_{t\rightarrow\infty}(x_\Gamma(t,s,x_0)-u_\Gamma(t))=0$ whenever $x_0>u_\Gamma(s)$:
the monotonicity of the solutions with respect to initial data guarantees the uniformity on bounded sets. According to Theorem \ref{th:3three}, these two properties hold in
\textsc{Case \hyperlink{A}{A}}. Let us check that they cannot hold simultaneously in \textsc{Case \hyperlink{C}{C}},
which will complete the proofs of (i) and (iii). In \textsc{Case \hyperlink{C1}{C1}}, we take $s\in\mathbb R$ for which there exists $m_\Gamma(s)$ and $x_0<m_\Gamma(s)<l_\Gamma(s)$. Theorem~\ref{th:3general}(ii) ensures that $\lim_{t\to\infty}(x_\Gamma(t,s,x_0)-\tilde l_{\gamma_+}(t))=0$, which combined with $\lim_{t\to\infty}(l_\Gamma(t)-\tilde u_{\gamma_+}(t))=0$ (guaranteed by Theorem \ref{th:2massobrecasos}(3)) and with the uniform separation of $\tilde l_{\gamma_+}$ and $\tilde u_{\gamma_+}$ precludes
$\lim_{t\rightarrow\infty}(x_\Gamma(t,s,x_0)-l_\Gamma(t))=0$. In \textsc{Case \hyperlink{C2}{C2}},
taking $x_0>m_\Gamma(s)>u_\Gamma(s)$ precludes $\lim_{t\rightarrow\infty}(x_\Gamma(t,s,x_0)-u_\Gamma(t))=0$. The proof is complete.
\end{proof}
In order to analyze the robustness of \textsc{Cases \hyperlink{A}{A}} and \hyperlink{C}{C} under small variations of $\Gamma$ for a fixed $f$, we need to describe the set $\mathfrak{C}_f$ of functions on which $\Gamma$ may vary to get pairs $(f,\Gamma)$ satisfying \hyperlink{h1}{\textbf{h1}}-\hyperlink{h6}{\textbf{h6}}.
Given $f\colon\mathbb{R}\times\mathbb{R}\times\mathbb{R}\rightarrow\mathbb{R}$ satisfying \hyperlink{h2}{\textbf{h2}}-\hyperlink{h4}{\textbf{h4}}, we define
\begin{align}
\mathcal{R}_f&=\{\gamma\in\mathbb{R}\colon\,\text{$f$ satisfies \hyperlink{h5}{\textbf{h5}$_\gamma$} and \hyperlink{h6}{\textbf{h6}$_\gamma$}}\}\,,\label{def.R}\\
\mathfrak{C}_f&=\{
\Gamma\colon\mathbb{R}\rightarrow\mathbb{R}\mid\text{$\Gamma$ is continuous and there exist
$\lim_{t\to\pm\infty}\Gamma(t)\in \mathcal{R}_f$}\}\,.\label{def.C}
\end{align}
\begin{proposition}\label{prop:3CaseCrobust}
Let $f$ satisfy {\rm\hyperlink{h2}{\textbf{h2}}-\hyperlink{h4}{\textbf{h4}}}.
The dynamics of \textsc{Cases \hyperlink{A}{A}}, \hyperlink{C1}{\textsc{C1}} and \hyperlink{C2}{\textsc{C2}}
persist under small perturbations of $\Gamma\in\mathfrak{C}_f$ with respect to $\|{\cdot}\|_\infty$.
\end{proposition}
\begin{proof}
In \textsc{Case \hyperlink{A}{A}}, the existence of three uniformly separated hyperbolic solutions and Theorem~\ref{th:hypcontinuation} prove the result. Now, assume that the dynamics fits in \textsc{Case \hyperlink{C1}{C1}} (described in Theorem \ref{th:2massobrecasos}(3)).
As at the beginning of the proof of Theorem \ref{th:3general}, we take
$\epsilon>0$ with $\epsilon<\epsilon_0=(1/3)\min\{\inf_{t\in\mathbb R}(\tilde u_{\gamma_\pm}-\tilde m_{\gamma_\pm}),
\inf_{t\in\mathbb R}(\tilde m_{\gamma_\pm}-\tilde l_{\gamma_\pm})\}$.
Theorem \ref{th:hypcontinuation} provides $\delta>0$ such that, if $\Delta\in\mathfrak{C}_f$ and $\|\Delta-\Gamma\|_{\infty}\le\delta$, then:
\begin{itemize}
\item[-] the equation $x'=f_\Delta(t,x)$ has two attractive hyperbolic solutions (so, it is
in \textsc{Cases \hyperlink{A}{A}} or \hyperlink{C}{\textsc{C}}) $\tilde l_\Delta<\tilde u_\Delta$, with $\|\tilde l_\Delta-\tilde l_\Gamma\|_\infty\le\epsilon/3$
and $\|\tilde u_\Delta-\tilde u_\Gamma\|_\infty\le\epsilon/3$;
\item[-] if $\delta_+=\lim_{t\to\infty}\Delta(t)$, then the three hyperbolic solutions $\tilde l_{\delta_+}<\tilde m_{\delta_+}<\tilde u_{\delta_+}$ of the future equation $x'=f(t,x,\delta_+)$ satisfy $\|\tilde l_{\gamma_+}-\tilde l_{\delta_+}\|_\infty<\epsilon/3$,
$\|\tilde m_{\gamma_+}-\tilde m_{\delta_+}\|_\infty<\epsilon/3$ and
$\|\tilde u_{\gamma_+}-\tilde u_{\delta_+}\|_\infty<\epsilon/3$.
\end{itemize}
Therefore, $|\tilde l_\Delta(t)-\tilde u_{\delta_+}(t)|\le \|\tilde l_\Delta-\tilde l_\Gamma\|_\infty+
|\tilde l_\Gamma(t)-\tilde u_{\gamma_+}(t)|+\|\tilde u_{\gamma_+}-\tilde u_{\delta_+}\|_\infty\le\epsilon$
if $t$ is large enough to ensure $|\tilde l_\Gamma(t)-\tilde u_{\gamma_+}(t)|\le\epsilon/3$. This precludes
the graph of $\tilde l_\Delta$ to approach those of $\tilde l_{\delta_+}$ or $\tilde m_{\delta_+}$, and
hence the dynamics fits in \textsc{Case \hyperlink{C1}{C1}}. The proof is analogous for \textsc{Case \hyperlink{C2}{C2}}.
\end{proof}


\subsection{Skewproduct formulation of the problem}\label{subsec:2skewproduct} In this subsection, we introduce the skewproduct formalism to deal with equation \eqref{eq:2completeproblem}.
The main purpose is to apply the results of \cite{dno1} in order to bound the number of uniformly separated solutions of the equations \eqref{eq:2completeproblem}, \eqref{eq:2pastproblem} and \eqref{eq:2futureproblem}, and to guarantee its hyperbolicity when there are three of them; that is, to prove Theorem~\ref{th:3three}.

In what follows we will manage the standard concepts of local and global flows, their corresponding orbits, invariant and ergodic measures for a given flow, and $\boldsymbol\upalpha$- and $\boldsymbol\upomega$-limit sets.
The required definitions and notation can be found in \cite{dno1}.

As previously, we call $f_\Gamma(t,x)=f(t,x,\Gamma(t))$ and $f_{\gamma_{\pm}}(t,x)=f(t,x,\gamma_{\pm})$.
We represent $h{\cdot}t(s,x)=h(t+s,x)$. Under hypothesis \hyperlink{h1}{\textbf{h1}} and \hyperlink{h2}{\textbf{h2}}, Theorem~I.3.1 of \cite{shenyi} and
Theorem IV.3 of \cite{selltopdyn}
ensure that the hull $\Omega_\Gamma$ of $f_\Gamma$, given by the closure of the set $\mathcal F_\Gamma=\{(f_\Gamma){\cdot}t\colon\, t\in\mathbb{R}\}$ on the set $C(\mathbb{R}\times\mathbb{R},\mathbb{R})$ provided with the compact-open topology, is a compact metric space contained in $C^{0,2}$, as well as the continuity of the time-translation flow $\sigma\colon\mathbb R\times\Omega_\Gamma\to\Omega_\Gamma,(t,\omega)\mapsto\omega{\cdot}t$ and the existence and continuity of the maps $\mathfrak f_\Gamma, (\mathfrak f_\Gamma)_x,(\mathfrak f_\Gamma)_{xx}\colon\Omega_\Gamma\times\mathbb R\to\mathbb R$, where $\mathfrak f_\Gamma(\omega,x)=\omega(0,x)$.

Let us consider the family of equations
\begin{equation}\label{eq:3eqhull}
 x'=\mathfrak f_\Gamma(\omega{\cdot}t,x)
\end{equation}
for $\omega\in\Omega_\Gamma$, let $\mathcal I_{\omega,x_0}\to\mathbb R,\,t\mapsto \mathfrak x_\Gamma(t,\omega,x_0)$ be the
maximal solution of \eqref{eq:3eqhull} satisfying $\mathfrak x_\Gamma(0,\omega,x_0)=x_0$, and let $\mathcal U_\Gamma=
\bigcup_{(\omega,x_0)\in \Omega_\Gamma\times \mathbb R}(\mathcal I_{\omega,x_0}\times\{(\omega,x_0)\})
\subseteq\mathbb R\times\Omega_\Gamma\times\mathbb R$. The map
\[
 \tau_\Gamma\colon\mathcal U_\Gamma\subseteq\mathbb R\times\Omega_\Gamma\times\mathbb R\to\Omega_\Gamma\times\mathbb R\,,\;
 (t,\omega,x_0)\mapsto (\omega{\cdot}t,\mathfrak x_\Gamma(t,\omega,x_0))
\]
defines a (possibly local) continuous flow on the \emph{bundle} $\Omega_\Gamma\times\mathbb R$ whose first component is given by the previously defined flow $\sigma$ on the \emph{base} $\Omega_\Gamma$: a {\em skewproduct} flow. Note that the equation \eqref{eq:2completeproblem} is included in this family: it is given by the element $f_\Gamma\in\Omega_\Gamma$, and the solutions are related by $x_\Gamma(t,s,x_0)=\mathfrak x_\Gamma(t-s,f_\Gamma{\cdot}s,x_0)$. This way of associating a flow to a nonautonomous equation is the {\em skewproduct formalism} which we have mentioned several times.

Our tools to prove Theorem \ref{th:3three} are the results of \cite{dno1}, which require two conditions on the function $\mathfrak f_\Gamma$ apart from those already mentioned. The first one is a coercivity property ensuring the dissipativity of the flow $\tau_\Gamma$ and hence the existence of global attractor; that is, of a compact $\tau_\Gamma$-invariant set $\mathcal A_{\Gamma}\subset\Omega_\Gamma\times\mathbb R$ such that $\lim_{t\rightarrow\infty} \text{dist}(C{\cdot}t,\mathcal{A})=0$ for every bounded set $C$, where $C{\cdot}t=\{(\omega{\cdot}t,\mathfrak x_\Gamma(t,\omega,x_0))\colon\, (\omega,x_0)\in C\}$ and
\begin{equation*}
\text{dist}(C_1,C_2)=\sup_{(\omega_1,x_1)\in C_1}\left(\inf_{(\omega_2,x_2)\in C_2}
\big(\mathrm{dist}_{\Omega_\Gamma\times\mathbb{R}}((\omega_1,x_1),(\omega_2,x_2))\big)\right)\,.
\end{equation*}
\begin{proposition}\label{prop:3coercivity}
Let $(f,\Gamma)$ satisfy {\rm\hyperlink{h1}{\textbf{h1}}}-{\rm\hyperlink{h3}{\textbf{h3}}}.
Then, $\mathfrak f_\Gamma$ satisfies the coercivity property
$\lim_{|x|\to\infty}\mathfrak f_\Gamma(\omega,x)/x=-\infty$ uniformly in $\omega\in\Omega_\Gamma$.
Hence, all the solutions of all the equations of the family \eqref{eq:3eqhull} are globally forward defined and bounded, and there exists
a global attractor $\mathcal A_\Gamma=\bigcup_{\omega\in\Omega_\Gamma}(\{\omega\}\times[\mathfrak l_\Gamma(\omega),\mathfrak u_\Gamma(\omega)])$
composed by the union of all the bounded $\tau_\Gamma$-orbits. In addition, the family
$\mathscr{A}_\Gamma=\{\mathscr{A}_\Gamma(t)\colon\, t\in\mathbb{R}\}$ given by
$\mathscr{A}_\Gamma(t)=[\mathfrak l_\Gamma(f_\Gamma{\cdot}t),\mathfrak u_\Gamma(f_\Gamma{\cdot}t)]$
is the pullback attractor of \eqref{eq:2completeproblem}.
\end{proposition}
\begin{proof}
Given $k>0$ there exists $\rho_k>0$ such that $f_\Gamma(t,x)/x\le -k$ if
$|x|\ge \rho$ and $t\in\mathbb R$. Given any $\omega\in\Omega_\Gamma$, we take $\{t_n\}$ with
$\omega(0,x)=\lim_{n\to\infty}f_\Gamma(t_n,x)$ and deduce that
$\mathfrak f_\Gamma(\w,x)/x=\omega(0,x)/x\le-k$ if $|x|\ge\rho$,
which proves the first assertion. For the existence, shape and composition of global
attractor, see e.g.~Theorem~5.1 of \cite{dno1} and references therein.
In particular,
\[
\begin{split}
 [\mathfrak l_\Gamma(f_\Gamma{\cdot}s),\mathfrak u_\Gamma(f_\Gamma{\cdot}s)]
 &=\{x_0\colon t\mapsto \mathfrak x_\Gamma(t,f_\Gamma{\cdot}s,x_0) \text{ is bounded}\}\\
 &=\{x_0\colon t\mapsto x_\Gamma(t+s,s,x_0) \text{ is bounded}\}=[l_\Gamma(s),u_\Gamma(s)]\,,
\end{split}
\]
which proves the last assertion (see Remark \ref{remark:sectionsboundedsolutionsopullbackattractor}).
\end{proof}
The second condition is a property on strict concavity of the derivative, which is introduced in terms of ergodic measures.
More precisely, we need the derivative
$x\mapsto(\mathfrak{f}_\Gamma)_x(\omega,\cdot)$ to be concave for all $\omega\in\Omega$, as well as
\lq\lq $m(\{\omega\in\Omega_\Gamma\colon\, (\mathfrak f_\Gamma)_{xx}(\omega,\cdot)$ is strictly decreasing on $J\})>0$ for every
compact set $J\subset\mathbb{R}$ and every $\sigma$-ergodic measure $m$ in $\Omega_\Gamma$": in the words of the paper
\cite{dno1} and according to its Proposition~3.9, these two properties are equivalent to the $\mathrm{(SDC)}_*$
character of $\mathfrak f_\Gamma$, so that we name the property in the same way. The next lemma
is used in the proof of this property (in Proposition~\ref{prop:3SDCstar}) and in that of Propositions~\ref{prop:3three} and \ref{prop:3forwardcon}.
\begin{lemma}\label{lemm:3union}
Let $(f,\Gamma)$ satisfy {\rm\hyperlink{h1}{\textbf{h1}}} and {\rm\hyperlink{h2}{\textbf{h2}}}.
Let $\mathcal F_\Gamma=\{f_\Gamma{\cdot}s\colon\, s\in\mathbb{R}\}$ be the $\sigma$-orbit of $f_\Gamma$,
and let $\Omega_-$ and $\Omega_+$ be the $\sigma$-invariant subsets of $\Omega_\Gamma$ given by
corresponding $\boldsymbol\upalpha$-limit and $\boldsymbol\upomega$-limit
sets. Then, $\Omega_\Gamma=\Omega_-\cup\mathcal F_\Gamma\cup\Omega_+$.
In addition, $\Omega_-$ and $\Omega_+$ are respectively contained in the hulls
$\Omega_{\gamma_-}$ and $\Omega_{\gamma_+}$ of the functions $f_{\gamma_-}$ and $f_{\gamma_+}$.
\end{lemma}
\begin{proof}
We take $\w\in\Omega_\Gamma$ and a sequence $\{t_n\}$ such that
$\w(t,x)=\lim_{n\to\infty}f_\Gamma(t+t_n,x)=\lim_{n\to\infty}f(t+t_n,x,\Gamma(t+t_n))$ uniformly on the
compact subsets of $\mathbb R\times\mathbb R$. If there exists a subsequence $\{t_m\}$ of $\{t_n\}$ with limit
$-\infty$, then $\omega$ belongs to $\Omega_-$, and it is easy to check that
$\omega(t,x)=\lim_{m\to\infty}f_{\gamma_-}(t+t_n,x)$ uniformly on the
compact subsets of $\mathbb R\times\mathbb R$, so that $\w\in\Omega_{\gamma_-}$.
The situation is analogous if there exists a subsequence $\{t_m\}$ with limit $\infty$. Otherwise,
there exists a subsequence $\{t_m\}$ with limit $t_0\in\mathbb R$, and it is easy to check that
$\omega$ coincides with $f_\Gamma{\cdot}t_0\in\mathcal F_{\Gamma}$.
\end{proof}

\begin{proposition}\label{prop:3SDCstar}
Let $(f,\Gamma)$ satisfy {\rm\hyperlink{h1}{\textbf{h1}}}, {\rm\hyperlink{h2}{\textbf{h2}}},
{\rm\hyperlink{h4}{\textbf{h4}}} and  {\rm\hyperlink{h5}{\textbf{h5}}}.
Then, $\mathfrak f_\Gamma$ is $\mathrm{(SDC)}_*$.
\end{proposition}
\begin{proof}
We take $\w\in\Omega$ and any sequence $\{t_n\}\rightarrow\infty$ such that
$\omega=\lim_{n\to\infty}f_{\gamma_+}{\cdot}t_n$ (in the compact-open topology). Since $\omega_x$ is the limit of any subsequence of $\{(f_{\gamma_+})_x{\cdot}t_n\}$
which uniformly converges on the compact subsets of $\mathbb R\times\mathbb R$,
$\w_x=\lim_{n\to\infty}(f_{\gamma_+})_x{\cdot}t_n$. For the same reason,
$\w_{xx}=\lim_{n\to\infty}(f_{\gamma_+})_{xx}{\cdot}t_n$.

It follows easily that \hyperlink{h4}{\textbf{h4}} implies that $x\mapsto \mathfrak (f_\Gamma)_x(\omega,x)$ is concave for all $\omega\in\Omega$. Let us check the rest.
Let $m$ be a $\sigma$-ergodic measure on $\Omega_\Gamma$. We deduce from
Lemma \ref{lemm:3union} that either $m(\Omega_-)=1$ or
that $m(\Omega_+)=1$ (or both if they are not disjoint), since
$\mathcal F_\Gamma=\bigcup_{n\in\mathbb{Z}}
\sigma_n(\{f_\Gamma(t+s,x)\colon\, s\in[0,1)\})$ is a countable union of disjoint sets of the same measure $m$.
Hence, it suffices to check that $x\mapsto(\mathfrak f_\Gamma)_{xx}(\omega,x)$ is strictly decreasing whenever
$\omega\in\Omega_-\cup\Omega_+$.
We work in the case $\w\in\Omega_+\subseteq\Omega_{\gamma_+}$, writing
$\omega=\lim_{n\to\infty}f_{\gamma_+}{\cdot}t_n$ for
a sequence $\{t_n\}\rightarrow\infty$.
Let us take $x_1<x_2$. Condition \hyperlink{h5}{h5$_{\gamma_+}$} (see
Remark \ref{rm:twovariables}\hyperlink{311}{.1}) ensures that $\inf_{t\in\mathbb R} (f_{xx}(t,x_1,\gamma_+)-f_{xx}(t,x_2,\gamma_+))>0$. Hence,
\begin{equation*}
\begin{split}
&(\mathfrak f_\Gamma)_{xx}(\omega,x_1)-(\mathfrak f_\Gamma)_{xx}(\omega,x_2)=\omega_{xx}(0,x_1)-\omega_{xx}(0,x_2)\\
&\qquad\qquad=\lim_{n\rightarrow\infty}(f_{xx}(t_n,x_1,\gamma_+)-f_{xx}(t_n,x_2,\gamma_+))
>0\,,
\end{split}
\end{equation*}
so $x\mapsto (\mathfrak f_\Gamma)_{xx}(\omega,x)$ is strictly decreasing. An analogous argument
for $\omega\in\Omega_-$ completes the proof.
\end{proof}
\begin{proof}[Proof of Theorem~$\mathrm{\ref{th:3three}}$]
Let us assume (a).
Since $\tau(t,f_\Gamma,l_{\Gamma}(0))=(f_{\Gamma}{\cdot}t,l_\Gamma(t))$,
the set ${K}_{l_\Gamma}=\mathrm{cl}_{\Omega_\Gamma\times\mathbb{R}}\{(f_\Gamma{\cdot}t,l_\Gamma(t))\colon\, t\in\mathbb{R}\}$ is
$\tau_\Gamma$-invariant and compact: it is the closure of a bounded orbit. It is also easy to check that
it projects over $\Omega_\Gamma$, meaning that there exists at least a point $x_\w$ such that $(\w,x_\w)\in K_{l_\Gamma}$
for any $\omega\in\Omega_\Gamma$. In the same way, we define ${K}_{m_\Gamma}$ and ${K}_{u_\Gamma}$,
and observe that the three sets are disjoint (due to the uniform separation of $l_{\Gamma}, m_{\Gamma}$ and
$u_{\Gamma}$) and satisfy $x_l<x_m<x_u$ whenever $(\omega,x_l)\in K_{l_\Gamma}$, $(\omega,x_m)\in K_{m_\Gamma}$ and
$(\omega,x_u)\in K_{u_\Gamma}$. Propositions~\ref{prop:3coercivity} and \ref{prop:3SDCstar}, combined with
Proposition 3.11 of \cite{dno1}, guarantee that
$\mathfrak f_\Gamma$ satisfies the hypotheses of Theorem~4.2(ii) of \cite{dno1}. This result ensures that
$K_{l_\Gamma}$, $K_{m_\Gamma}$ and $K_{u_\Gamma}$ are ``hyperbolic copies of $\Omega_\Gamma$":
the $\tau_\Gamma$-invariant graphs of three real continuous maps
$\tilde{\mathfrak l}_\Gamma<\tilde{\mathfrak m}_\Gamma<\tilde{\mathfrak u}_\Gamma$ with
$K_{l_\Gamma}$ and $K_{u_\Gamma}$ uniformly asymptotically exponentially stable at $\infty$, and
$K_{m_\Gamma}$ uniformly asymptotically exponentially stable at $-\infty$. In the case of $K_{l_\Gamma}$,
this property means the existence of $k_l\ge 1$, $\beta_l>0$ and $\rho_l>0$ such that if $\omega\in\Omega_\Gamma$ and
$|\tilde{\mathfrak l}_\Gamma(\omega)-x_0|\le\rho_l$, then $|\tilde{\mathfrak l}_\Gamma(\omega{\cdot}t)-\mathfrak x_\Gamma(t,\omega,x_0)|\le
k_l e^{-\beta_l t}|\tilde{\mathfrak l}_{\Gamma}(\omega)-x_0|$ if $t\ge 0$. We call $\tilde l_\Gamma(t)=\tilde{\mathfrak l}_\Gamma(f_\Gamma{\cdot}t)$ for $t\in\mathbb R$,
so that $x_\Gamma(t,s,\tilde l_\Gamma(s))=\mathfrak x_\Gamma(t-s,f_\Gamma{\cdot}s,\tilde{\mathfrak l}_\Gamma(f_\Gamma{\cdot}s))=
\mathfrak x_\Gamma(t,f_\Gamma,\tilde{\mathfrak l}_\Gamma(f_\Gamma))=\tilde{\mathfrak l}_\Gamma(f_\Gamma{\cdot}t)=\tilde l_\Gamma(t)$:
the map $t\mapsto \tilde l_\Gamma(t)$ is a bounded solution of \eqref{eq:2completeproblem} and satisfies
$|\tilde l_\Gamma(t)-x_\Gamma(t,s,x_0)|=
|\tilde{\mathfrak l}_\Gamma(f_\Gamma{\cdot}s{\cdot}(t-s))-\mathfrak x_\Gamma(t-s,f_\Gamma{\cdot}s,x_0)|
\le k_l e^{-\beta_l(t-s)}|\tilde{\mathfrak l}_\Gamma(f_\Gamma{\cdot}s)-x_0|=k_l e^{-\beta_l(t-s)}|\tilde l_\Gamma(s)-x_0|$
for $t\ge s$ if $|\tilde l_\Gamma(s)-x_0|\le\rho_l$. Rewriting this inequality as
$|\tilde l_\Gamma(t)-x_\Gamma(t,s,x_0)|/|\tilde l_\Gamma(s)-x_0|\le k_l e^{-\beta_l(t-s)}$ and
taking limit as $x_0\to l_\Gamma(s)$, we get $\exp(\int_s^t f_x(r,\tilde l_\Gamma(r))\,dr)=
(\partial/\partial y) x_\Gamma(t,s,y)|_{y=\tilde l_\Gamma(s)}\le k_l e^{-\beta_l(t-s)}$ if $t\ge s$,
which shows the attractive hyperbolic character of $\tilde l_\Gamma$. Analogous definitions and arguments
apply to $\tilde m_\Gamma$ and $\tilde u_\Gamma$. Consequently, (c) holds.

It is obvious that (c) ensures (b). Let us assume (b) and deduce (a).
It is easy to deduce from Proposition \ref{prop:2hiperbolicas} and
Lemma \ref{lemm:3limit}(ii) that a solution cannot approach an attractive (resp. repulsive) hyperbolic solution
of the same equation as time decreases (resp. time increases),
so that an attractive hyperbolic solution and a
repulsive hyperbolic solution are uniformly separated.
Let us check that, if there exists three hyperbolic solutions $\tilde x_1<\tilde x_2<\tilde  x_3$ and
(at least) two of them are hyperbolic attractive, then $\tilde x_2$ is repulsive (which
combined with the previous assertion shows the uniform separation of the three of them).
So, for contradiction, we assume that (for
instance) $\tilde x_1$ and $\tilde x_2$ are attractive, and take
$x_0\in (\tilde x_1(0),\tilde x_2(0))$. Then, the graphs of the four solutions
$\tilde x_1$, $x_\Gamma(\cdot,0,x_0)$, $\tilde x_2$ and $\tilde x_3$ cannot approach
each other as time decreases, since any pair of consecutive ones contains a hyperbolic attractive solution.
It is easy to deduce that the $\boldsymbol\upalpha$-limit sets of $(f_\Gamma,\tilde x_1(0))$, $(f_\Gamma,x_0)$,
$(f_\Gamma,\tilde x_2(0))$ and $(f_\Gamma,\tilde x_3(0))$ are disjoint. In addition, all of them project over
the $\sigma$-invariant subset $\Omega_-\subseteq\Omega$ of Lemma \ref{lemm:3union}, and this
contradicts Theorem~4.2(ii) of \cite{dno1} (which establishes the existence of at most
three such compact sets). The case of attractive $\tilde x_2$
forces, as just seen, $\tilde x_1$ and $\tilde x_3$ to be repulsive and hence shows the
asserted uniform separation. And the existence of two consecutive repulsive hyperbolic solutions
is precluded with the previous argument, working with $\boldsymbol\upomega$-limit sets.
(In fact, these last two situations are impossible, as Theorem~4.2(ii) of \cite{dno1}
ensures.) Hence, (a) holds.

Let us check that the $\tau_\Gamma$-global attractor $\mathcal A_\Gamma$ (see Proposition
\ref{prop:3coercivity}) takes the form
$\mathcal A_\Gamma=\bigcup_{\omega\in\Omega_\Gamma}(\{\omega\}\times[\tilde{\mathfrak l}_\Gamma(\omega),
\tilde{\mathfrak u}_\Gamma(\omega)])$. We assume for contradiction the existence of $(\w_0,x_0)\in \mathcal A_\Gamma$
with $x_0>\tilde{\mathfrak u}_\Gamma(\omega_0)$ and deduce from the attractive hyperbolicity of $\tilde{\mathfrak u}_\Gamma$
that the $\boldsymbol\upalpha$-limit set of $(\w_0,x_0)$, which projects on a
$\sigma$-invariant compact subset $\Omega_0\subseteq\Omega$, is strictly above
$\{(\omega,\tilde{\mathfrak u}_\Gamma(\omega)):\,\omega\in\Omega_0\}$.
This means the existence of four different $\tau$-invariant compact sets projecting on $\Omega_0$, which contradicts
Theorem 4.2 of \cite{dno1}. An analogous argument shows that $x_0\ge \tilde{\mathfrak l}_\Gamma(\omega_0)$, which completes the proof of the equality (since the inverse contention is obvious).
Therefore, the last assertion in Proposition \ref{prop:3coercivity} shows that
$\tilde l_\Gamma=l_\Gamma$ and $\tilde u_\Gamma=u_\Gamma$.
Let us take $(s,x_0)$ with $x_0< \tilde m_\Gamma(s)=
\tilde{\mathfrak m}_\Gamma(f_\Gamma{\cdot}s)$. The $\boldsymbol\upomega$-limit set of
$(f_\Gamma{\cdot}s,x_0)$ cannot be contained into the
repulsive set $K_{m_\Gamma}$ and it is below it, so that it is contained into $K_{l_\Gamma}$.
This fact and the hyperbolic character of $K_{l_\Gamma}$ show that $\lim_{t\to\infty}
(x_\Gamma(t,s,x_0)-\tilde l_{\Gamma}(t))=\lim_{t\to\infty}
(\mathfrak x_\Gamma(t-s,f_{\Gamma}{\cdot}s,x_0)-\tilde{\mathfrak l}_{\Gamma}(f_\Gamma{\cdot}s{\cdot}(t-s)))=0$.
Analogous arguments (working with the $\boldsymbol\upalpha$-limit set in the case
$x_0\in(\tilde l_\Gamma(s),\tilde u_\Gamma(s))$) show the remaining limiting properties.
The last assertion of the theorem follows easily from the previous ones.
\end{proof}
\begin{remark}\label{rm:3siempre-dos}
A map $h\colon\mathbb{R}\times\mathbb{R}\rightarrow\mathbb{R}$ satisfying  \hyperlink{h2t}{\textbf{h2}$_*$}-\hyperlink{h5t}{\textbf{h5}$_*$} defines a skewproduct flow
given by the family of equations $x'=\mathfrak h(\omega{\cdot}t,x)$ for $\omega$ in the hull of $h$. The arguments of Propositions \ref{prop:3coercivity} and \ref{prop:3SDCstar} are enough to
repeat step by step the proof of Theorem \ref{th:3three} in order to get analogous conclusions for $x'=h(t,x)$.
\end{remark}
\begin{proof}[Proof of Proposition~$\mathrm{\ref{prop:3three}}$]
(i) As in the proof of Theorem~\ref{th:3three}, the three uniformly separated hyperbolic solutions of
$f_{\gamma_+}$ given by \hyperlink{h6t}{\textbf{h6}$_*$} define the three unique
$\tau_\gamma$-invariant compact sets $K_{l_{\gamma_+}}$, $K_{m_{\gamma_+}}$ and $K_{u_{\gamma_+}}$
projecting over the hull $\Omega_{\gamma_+}$ of $f_{\gamma_+}$,
which are hyperbolic copies of $\Omega_{\gamma_+}$: the graphs of
$\tilde{\mathfrak l}_{\gamma_+}<\tilde{\mathfrak m}_{\gamma_+}<\tilde{\mathfrak u}_{\gamma_+}$,
respectively, with $\tilde{\mathfrak l}_{\gamma_+}(f_{\gamma_+}{\cdot}s)=\tilde l_{\gamma_+}(s)$,
$\tilde{\mathfrak m}_{\gamma_+}(f_{\gamma_+}{\cdot}s)=\tilde m_{\gamma_+}(s)$
and $\tilde{\mathfrak u}_{\gamma_+}(f_{\gamma_+}{\cdot}s)=\tilde u_{\gamma_+}(s)$.
The uniform separation of $\tilde l_\Gamma<\tilde m_\Gamma<\tilde u_\Gamma$ ensures that the $\boldsymbol\upomega$-limit sets of $(f_\Gamma,\tilde l_\Gamma(0))$, $(f_\Gamma,\tilde m_\Gamma(0))$, and $(f_\Gamma,\tilde u_\Gamma(0))$ are disjoint and project onto $\Omega_+\subseteq\Omega_{\gamma_+}$
(see Lemma~\ref{lemm:3union}), and hence they are contained into $K_{l_{\gamma_+}}$, $K_{m_{\gamma_+}}$ and $K_{u_{\gamma_+}}$, respectively.
The final arguments of the proof of Theorem~\ref{th:3three} allow us to check that
$\lim_{t\rightarrow\infty}(x_\Gamma(t,s,x_0)-\tilde u_{\gamma_+}(t))=0$ for $x_0>\tilde m_\Gamma(s)$ and
$\lim_{t\rightarrow\infty}(x_\Gamma(t,s,x_0)-\tilde l_{\gamma_+}(t))=0$ for $x_0<\tilde m_\Gamma(s)$.

(ii) As in (i), we define $\tilde{\mathfrak l}_{\gamma_-}<\tilde{\mathfrak m}_{\gamma_-}<\tilde{\mathfrak u}_{\gamma_-}$ on $\Omega_{\gamma_-}$
and $K_{l_{\gamma_-}}$, $K_{m_{\gamma_-}}$ and $K_{u_{\gamma_-}}$.
The existence of a bounded solution of \eqref{eq:2completeproblem} strictly above $u_\Gamma$ or below $l_\Gamma$ would give rise to a fourth
$\tau$-invariant compact set projecting onto $\Omega_-$ (given by the $\boldsymbol\upomega$-limit set of the corresponding $\tau_\Gamma$-orbit), which contradicts Theorem 4.2 of \cite{dno1}. A similar argument shows the remaining assertion in (ii).
\end{proof}
We recall that Theorem~5.1(i) of \cite{dno1} ensures that $\mathfrak l_\Gamma$ and $\mathfrak u_\Gamma$ are respectively lower and upper semicontinuous maps.
The last result in this section provides a new characterization of {\sc Case \hyperlink{A}{A}} in terms of the global
attractor $\mathcal A_\Gamma$.
\begin{proposition} \label{prop:3forwardcon}
Let $(f,\Gamma)$ satisfy {\rm\hyperlink{h1}{\textbf{h1}}-\hyperlink{h6}{\textbf{h6}}}. The equation \eqref{eq:2completeproblem}
is in \textsc{Case \hyperlink{A}{A}} if and only if the delimiters
$\mathfrak l_\Gamma,\mathfrak u_\Gamma\colon\Omega_\Gamma\rightarrow\mathbb{R}$
of the global attractor $\mathcal{A}_\Gamma$ are continuous.
\end{proposition}
\begin{proof} Necessity has been established in the proof of Theorem \ref{th:3three}: \textsc{Case \hyperlink{A}{A}} ensures
that $\mathcal A_\Gamma=\bigcup_{\omega\in\Omega_\Gamma}(\{\omega\}\times[\tilde{\mathfrak l}_\Gamma(\omega),
\tilde{\mathfrak u}_\Gamma(\omega)])$ with
$\tilde{\mathfrak l}_\Gamma,\tilde{\mathfrak u}_\Gamma\colon\Omega_\Gamma\to\mathbb R$ continuous.

To show sufficiency, we assume the continuity of the delimiters. As in the proof of Theorem \ref{th:3three},
we define $\tilde{\mathfrak l}_{\gamma_+}$ and
$\tilde{\mathfrak u}_{\gamma_+}$ as the continuous functions defined on the hull of $f_{\gamma_+}$ which
give rise to the corresponding global attractor, and which are uniformly separated.
The characterization of these sets in terms of
bounded orbits shows that $[\mathfrak l_\Gamma(\w),\mathfrak u_\Gamma(\w)]=
[\tilde{\mathfrak l}_{\gamma_+}(\omega),\tilde{\mathfrak u}_{\gamma_+}(\omega)]$ for $\omega\in\Omega_+$.
We fix $\omega\in\Omega_+$ and
$\{s_n\}\uparrow\infty$ with $\omega=\lim_{n\to\infty}f_\Gamma{\cdot}s_n$.
It follows from $\lim_{n\to\infty}l_\Gamma(s_n)=\lim_{n\to\infty}\mathfrak l_\Gamma(f_\Gamma{\cdot}s_n)=
\mathfrak l_\Gamma(\w)=\tilde{\mathfrak l}_{\gamma_+}(\omega)$ and from the hyperbolicity of
the graph of $\tilde{\mathfrak l}_{\gamma_+}$ that there exists $s>0$ such that $\lim_{t\to\infty}(l_\Gamma(t)-
\tilde{\mathfrak l}_{\gamma_+}((\omega{\cdot}(-s)){\cdot}t))=0$. Theorem \ref{th:3general}(i) and Lemma \ref{lemm:3limit}(iii)
ensure that $l_\Gamma$ is hyperbolic attractive. The same holds for $u_\Gamma$. Hence, there
exist two attractive hyperbolic solutions, and their limiting properties show that they
are uniformly separated, which means {\sc Case \hyperlink{A}{A}}.
\end{proof}
Notice that the proof of Proposition~\ref{prop:3forwardcon} also shows that $\mathfrak u_\Gamma$ (resp. $\mathfrak l_\Gamma$) is continuous
in \textsc{Cases \hyperlink{B1}{B1}} and \textsc{\hyperlink{C1}{C1}} (resp. \textsc{Cases \hyperlink{B2}{B2}} and \textsc{\hyperlink{C2}{C2}}).
\section{Some mechanisms producing critical transitions}\label{sec:mechanisms}
Different mechanisms have been described by applied scientists as sources of critical transitions. In this section, we assume that $(f,\Gamma)$ satisfies \hyperlink{h1}{\textbf{h1}}-\hyperlink{h6}{\textbf{h6}} (in the sense explained in Remark \ref{rm:twovariables}\hyperlink{311}{.1}), and we study these three different types of one-parametric perturbations of \eqref{eq:2completeproblem}:
\begin{itemize}
\item \emph{Rate-induced tipping/tracking:} we study the values of the rate $c>0$ for which $x'=f(t,x,\Gamma(ct))$ presents tipping or tracking. The parameter $c$ determines the speed of the transition: larger values of $c$ mean a faster transition from the past $\gamma_-$ to the future $\gamma_+$ (see \cite{aspw}, \cite{lno1}).
\item \emph{Phase-induced tipping/tracking:} we study the values of the phase $c\in\mathbb{R}$ for which $x'=f(t,x,\Gamma(t+c))$ presents tipping or tracking. The parameter $c$ represents the phase of the initial state from which the transition departs (see \cite{alas}, \cite{altw}, \cite{lno2}).
\item \emph{Size-induced tipping/tracking:} we study the values of the size $d\in\mathbb{R}$ for which $x'=f(t,x,d\,\Gamma(t))$ presents tipping or tracking.
It is usually assumed that $\lim_{t\rightarrow-\infty} \Gamma(t)=0$ and $\lim_{t\rightarrow\infty}\Gamma(t)=1$. The parameter $d$ represents the distance between the steady states of the past and future equations, that is, the size of the transition. Special attention has been paid to problems of the type $x'=h(t,x-d\,\Gamma(t))$ (see \cite{lno2}). Notice that, in this particular case, the phase space of the future equation $x'=h(t,x-d)$ is a lift of that of the past one $x'=h(t,x)$, and hence they have the same number and type of uniformly separated hyperbolic solutions.
\end{itemize}
\par
It is interesting to observe that, while $\lim_{t\rightarrow\pm\infty} \Gamma(ct)=\gamma_\pm$ and $\lim_{t\rightarrow\pm\infty} \Gamma(t+c)=\gamma_\pm$ for all $c$, this is not the case
for $\lim_{t\rightarrow\pm\infty} d\,\Gamma(t)$ unless $\gamma_{\pm}=0$. Therefore,
to use the results of Subsection~\ref{subsec:2transition} in the size-induced tipping/tracking analysis, we must restrict ourselves to the set of parameters $d$ such that $(f,d\,\Gamma)$ satisfies \hyperlink{h1}{\textbf{h1}}-\hyperlink{h6}{\textbf{h6}}.
However, in this section, we will just focus on the size-induced problem given by $x'=h(t,x-d\,\Gamma(t))$ assuming that $h$ satisfies
\hyperlink{h2t}{\textbf{h2}$_*$}-\hyperlink{h5t}{\textbf{h5}$_*$}, which guarantees
the same properties for the all the maps $(t,x)\mapsto h(t,x-d\Gamma(t))$.

\begin{remark}\label{rm:nonincreasingchange}
The main results of this section are formulated for the particular case in which $\gamma\mapsto f(t,x,\gamma)$ is nondecreasing for all $(t,x)\in\mathbb{R}\times\mathbb{R}$.
We remark that, to get conclusions for the equation $x'=f(t,x,\Gamma(t))$ in the case of nonincreasing $\gamma\mapsto f(t,x,\gamma)$ (which, of course, does not exhaust the possibilities), we can rewrite it as $x'=\tilde f(t,x,-\Gamma(t))$
for $\tilde f(t,x,\gamma)=f(t,x,-\gamma)$, having in mind that the pair $(\tilde f,-\Gamma)$ satisfies \hyperlink{h1}{\textbf{h1}}-\hyperlink{h6}{\textbf{h6}} if $\Gamma$ satisfies \hyperlink{h1}{\textbf{h1}} and $f$ satisfies \hyperlink{h2}{\textbf{h2}}-\hyperlink{h4}{\textbf{h4}}, \hyperlink{h5}{\textbf{h5}$_{-\gamma_\pm}$} and \hyperlink{h6}{\textbf{h6}$_{-\gamma_\pm}$}.
\end{remark}

The main results of this section are the following: Theorem~\ref{th:4compacity} precludes tipping in the case in which $x'=f(t,x,\gamma)$ has three uniformly separated hyperbolic solutions for every value of the parameter $\gamma$ in the image of $\Gamma$;
Theorem~\ref{th:4Gammad} includes certain transition equations as a part of a size-induced tipping phenomena;
and Theorems~\ref{th:4bifurcationfunction} and \ref{th:4bifurcationfunction2} construct continuous bifurcation functions given by the unique tipping point for certain rate-induced and phase-induced tipping problems.

The next result plays a fundamental role in the proof of Theorem~\ref{th:4compacity}.
\begin{proposition}\label{prop:3monotonicity} Let $h_1,h_2\colon\mathbb{R}\times\mathbb{R}\rightarrow\mathbb{R}$ satisfy {\rm\hyperlink{h2t}{\textbf{h2}$_*$}}-{\rm\hyperlink{h3t}{\textbf{h3}$_*$}} and $h_1(t,x)\leq h_2(t,x)$ for all $(t,x)\in\mathbb{R}\times\mathbb{R}$.
\begin{enumerate}[label=\rm{(\roman*)}]
\item Let $l_i$ (resp. $u_i$) be the lower (resp. upper) bounded solution of $x'=h_i(t,x)$ for $i=1,2$. Then, $l_1\leq l_2$ and $u_1\leq u_2$.
\item Assume that $h_1$ and $h_2$ also satisfy {\rm\hyperlink{h4t}{\textbf{h4}$_*$}}-{\rm\hyperlink{h6t}{\textbf{h6}$_*$}}, and let $\tilde l_i<\tilde m_i<\tilde u_i$ be the hyperbolic solutions of $x'=h_i(t,x)$ for $i=1,2$ provided by {\rm\hyperlink{h6t}{\textbf{h6}$_*$}}. Then, $\inf_{t\in\mathbb R}(\tilde m_1(t)-\tilde l_2(t))>0$ if and only if $\inf_{t\in\mathbb R}(\tilde u_1(t)-\tilde m_2(t))>0$, in which case $\tilde l_1\leq \tilde l_2<\tilde m_2\leq \tilde m_1<\tilde u_1\leq \tilde u_2$.
    If, in addition, $h_1(t,x)< h_2(t,x)$ for all $(t,x)\in\mathbb{R}\times\mathbb{R}$, then all the inequalities are strict.
\end{enumerate}
\end{proposition}
\begin{proof} (i) Since $h_1(s,u_1(s))\le h_2(s,u_1(s))$ for all $s\in\mathbb R$, a standard argument of comparison of solutions shows that $x_2(t,s,u_1(s))\le u_1(t)$ whenever $s\in\mathbb R$ and $t\le s$. Hence, $x_2(t,s,u_1(s))$ remains bounded from above as time decreases, which means that $u_1(t)\le u_2(t)$ for all $t\in\mathbb R$: see Section~\ref{sec:2tipping}. Analogously, $l_1\le l_2$.

(ii) Let us assume $\inf_{t\in\mathbb R}(\tilde m_1(t)-\tilde l_2(t))>0$: the argument is analogous if $\inf_{t\in\mathbb R}(\tilde u_1(t)-\tilde m_2(t))>0$. We define $m_s(t)=x_2(t,s,\tilde m_1(s))$ (which always exists since $\tilde l_2<\tilde m_1\le\tilde u_1\le\tilde u_2$) and observe that $m_s(t)>\tilde l_2(t)$ for all $s,t\in\mathbb R$.
A standard argument of comparison of solutions shows that $x_2(t,s,\tilde m_1(s))\le\tilde m_1(t)$ if $t\le s$. In addition, if $s_1<s_2$, then
$m_{s_1}\ge m_{s_2}$: $m_{s_2}(t)=x_2(t,s_2,\tilde m_1(s_2))=x_2(t,s_1,x_2(s_1,s_2,\tilde m_1(s_2)))\le x_2(t,s_1,\tilde m_1(s_1))=m_{s_1}(t)$.
Therefore, there exists $m_\infty(t)=\lim_{s\to\infty}m_s(t)\in[\tilde l_2(t),\tilde m_1(t)]$ for all $t\in\mathbb R$. It is easy to check that
$t\mapsto m_\infty(t)$ is a bounded solution of $x'=h_2(t,x)$. Our goal is to check that $\lim_{t\to\infty}|m_\infty(t)-\tilde m_2(t)|=0$, which, according
to Theorem \ref{th:3three} and Remark \ref{rm:3siempre}, ensures that $m_\infty=\tilde m_2$ and hence that $\tilde m_2\le\tilde m_1$. This inequality combined with (i) proves that $\tilde l_1\leq \tilde l_2<\tilde m_2\leq \tilde m_1<\tilde u_1\leq \tilde u_2$. In turn, this
chain of inequalities combined with the uniform separation of the hyperbolic solutions $\tilde l_i<\tilde m_i<\tilde u_i$ (see Theorem~\ref{th:3three}) ensures that $\inf_{t\in\mathbb R}(\tilde u_1(t)-\tilde m_2(t))>0$.

The inequalities $m_\infty\le\tilde m_1<\tilde u_1\le\tilde u_2$ and $\inf_{t\in\mathbb R}(\tilde u_1(t)-\tilde m_1(t))>0$
preclude $\lim_{t\to\infty}|m_\infty(t)-\tilde u_2(t)|=0$.
Hence, the unique possibility to be excluded (see again Theorem \ref{th:3three})
is $\lim_{t\to\infty}|m_\infty(t)-\tilde l_2(t)|=0$, which we assume for contradiction. Let $(k,\beta)$ be a dichotomy constant pair for
the hyperbolic solution $\tilde l_2$ of $x'=h_2(t,x)$, and let $\rho>0$ be the constant associated by Proposition \ref{prop:2hiperbolicas}
to $\beta/2$. We look for $t_0\in\mathbb R$ such that $|m_\infty(t_0)-\tilde l_2(t_0)|\le \rho/2$, and $s_0\ge t_0$ such that
$|\tilde l_2(t_0)-x_2(t_0,s,\tilde m_1(s)))|\le\rho$ for all $s\ge s_0$. Then,
$|\tilde l_2(s)-\tilde m_1(s)|=|\tilde l_2(s)-x_2(s,t_0,x_2(t_0,s,\tilde m_1(s)))|\le ke^{-\beta/2(s-t_0)}\rho$ for all $s\ge s_0$.
Taking limit as $s\to\infty$ yields $\inf_{t\in\mathbb R}(\tilde m_1(t)-\tilde l_2(t))=0$, which is the sought-for contradiction.

The last assertion follows easily from contradiction and comparison. For instance, if
$m_1(s)=m_2(s)$, then $m_1(t)=x_1(t,s,m_1(s))<x_2(t,s,m_2(s))=m_2(t)$ for $t>s$, which is not the case. \end{proof}
Recall that $t\mapsto x_\Gamma(t,s,x_0)$ is the solution of \eqref{eq:2completeproblem} satisfying $x_\Gamma(s,s,x_0)=x_0$ and that $\tilde l_\gamma<\tilde m_\gamma<\tilde u_\gamma$ represent the three hyperbolic solutions of $x'=f(t,x,\gamma)$ provided by {\rm\hyperlink{h6}{\textbf{h6}$_{\gamma}$}}.

\begin{proposition}\label{rmk:morerange} Let $(f,\Gamma)$ satisfy {\rm\hyperlink{h1}{\textbf{h1}}-\hyperlink{h6}{\textbf{h6}}} with
$\gamma\mapsto f(t,x,\gamma)$ nondecreasing for all $(t,x)\in\mathbb{R}\times\mathbb{R}$, and define $[\bar\gamma_m,\bar\gamma_M]=\mathrm{cl}(\Gamma(\mathbb{R}))$. Assume also that $f$ satisfies {\rm\hyperlink{h6}{\textbf{h6}$_{\bar\gamma_m}$}} and {\rm\hyperlink{h6}{\textbf{h6}$_{\bar\gamma_M}$}}, that $\inf_{t\in\mathbb{R}} (\tilde m_{\gamma_+}(t)-\tilde l_{\bar\gamma_M}(t))>0$ (or, equivalently, $\inf_{t\in\mathbb{R}} (\tilde u_{\gamma_+}(t)-\tilde m_{\bar\gamma_M}(t))>0$), and that
  $\inf_{t\in\mathbb{R}} (\tilde m_{\bar\gamma_m}(t)-\tilde l_{\gamma_+}(t))>0$ (or, equivalently, $\inf_{t\in\mathbb{R}} (\tilde u_{\bar\gamma_m}(t)-\tilde m_{\gamma_+}(t))>0$).
Then, \eqref{eq:2completeproblem} is in \textsc{Case \hyperlink{A}{A}}.
\end{proposition}
\begin{proof} The fact that $\inf_{t\in\mathbb{R}} (\tilde m_{\gamma_+}(t)-\tilde l_{\bar\gamma_M}(t))>0$ allows us to apply Proposition~\ref{prop:3monotonicity}(ii) to $f_{\gamma_+}\leq f_{\bar\gamma_M}$ and conclude that $\tilde l_{\gamma_+}\leq\tilde l_{\bar\gamma_M}<\tilde m_{\bar\gamma_M}\leq\tilde m_{\gamma_+}<\tilde u_{\gamma_+}\leq \tilde u_{\bar\gamma_M}$.
Proposition~\ref{prop:3monotonicity}(i) applied to $f_\Gamma\leq f_{\bar\gamma_M}$ ensures that $l_\Gamma\leq \tilde l_{\bar\gamma_M}$. The previous facts and $\inf_{t\in\mathbb{R}}(\tilde m_{\bar\gamma_M}(t)-\tilde l_{\bar\gamma_M}(t))>0$ ensure that $\inf_{t\in\mathbb{R}}(\tilde m_{\gamma_+}(t)-l_\Gamma(t))>0$. Hence,
Theorem~\ref{th:3general}(ii) ensures that $\lim_{t\rightarrow\infty}(l_\Gamma(t)-\tilde l_{\gamma_+}(t))=0$. Similarly, $\inf_{t\in\mathbb{R}} (\tilde m_{\bar\gamma_m}(t)-\tilde l_{\gamma_+}(t))>0$, Proposition~\ref{prop:3monotonicity}(ii) applied to $f_{\bar\gamma_m}\leq f_{\gamma_+}$ and Proposition~\ref{prop:3monotonicity}(i) applied to $f_{\bar\gamma_m}\leq f_\Gamma$ yield $\lim_{t\rightarrow\infty}(u_\Gamma(t)-\tilde u_{\gamma_+}(t))=0$. Hence, the result follows from Theorem \ref{th:2massobrecasos}(1).
\end{proof}
The following relevant theorem is a consequence of the previous results.
\begin{theorem}\label{th:4compacity} Let $(f,\Gamma)$ satisfy {\rm\hyperlink{h1}{\textbf{h1}}-\hyperlink{h4}{\textbf{h4}}} with
$\gamma\mapsto f(t,x,\gamma)$ nondecreasing for all $(t,x)\in\mathbb{R}\times\mathbb{R}$. Assume also that $f$ satisfies {\rm\hyperlink{h5}{\textbf{h5}$_\gamma$} and \hyperlink{h6}{\textbf{h6}$_\gamma$}} for all $\gamma\in\mathrm{cl}(\Gamma(\mathbb{R}))$. Then, \eqref{eq:2completeproblem} is in \textsc{Case \hyperlink{A}{A}}. In particular, under these hypotheses, equations $x'=f(t,x,\Gamma(ct))$ and $x'=f(t,x,\Gamma(t+c))$ are in \textsc{Case \hyperlink{A}{A}} for all $c\in\mathbb R$: there is neither rate-induced nor phase-induced tipping under these hypotheses.
\end{theorem}
\begin{proof} The hypotheses on the strict concavity of $f_x$ and the number of hyperbolic minimal sets ensure that the dynamics of $x'=f(t,x,\gamma)$ is that of Theorem~\ref{th:3three} for all $\gamma\in\mathrm{cl}(\Gamma(\mathbb{R}))$. Given $\gamma_1\in\mathrm{cl}(\Gamma(\mathbb{R}))$, we deduce from Theorem \ref{th:hypcontinuation} that $\inf_{t\in\mathbb{R}}(\tilde u_{\gamma_1}(t)-\tilde m_{\gamma_2}(t))>0$ and $\inf_{t\in\mathbb{R}}(\tilde m_{\gamma_2}(t)-\tilde l_{\gamma_1}(t))>0$ if $\gamma_2$ is close enough to $\gamma_1$, and hence, applying Proposition \ref{prop:3monotonicity}(ii), $\tilde l_{\gamma_1}\leq \tilde l_{\gamma_2}<\tilde m_{\gamma_2}\leq \tilde m_{\gamma_1}<\tilde u_{\gamma_1}\leq \tilde u_{\gamma_2}$ if $\gamma_1<\gamma_2$. Since $\mathrm{cl}(\Gamma(\mathbb{R}))$ is compact and connected, the previous inequalities are valid for all $\gamma_1,\gamma_2\in\mathrm{cl}(\Gamma(\mathbb{R}))$ with $\gamma_1<\gamma_2$. Note also that the strict inequalities correspond to uniform separation. It follows easily that the hypotheses Proposition~\ref{rmk:morerange} hold, so \eqref{eq:2completeproblem} is in \textsc{Case \hyperlink{A}{A}}. The last assertion is an immediate consequence.
\end{proof}

The proof of the following lemma, used in the proof of Theorem~\ref{th:4Gammad}, follows easily from the description of \textsc{Cases \hyperlink{A}{A}}, \hyperlink{B}{\textsc{B}} and \hyperlink{C}{\textsc{C}} made in the previous section.
\begin{lemma}\label{lemma:4boundsandlimits} Let $(f,\Gamma)$ satisfy {\rm\hyperlink{h1}{\textbf{h1}}-\hyperlink{h6}{\textbf{h6}}}. Let $(s,x_0)\in\mathbb{R}\times\mathbb{R}$, and let $\bar x\colon\mathbb{R}\rightarrow\mathbb{R}$ be a continuous function. Assume that $\bar x(t)\leq x_\Gamma(t,s,x_0)$ (resp. $\bar x(t)\geq x_\Gamma(t,s,x_0)$) for all $t\geq s$.
\begin{enumerate}[label=\rm{(\roman*)}]
\item If $\lim_{t\rightarrow\infty} (\bar x(t)-\tilde u_{\gamma_+}(t))=0$ (resp. $\lim_{t\rightarrow\infty} (\bar x(t)-\tilde l_{\gamma_+}(t))=0$), then $\lim_{t\rightarrow\infty} (x_\Gamma(t,s,x_0)-\tilde u_{\gamma_+}(t))=0$ (resp. $\lim_{t\rightarrow\infty} (x_\Gamma(t,s,x_0)-\tilde l_{\gamma_+}(t))=0$).
\item If $\lim_{t\rightarrow\infty} (\bar x(t)-\tilde m_{\gamma_+}(t))=0$, then $\lim_{t\rightarrow\infty} (x_\Gamma(t,s,x_0)-\tilde u_{\gamma_+}(t))=0$ (resp. $\lim_{t\rightarrow\infty} (x_\Gamma(t,s,x_0)-\tilde l_{\gamma_+}(t))=0$) or $\lim_{t\rightarrow\infty} (x_\Gamma(t,s,x_0)-\tilde m_{\gamma_+}(t))=0$.
\end{enumerate}
The same statements hold for $t\rightarrow-\infty$ and the hyperbolic solutions $\tilde l_{\gamma_-}$, $\tilde m_{\gamma_-}$ and $\tilde u_{\gamma_-}$ if we assume that $t\mapsto x_\Gamma(t,s,x_0)$ is a bounded solution and that $\bar x(t)\leq x_\Gamma(t,s,x_0)$ (resp. $\bar x(t)\geq x_\Gamma(t,s,x_0)$) for all $t\leq s$.
\end{lemma}
In the next result, taking as starting point a pair $(f,\Gamma)$ which does not necessarily  satisfy the hypotheses of Theorem~\ref{th:4compacity}, we design a one-parametric perturbation of \eqref{eq:2completeproblem} for which all the dynamical possibilities occur as the parameter
varies in $\mathbb R$. The sets $\mathcal{R}_f$ and $\mathfrak C_f$ appearing in the statement are defined by \eqref{def.R} and \eqref{def.C}.
\begin{theorem} \label{th:4Gammad}
Let $(f,\Gamma)$ satisfy {\rm\hyperlink{h1}{\textbf{h1}}-\hyperlink{h6}{\textbf{h6}}}, define
\begin{equation}\label{def:gama+-}
\gamma_m=\min\{\gamma_-,\gamma_+\}\quad\text{and}\quad\gamma_M=\max\{\gamma_-,\gamma_+\}\,,
\end{equation}
and assume that
\begin{enumerate}[label=\rm{(\arabic*)}]
\item $\sup_{t\in\mathbb{R}}\Gamma(t)>\gamma_M$ (resp. $\inf_{t\in\mathbb{R}}\Gamma(t)<\gamma_m$),
\item $[\inf_{t\in\mathbb R}\Gamma(t),\gamma_M]\subseteq\mathcal R_f$ (resp.~$[\gamma_m,\sup_{t\in\mathbb R}\Gamma(t)]\subseteq\mathcal R_f$),
\item the map $\gamma\mapsto f(t,x,\gamma)$ is strictly increasing for all $(t,x)\in\mathbb R\times\mathbb R$, and $\lim_{\gamma\rightarrow\pm\infty}f(t,x,\gamma)=\pm\infty$
    uniformly on compact sets of $\mathbb{R}\times\mathbb{R}$.
\end{enumerate}
Let $\Delta_1(t)=\min\{\gamma_M,\Gamma(t)\}$ (resp. $\Delta_1(t)=\max\{\gamma_m,\Gamma(t)\}$), let $\Delta_2(t)=\Gamma(t)-\Delta_1(t)$,
and let $\Gamma_d(t)=\Delta_1(t)+d\,\Delta_2(t)$. Then,
\begin{itemize}
\item[\rm(i)] there exist $d_\Gamma^-<0<d_\Gamma^+$ such that the equation
\begin{equation}\label{eq:4fconGammad}
x'=f(t,x,\Gamma_d(t))
\end{equation}
is in \textsc{Case \hyperlink{A}{A}} for $d\in (d_\Gamma^-,d_\Gamma^+)$, in \textsc{Case \hyperlink{B1}{B1}} (resp. \textsc{\hyperlink{B2}{B2}}) for $d=d_\Gamma^+$, in \textsc{Case \hyperlink{B2}{B2}} (resp. \textsc{\hyperlink{B1}{B1}}) for $d=d_\Gamma^-$, in \textsc{Case \hyperlink{C1}{C1}} (resp. \textsc{\hyperlink{C2}{C2}}) for $d\in(d_\Gamma^+,\infty)$, and in \textsc{Case \hyperlink{C2}{C2}} (resp. \textsc{\hyperlink{C1}{C1}}) for $d\in(-\infty,d_\Gamma^-)$.
\item[\rm(ii)]
Let $\mathfrak{C}_f$ be endowed with the $\|\cdot\|_\infty$ norm.
The maps $\mathcal C\to\mathbb R,\;\Lambda \mapsto d_\Lambda^\pm$ are well-defined and continuous on the set
$\mathcal C=\{\Lambda\in\mathfrak C_f\colon \text{$\Lambda$ satisfies {\rm(1)} and {\rm(2)}}\}$.
\item[\rm(iii)] If $\Lambda^1,\Lambda^2\in\mathcal C$ and $\Lambda^1\leq\Lambda^2$ (resp. $\Lambda^2\leq\Lambda^1$), then $d_{\Lambda^1}^-\leq d_{\Lambda^2}^-<0<d_{\Lambda^2}^+\leq d_{\Lambda^1}^+$; and if, in addition, there exists $t_*\in\mathbb{R}$ such that $\max(\lim_{t\to\pm\infty}\Lambda_1(t))<\Lambda^1(t_*)<\Lambda^2(t_*)$ (resp. $\Lambda^2(t_*)<\Lambda^1(t_*)<\min(\lim_{t\to\pm\infty}\Lambda_1(t))$), then $d_{\Lambda^1}^-<d_{\Lambda^2}^-<0<d_{\Lambda^2}^+<d_{\Lambda^1}^+$.
\end{itemize}
\end{theorem}
\begin{proof}
Since $\lim_{t\rightarrow\pm\infty}\Delta_1(t)=\gamma_\pm$ and $\lim_{t\rightarrow\pm\infty}\Delta_2(t)=0$, $\lim_{t\rightarrow\pm\infty}\Gamma_d(t)=\gamma_\pm$ for all $d\in\mathbb{R}$. That is,
all the pairs $(f,\Gamma_d)$ satisfy conditions \hyperlink{h1}{\textbf{h1}}-\hyperlink{h6}{\textbf{h6}}.
We write \eqref{eq:4fconGammad}$_d$ to make reference to equation \eqref{eq:4fconGammad} for the particular value of the parameter $d$, and call
$l_{\Gamma_d},m_{\Gamma_d}$ and $u_{\Gamma_d}$ the solutions of \eqref{eq:4fconGammad}$_d$ provided by Theorem \ref{th:3general} (or $\tilde l_{\Gamma_d},\tilde m_{\Gamma_d}$ and $\tilde u_{\Gamma_d}$ when we know that they are hyperbolic).

We will reason in the case $\sup_{t\in\mathbb{R}}\Gamma(t)>\gamma_M$. The arguments can be easily adapted to the case $\inf_{t\in\mathbb{R}}\Gamma(t)<\gamma_m$ having in mind this fundamental difference: if $\sup_{t\in\mathbb{R}}\Gamma(t)>\gamma_M$, then $\Delta_2(t)\geq 0$ for all $t\in\mathbb{R}$, so that $d\to\Gamma_d$ is a nondecreasing map, and hence the increasing character of $f$ guaranteed by (3) ensures that $d\mapsto f(t,x,\Gamma_d(t))$ is nondecreasing for all $(t,x)\in\mathbb{R}$; whereas, in the case $\inf_{t\in\mathbb{R}}\Gamma(t)<\gamma_m$, the map $d\mapsto f(t,x,\Gamma_d(t))$ is nonincreasing for all $(t,x)\in\mathbb{R}$.

(i) Observe first that condition (2) allows us to deduce from Theorem~\ref{th:4compacity} that \eqref{eq:4fconGammad}$_0$ is in \textsc{Case \hyperlink{A}{A}}. Let us fix $d\ge 0$. Proposition~\ref{prop:3monotonicity}(i) and Theorem \ref{th:3three} yield $\tilde u_{\Gamma_0}\leq u_{\Gamma_d}$ and $\lim_{t\to\infty}(\tilde u_{\Gamma_0}(t)-\tilde u_{\gamma_+}(t))=0$, and hence Lemma~\ref{lemma:4boundsandlimits}(i) applied to $\bar x=\tilde u_{\Gamma_0}$ guarantees that $\lim_{t\rightarrow\infty}(u_{\Gamma_d}(t)-\tilde u_{\gamma_+}(t))=0$.
Therefore, Theorem~\ref{th:2massobrecasos} ensures that \eqref{eq:4fconGammad}$_d$ is in \textsc{Cases} \hyperlink{A}{A}, \hyperlink{B1}{B1} or \hyperlink{C1}{C1} if $d\ge 0$.
Analogously, if $d\le 0$, since $\tilde l_{\Gamma_d}\leq \tilde l_{\Gamma_0}$ and $\lim_{t\to-\infty}(\tilde l_{\Gamma_0}(t)-\tilde l_{\gamma_-}(t))=0$, we have $\lim_{t\rightarrow-\infty}(\tilde l_{\Gamma_d}(t)-\tilde l_{\gamma_-}(t))=0$, and Theorem~\ref{th:2massobrecasos-2} ensures that \eqref{eq:4fconGammad}$_d$ is in \textsc{Cases} \hyperlink{A}{A}, \hyperlink{B2}{B2} or \hyperlink{C2}{C2}.

The condition $\sup_{t\in\mathbb R}\Gamma(t)>\gamma_M$ provides a nondegenerate interval $[a,b]$ and a constant $\delta>0$ such that $\Delta_2(t)\ge \delta$ for all $t\in[a,b]$. In particular, $\Delta_1(t)=\gamma_M$ for all $t\in[a,b]$. We define $d_\Gamma^+=\sup\{d\geq0\colon\, $\eqref{eq:4fconGammad}$_{\bar d}$ is in \textsc{Case \hyperlink{A}{A}} for all $\bar d\in[0,d]\}\ge 0$ and deduce from Proposition \ref{prop:3CaseCrobust} that $d_\Gamma^+>0$. Our first goal is to check that $d_\Gamma^+$ is finite. More precisely, we will find $d_0=d_0(\delta,a,b,\gamma_M)$ such that $d_\Gamma^+\le d_0$. The arguments of the proof of
Theorem~\ref{th:4compacity} show that, if $d\in[0,d_\Gamma^+)$, then $\tilde l_{\Gamma_0}\leq\tilde l_{\Gamma_d}<\tilde m_{\Gamma_d}\leq\tilde m_{\Gamma_0}<\tilde u_{\Gamma_0}\leq\tilde u_{\Gamma_d}$. We take $k>0$ such that $-k\le \tilde l_{\Gamma_0}(t)<\tilde m_{\Gamma_0}(t)\le k$ for all $t\in\mathbb R$.
According to the Mean Value Theorem, as long as $\tilde m_{\Gamma_d}$ exists on $[a,b]$ (which is the case if $d\in[0,d_\Gamma^+)$), there exists
$c_d\in[a,b]$ such that
\[
 \tilde m_{\Gamma_d}(a)
 =\tilde m_{\Gamma_d}(b)-(b-a)\tilde m_{\Gamma_d}'(c_d)
 =\tilde m_{\Gamma_d}(b)-(b-a) f(c_d,\tilde m_{\Gamma_d}(c_d),\gamma_M+d\Delta_2(c_d))\,.
\]
Condition (3) provides $d_0>0$ large enough to ensure that $f(t,x,\gamma_M+d_0\delta)> 2k/(b-a)$ for all $(t,x)\in[a,b]\times[-k,k]$. To check that $d_\Gamma^+\le d_0$, we assume  that this is not the case, to get $\tilde m_{\Gamma_{d_0}}(a)\le k-(b-a)f(c_{d_0},\tilde m_{\Gamma_{d_0}}(c_d),\gamma_M+d_0\delta)<-k\le\tilde l_{\Gamma_0}(a)\le \tilde l_{\Gamma_{d_0}}(a)$. Corollary~\ref{coro:3exhaust} ensures \textsc{Case} \hyperlink{C}{C} for \eqref{eq:4fconGammad}$_{d_0}$, which is a contradiction. To check that $d_\Gamma^-=\inf\{d\leq0\colon\text{\eqref{eq:4fconGammad}$_{\bar d}$ is in \textsc{Case \hyperlink{A}{A}} for all $\bar d\in[d,0]$}\}$ is bounded from below, we work with the limit of $f(t,x,\gamma)$ as $\gamma\to-\infty$.
\par
For the sake of simplicity, in what follows we will replace the subscripts $\Gamma_{d_\Gamma^+}$ and $\Gamma_{d_\Gamma^-}$ by $d_\Gamma^+$ and $d_\Gamma^-$ respectively.
Now, notice that the robustness of \textsc{Cases \hyperlink{A}{A}}, \hyperlink{C1}{\textsc{C1}} and \hyperlink{C2}{\textsc{C2}} guaranteed by Proposition~\ref{prop:3CaseCrobust} ensures that \eqref{eq:4fconGammad}$_{d_\Gamma^+}$ is in \textsc{Case \hyperlink{B1}{B1}} and
\eqref{eq:4fconGammad}$_{d_\Gamma^-}$ is in \textsc{Case \hyperlink{B2}{B2}}.
In particular, $m_{d_\Gamma^+}=l_{d_\Gamma^+}$ and $m_{d_\Gamma^-}=u_{d_\Gamma^-}$.
Let us see that \eqref{eq:4fconGammad}$_d$ is in \textsc{Case \hyperlink{C1}{C1}} for all $d>d_\Gamma^+$. We fix $d>d_\Gamma^+$ and observe that (3) guarantees $f(t,x,\Gamma_{d_\Gamma^+}(t))<f(t,x,\Gamma_d(t))$ for all $t\in[a,b]$. Since
Proposition~\ref{prop:3monotonicity}(i) ensures that $l_{d_\Gamma^+}\leq l_{\Gamma_d}$, an easy contradiction argument shows that $l_{\Gamma_d}(t)>l_{d_\Gamma^+}(t)=m_{d_\Gamma^+}(t)$ for all $t\in[a,b]$. Therefore, Theorem~\ref{th:3general}(ii) yields $\lim_{t\rightarrow\infty}(x_{d_\Gamma^+}(t,a,l_{\Gamma_d}(a))-\tilde u_{\gamma_+}(t))=0$.
A new comparison argument ensures that $x_{d_\Gamma^+}(t,a,l_{\Gamma_d}(a))\leq x_{\Gamma_d}(t,a,l_{\Gamma_d}(a))=l_{\Gamma_d}(t)$ for all $t\geq a$, and hence Lemma~\ref{lemma:4boundsandlimits}(i) for $\bar x=x_{d_\Gamma^+}(\cdot,a,l_{\Gamma_d}(a))$ ensures that $\lim_{t\rightarrow\infty}(l_{\Gamma_d}(t)-\tilde u_{\gamma_+}(t))=0$. According to Theorem~\ref{th:2massobrecasos}, $\eqref{eq:4fconGammad}_d$ is in \textsc{Case \hyperlink{C1}{C1}}. An analogous argument shows that $\eqref{eq:4fconGammad}_d$ is in \textsc{Case \hyperlink{C2}{C2}} for all $d<d_\Gamma^-$.

(ii) It is clear that all the hypotheses on $(f,\Gamma)$ are also satisfied by $(f,\Lambda)$ for all $\Lambda\in\mathcal C$, which proves the first assertion in (ii). We fix $\Lambda\in\mathcal C$ and a sequence $\{\Lambda^n\}_n$ in the same set with $\lim_{n\to\infty}\|\Lambda^n-\Lambda\|_\infty=0$. We observe, first, that the corresponding asymptotic limits also converge to those of $\Lambda$, from where it is easy to deduce that
$\lim_{n\to\infty}\|(\Lambda^n)_d-\Lambda_d\|_\infty=0$ for all $d\in\mathbb R$. Second, that by restricting ourselves to large values of $n$, we can assume that the interval $[a,b]$ and the constant $\delta$ of the proof of (i) are common for all $\Lambda^n$. And third, as a consequence of Theorem \ref{th:hypcontinuation}, that there exists $k$ such that $-k\le\tilde l_{(\Lambda^n)_0}(t)\le\tilde m_{(\Lambda^n)_0}(t)\le k$ for all $t\in\mathbb R$. These three properties show that the constant $d_0>0$ of the proof of (i) can be chosen to be an upper-bound for all $d_{\Lambda^n}^+$. Let us assume for contradiction that $\{d_{\Lambda^n}^+\}_n$ does not tend to $d_\Lambda^+$. Hence, there exits a subsequence $\{d_{\Lambda^m}^+\}_m$ with finite limit $d_*\ne d_\Lambda^+$. That is, $x'=f(t,x,\Lambda_{d_*}(t))$ is in \textsc{Case} \hyperlink{A}{A} or \hyperlink{C1}{C1}. But this contradicts the robustness of these cases ensured by Proposition \ref{prop:3CaseCrobust}, since $x'=f(t,x,(\Lambda^m)_{d_{\Lambda^m}^+}(t))$ is in \textsc{Case} \hyperlink{B1}{B1} for all $m\in\mathbb N$. An analogous argument shows that $\Gamma\mapsto d_\Gamma^-$ is continuous.

(iii) Let us take $\Lambda^1,\Lambda^2\in\mathcal C$ with $\Lambda^1\leq\Lambda^2$ and check that $d_{\Lambda^2}^+\leq d_{\Lambda^1}^+$. We take any $d\in(0,d_{\Lambda^2}^+)$ (so that $x'=f(t,x,(\Lambda^2)_d(t))$ is in \textsc{Case \hyperlink{A}{A}}), and will check that $x'=f(t,x,(\Lambda^1)_d(t))$ is also in \textsc{Case \hyperlink{A}{A}} (which means that $d\le d_{\Lambda^1}^+$ and hence proves our assertion). As seen in the proof of (i),
$\lim_{t\rightarrow\infty}(u_{(\Lambda^1)_d}(t)-\tilde u_{\gamma_+}(t))=0$, and hence
it suffices to check that $\lim_{t \rightarrow\infty}(l_{(\Lambda^1)_d}(t)-\tilde l_{\gamma_+}(t))=0$ (see Theorem~\ref{th:2massobrecasos}). It is not hard to check that $(\Lambda^1)_d\le(\Lambda^2)_d$ for all $d\ge 0$. Hence, by (3), $f(t,x,(\Lambda^1)_d(t))\le f(t,x,(\Lambda^2)_d(t))$, which according to Proposition ~\ref{prop:3monotonicity}(i) ensures that $l_{(\Lambda^1)_d}\leq l_{(\Lambda^2)_d}$. Since $\lim_{t\rightarrow\infty}(l_{(\Lambda^2)_d}(t)-\tilde l_{\gamma_+}(t))=0$,
Lemma~\ref{lemma:4boundsandlimits}(i) applied to $\bar x=l_{(\Lambda^2)_d}$ proves the claim. An analogous argument proves that $d_{\Lambda^1}^-\le d_{\Lambda^2}^-$.

Finally, let us check the strict monotonicity statement. The new condition provides $b>a$ such that $\Lambda^2(t)>\Lambda^1(t)>\max(\lim_{t\to\pm\infty}\Lambda_1(t))$ for all $t\in[a,b]$. It is not hard to check that $(\Lambda^2)_d(t)>(\Lambda^1)_d(t)$ for all $d>0$ and $t\in[a,b]$, and hence, by (3), $f(t,x,(\Lambda^1)_{d_{\Lambda^1}^+}(t))<f(t,x,(\Lambda^2)_{d_{\Lambda^1}^+}(t))$ for all $t\in[a,b]$. The argument used to complete the proof of (i) shows that $x'=f(t,x,(\Lambda^2)_{d_{\Lambda^1}^+}(t))$ is in \textsc{Case \hyperlink{C1}{C1}}, which proves that $d_{\Lambda^2}^+< d_{\Lambda^1}^+$. Analogously, $d_{\Lambda^1}^-< d_{\Lambda^2}^-$.
\end{proof}
\begin{remark} Consider the framework of Theorem~\ref{th:4Gammad} in the case of $\sup_{t\in\mathbb{R}}\Gamma(t)>\gamma_M$. Let $K_l\subset\mathbb{R}$ be a compact subset for which there exists $\delta>0$ such that $[\tilde l_{\gamma_m}(t)-\delta,\tilde m_{\gamma_m}(t)+\delta]\subseteq K_l$ for all $t\in\mathbb{R}$, and let $K_u\subset\mathbb{R}$ be a compact subset such that $[\tilde m_{\gamma_M}(t)-\delta,\tilde u_{\gamma_M}(t)+\delta]\subseteq K_u$ for all $t\in\mathbb{R}$. We recall that Proposition~\ref{prop:3monotonicity}(ii) guarantees that $[\tilde l_{\Gamma_0}(t)-\delta,\tilde m_{\Gamma_0}(t)+\delta]\subseteq K_l$ and $[\tilde m_{\Gamma_0}(t)-\delta,\tilde u_{\Gamma_0}(t)+\delta]\subseteq K_u$ for all $t\in\mathbb{R}$. This condition on $K_u$ and $K_l$ is what we need in the proof, but we define $K_l$ and $K_u$ in terms of $\tilde l_{\gamma_m}$, $\tilde m_{\gamma_m}$, $\tilde m_{\gamma_M}$ and $\tilde u_{\gamma_M}$ because they are the a priori known hyperbolic solutions of the past and future equations. A detailed look into the proof of Theorem~\ref{th:4Gammad} shows that we can still guarantee the existence of $d_\Gamma^+>0$ and the continuity and monotonic variation of $\Lambda\mapsto d_\Lambda^+$ described in (ii) and (iii) if we replace (3) of Theorem~\ref{th:4Gammad} by
\begin{enumerate}
\item[(4)] the map $\gamma\mapsto f(t,x,\gamma)$ is nondecreasing for all $(t,x)\in\mathbb{R}\times\mathbb{R}$, there exists a dense $D\subseteq\mathbb{R}$ such that $\gamma\mapsto f(t,x,\gamma)$ is strictly increasing for all $(t,x)\in D\times K_l$, and $\lim_{\gamma\rightarrow\infty} f(t,x,\gamma)=\infty$ uniformly on compact sets of $\mathbb{R}\times K_l$.
\end{enumerate}
Analogously, we can say the same for $d_\Gamma^-<0$ if we replace (3) of Theorem~\ref{th:4Gammad} by
\begin{enumerate}
\item[(5)] the map $\gamma\mapsto f(t,x,\gamma)$ is nondecreasing for all $(t,x)\in\mathbb{R}\times\mathbb{R}$, there exists a dense $D\subseteq\mathbb{R}$ such that $\gamma\mapsto f(t,x,\gamma)$ is strictly increasing for all $(t,x)\in D\times K_u$, and $\lim_{\gamma\rightarrow-\infty} f(t,x,\gamma)=-\infty$ uniformly on compact sets of $\mathbb{R}\times K_u$.
\end{enumerate}
In the case of $\inf_{t\in\mathbb{R}}\Gamma(t)<\gamma_m$: condition (5) ensures the existence, continuity and monotonicity of $d_\Gamma^+>0$ and (4) ensures that of $d_\Gamma^-<0$.
\end{remark}
Theorem~\ref{th:4Gammad} allows us to define continuous bifurcation maps for the $c$-parametric rate-induced and phase-induced tipping problems described at the beginning of Section \ref{sec:mechanisms}, as the next result shows: continuous functions of $c$ whose signs determine the particular dynamical case.
\begin{theorem}\label{th:4bifurcationfunction} Let $(f,\Gamma)$ satisfy {\rm\hyperlink{h1}{\textbf{h1}}-\hyperlink{h6}{\textbf{h6}}} and the conditions {\rm(1)-(3)} of Theorem~{\rm\ref{th:4Gammad}} in the situation $\sup_{t\in\mathbb{R}}\Gamma(t)>\gamma_M$. Let $\Gamma^c$ be defined by either $t\mapsto \Gamma(ct)$ for $c\in C=(0,\infty)$ or by $t\mapsto \Gamma(t+c)$ for $c\in C=\mathbb R$. Then,
\begin{itemize}
\item[\rm(i)] there exists a continuous bifurcation function $\varphi\colon C\rightarrow\mathbb{R}$ such that equation
\begin{equation}\label{eq:4fconGammac}
 x'=f(t,x,\Gamma^c(t))
\end{equation}
is in \textsc{Case \hyperlink{A}{A}} if $\varphi(c)<0$, \textsc{\hyperlink{B1}{B1}} if $\varphi(c)=0$, and \textsc{\hyperlink{C1}{C1}} if $\varphi(c)>0$.
\item[\rm(ii)] In the case of $\Gamma^c(t)=\Gamma(ct)$, if $t\mapsto \Gamma(t)$ is nondecreasing on $(-\infty,0]$ and nonincreasing on $[0,\infty)$, then $\varphi$ is strictly decreasing on~$(0,\infty)$.
\item[\rm(iii)] The same properties hold without the hypothesis about $\lim_{\gamma\rightarrow\pm\infty}f(t,x,\gamma)$ in condition {\rm (3)} of Theorem {\rm\ref{th:4Gammad}}.
\end{itemize}
In the situation $\inf_{t\in\mathbb{R}}\Gamma(t)<\gamma_m$: {\rm(i)} holds with {\rm \hyperlink{B1}{B1}} and {\rm \hyperlink{C1}{C1}} respectively replaced by {\rm\hyperlink{B2}{B2}} and {\rm \hyperlink{C2}{C2}}; the bifurcation function is also strictly decreasing if $t\mapsto \Gamma(t)$ is nonincreasing on $(-\infty,0]$ and nondecreasing on $[0,\infty)$; and {\rm(iii)} also holds.
\end{theorem}
\begin{proof}
(i) Let us fix $c\in C$ and define $(\Gamma^c)_d$ for $d\ge 0$ as in the statement of  Theorem~\ref{th:4Gammad}, which provides $d_{\Gamma^c}^+>0$ such that
$x'=f(t,x,(\Gamma^c)_d(t))$ is in \textsc{Case \hyperlink{A}{A}} if $d\in[0,d_{\Gamma^c}^+)$, in \textsc{\hyperlink{B1}{B1}} if $d=d_{\Gamma^c}^+$ and in \textsc{\hyperlink{C1}{C1}} if $d\in(d_{\Gamma^c}^+,\infty)$. This means that \eqref{eq:4fconGammac}$_c$ is in \textsc{Case \hyperlink{A}{A}} if $1\in[0,d_{\Gamma^c}^+)$, in \textsc{\hyperlink{B1}{B1}} if $d_{\Gamma^c}^+=1$, and in \textsc{\hyperlink{C1}{C1}} if $1\in(d_{\Gamma^c}^+,\infty)$.
Therefore, the map $\varphi\colon C\to\mathbb R$ given by $\varphi(c)=1-d_{\Gamma^c}^+$ satisfies the statements concerning the dynamical cases.
Notice that $\Gamma^c$ belongs to the set $\mathcal C$ of Theorem~\ref{th:4Gammad}(ii) for all $c\in C$. The continuity of $C\rightarrow \mathcal C$, $c\mapsto\Gamma^c$ in the $\|\cdot\|_\infty$ norm (which, in both cases, follows from the existence of asymptotic limits of $\Gamma$ and its uniform continuity) and Theorem~\ref{th:4Gammad}(ii) ensure the continuity of $\varphi$ on $C$.

(ii) It is easy to check that $t\mapsto \Gamma(t)$ is nondecreasing on $(-\infty,0]$ and nonincreasing on $[0,\infty)$ if and only if $c\mapsto\Gamma^c(t)= \Gamma(ct)$ is nonincreasing for all $t\in\mathbb{R}$. Notice that $\Gamma(0)=\sup_{t\in\mathbb{R}}\Gamma(t)>\gamma_M\geq \gamma_+=\lim_{t\rightarrow\infty}\Gamma(t)$. We take $\gamma_0\in(\gamma_M,\Gamma(0))$ and
$t_0=\inf\{t>0\colon \Gamma(t)=\gamma_0\}$, so that $\Gamma(t_0)<\Gamma(s)$ for every $0\leq s<t_0$. Given $0<c_1<c_2$, we take $t_*=t_0/c_2$. Since $0<c_1t_*<c_2t_*=t_0$, we get $\gamma_M<\Gamma^{c_2}(t_*)<\Gamma^{c_1}(t_*)$. Theorem~\ref{th:4Gammad}(iii) guarantees that $\varphi(c_2)=1-d_{\Gamma^{c_2}}^+<1-d_{\Gamma^{c_1}}^+=\varphi(c_1)$.

(iii) Let us define $[\bar\gamma_m,\bar\gamma_M]=\mathrm{cl}(\Gamma(\mathbb{R}))$ and
\begin{equation*}
g(t,x,\gamma)=\left\{\begin{array}{ll}
f(t,x,\bar\gamma_m)-(\gamma-\bar\gamma_m)^2&\text{if }\gamma<\bar\gamma_m\,,\\
f(t,x,\gamma)&\text{if }\gamma\in [\bar\gamma_m,\bar\gamma_M]\,,\\
f(t,x,\bar\gamma_M)+(\gamma-\bar\gamma_M)^2&\text{if }\gamma>\bar\gamma_M\,.\\
\end{array}\right.
\end{equation*}
Then $g(t,x,\Gamma^c(t))=f(t,x,\Gamma^c(t))$, and it is easy to check that $(g,\Gamma)$ satisfies all the initial hypotheses, from where (iii) follows.

The final statements can be proved as the previous ones from the corresponding assertions in Theorem \ref{th:4Gammad}.
\end{proof}

Note that Theorem~\ref{th:4bifurcationfunction}(ii) does not refer to the phase-induced tipping problem since the monotonicity of $c\mapsto\Gamma(t+c)$ on $\mathbb{R}$ for all $t\in\mathbb{R}$ is equivalent to the monotonicity of $t\mapsto\Gamma(t)$ on $\mathbb{R}$, in which case Theorem \ref{th:4compacity} precludes the possibility of tipping. Recall that the analogues of Theorems~\ref{th:4Gammad} and \ref{th:4bifurcationfunction} can be formulated in the case of nonincreasing $\gamma\mapsto f(t,x,\gamma)$: see Remark \ref{rm:nonincreasingchange}.

To end this section, we use Theorem~\ref{th:4Gammad} to study the size-induced, rate-induced and phase-induced tipping problems $x'=h(t,x-d\, \Gamma(t))$, $x'=h(t,x- \Gamma(ct))$ and $x'=h(t,x-\Gamma(t+c))$ for functions $h$ of two variables and nonconstant monotonic transition
functions $\Gamma$. We will prove that there exist $d_1<0<d_2$ such that the equations present tracking when the parameter belongs to $(d_1,d_2)$ and tipping when it is outside its closure.
Hence, $d_1$ and $d_2$ are nonautonomous bifurcation points, which will be called \emph{tipping points}.

Before presenting the result, we introduce suitable spaces of transition functions:
\begin{equation*}
\begin{split}
\mathfrak{D}_+&=\{\Gamma\in C^1(\mathbb R,\mathbb R)\colon\text{$\Gamma$ is nondecreasing and nonconstant,}\\
&\qquad\text{and there exist $\lim_{t\rightarrow\pm\infty}\Gamma(t)\in\mathbb R$ and
$\lim_{t\rightarrow\pm\infty}\Gamma'(t)=0$}\},\;\\
\mathfrak{D}_-&=\{\Gamma\in C^1(\mathbb R,\mathbb R)\colon-\Gamma\in\mathfrak D_+\}\,.
\end{split}
\end{equation*}
\begin{theorem}\label{th:applshift} Let $h$ satisfy {\rm\hyperlink{h2t}{\textbf{h2}$_*$}-\hyperlink{h6t}{\textbf{h6}$_*$}}, and let $\Gamma \in \mathfrak{D}_{+}$ (resp. $\Gamma\in\mathfrak{D}_{-}$).
Then,
\begin{itemize}
\item[\rm(i)] there exist $d_\Gamma^-<0<d_\Gamma^+$ such that
\begin{equation}\label{eq:4shift}
x'=h(t,x-d\, \Gamma(t))
\end{equation}
is in \textsc{Case \hyperlink{A}{A}} for $d\in (d_\Gamma^-,d_\Gamma^+)$, in \textsc{Case \hyperlink{B1}{B1}} (resp. \textsc{\hyperlink{B2}{B2}}) for $d=d_\Gamma^-$, in \textsc{Case \hyperlink{B2}{B2}} (resp. \textsc{\hyperlink{B1}{B1}}) for $d=d_\Gamma^+$, in \textsc{Case \hyperlink{C1}{C1}} (resp. \textsc{\hyperlink{C2}{C2}}) for $d\in(-\infty,d_\Gamma^-)$, and in \textsc{Case \hyperlink{C2}{C2}} (resp. \textsc{\hyperlink{C1}{C1}}) for $d\in(d_\Gamma^+,\infty)$.
\item[\rm(ii)] Let $\mathfrak{D}_\pm$ be endowed with the $\|\cdot\|_\infty$ norm. The maps $\Gamma\mapsto d_\Gamma^-,d_\Gamma^+$ are continuous on $\mathfrak{D}_\pm$.
\end{itemize}
\end{theorem}
\begin{proof} (i) The change of variables $x=y+d\,\Gamma(t)$ takes \eqref{eq:4shift} to
\begin{equation}\label{eq:cambio}
y'=h(t,y)-d\, \Gamma'(t)\,.
\end{equation}
We will reason in the case $\Gamma \in \mathfrak{D}_{+}$: the other one is analogous. Let us define $f(t,y,\gamma)=h(t,y)+\gamma$ and check that the pair $(f,-\Gamma')$ satisfies all the hypotheses of Theorem~\ref{th:4Gammad}. Since $\lim_{t\rightarrow\pm\infty} \Gamma'(t)=0$, the constants $\gamma_m$ and $\gamma_M$ defined by \eqref{def:gama+-} for $\Gamma'$ are 0. It follows easily that $(f,-\Gamma')$ satisfies {\rm\hyperlink{h1}{\textbf{h1}}-\hyperlink{h6}{\textbf{h6}}}. On the other hand, since $\Gamma'$ is not identically 0 (because $\lim_{t\to-\infty}\Gamma(t)<\lim_{t\to\infty}\Gamma(t)$), condition (1) holds, with
$\inf_{t\in\mathbb R}(-\Gamma')(t)<0$. Condition \hyperlink{h6}{\textbf{h6}$_0$} ensures (2), and (3) is obvious. Theorem~\ref{th:4Gammad} applied to this pair ensures the existence of $d_\Gamma^-<0<d_\Gamma^+$ such that \eqref{eq:cambio} is in \textsc{Case \hyperlink{A}{A}} for $d\in(d^-_\Gamma,d^+_\Gamma)$, \textsc{\hyperlink{B1}{B1}} for $d=d_\Gamma^-$, \textsc{\hyperlink{B2}{B2}}, for $d=d_\Gamma^+$, \textsc{\hyperlink{C1}{C1}} for $d\in(-\infty,d_\Gamma^-)$, and \textsc{ \hyperlink{C2}{C2}} for $d\in(d_\Gamma^+,\infty)$.
Having in mind how the change of variables $x=y+d\,\Gamma(t)$ transforms the set of bounded solutions, we conclude that \eqref{eq:4shift}$_d$ presents the same cases.

(ii) The arguments of Theorem~\ref{th:4Gammad}(ii) show the continuity of $\Gamma\mapsto d_\Gamma^-,d_\Gamma^+$.
\end{proof}
\begin{theorem}\label{th:4bifurcationfunction2} Let $h$ satisfy {\rm\hyperlink{h2t}{\textbf{h2}$_*$}-\hyperlink{h6t}{\textbf{h6}$_*$}}, and let $\Gamma \in \mathfrak{D}_{+}$ (resp. $\Gamma\in\mathfrak{D}_{-}$). Let $\Gamma^c$ be defined by either $t\mapsto \Gamma(ct)$ for $c\in C=(0,\infty)$ or by $t\mapsto \Gamma(t+c)$ for $c\in C=\mathbb R$. Then, there exists a continuous bifurcation function $\varphi\colon C\rightarrow\mathbb{R}$ such that
\begin{equation}\label{eq:ratephase}
x'=h(t,x-\Gamma^c(t))
\end{equation}
is in \textsc{Case \hyperlink{A}{A}} if $\varphi(c)<0$, in \textsc{Case \hyperlink{B2}{B2}} (resp.~{\rm\hyperlink{B1}{B1}}) if $\varphi(c)=0$, and in \textsc{Case \hyperlink{C2}{C2}} (resp.~{\rm\hyperlink{C1}{C1}}) if $\varphi(c)>0$.
\end{theorem}
\begin{proof} As in the proof of Theorem~\ref{th:applshift}, we make the change of variables $y=x-\Gamma^c(t)$, which takes \eqref{eq:ratephase} to
\[
 y=h(t,y)-(\Gamma^c)'(t)\,,
\]
define $f(t,y,\gamma)=h(t,y)+\gamma$, and check that $(f,-(\Gamma^c)')$ satisfies the conditions of the second situation analyzed in Theorem~\ref{th:4Gammad}. Hence, we can repeat the arguments proving Theorem \ref{th:4bifurcationfunction}(i) to find the continuous bifurcation function $\varphi\colon C\to\mathbb R$ for the transformed equation. And, as in the proof of Theorem \ref{th:applshift}, it is clear that the change of variables preserves the dynamical situation.
\end{proof}
\section{A single species population subject to Allee effect}\label{section:singlespeciesAlleeeffect}
The Allee effect (see \cite{courchamp1}, \cite{kramer1} and the references therein)
consists on a positive correlation between the size of a population and its fitness, that is, the population growth rate per individual.
Since the Allee effect has been found to be responsible for an increase in the risk of extinction for low density populations,
it has attracted the interest of both biologists and mathematicians.
There exist several biological mechanisms related to survival and reproduction which can justify the appearance of the Allee effect on different biological systems (see \cite{berec1}, \cite{courchamp1}). The study of such mechanisms is an active area of both theoretical and experimental research.

In this section, we will use the results of \cite{dno1,dno2} and of the previous sections to explain how scalar nonautonomous equations with concave derivative are a natural way to model populations subject to Allee effect, and how critical transitions in nonautonomous population dynamics can be triggered by Allee effect; that is, how abrupt and sudden changes in the number of individuals can arise from small changes on external parameters of these models.
In particular, concerning critical transitions, we will discuss the extinction of a native species \cite{murray1} and the invasion of a habitat \cite{lewis1}, that is, the uncontrolled proliferation of a non-native species.

In the classical autonomous logistic population equation, $x'=r\,x\,(1-x/K)$, $r$ and $K$ are strictly positive constants representing the intrinsic rate of increase with unlimited resources and the maximum population size with positive growth rate, respectively. Our starting point will be the nonautonomous counterpart of this equation, $x'=r(t)\,x\,(1-x/K(t))$, where $r$ and $K$ are bounded uniformly continuous functions positively bounded from below. The role of $r$ in this time-dependent setting is similar, but the time-dependent map $K$ does no longer represents the maximum population size: its role is replaced by an a priori unknown strictly positive solution which is hyperbolic and attractive (see e.g.~\cite{lno1}).
This model appears, for instance, in Chapter 8 of \cite{renshaw1}. Allee effect can be added to the nonautonomous logistic model either in a general multiplicative form (\cite{boukal1}, \cite{ye1}),
\begin{equation}\label{eq:alleemultiplicativo}
x'=r(t)\, x\, \left(1-\frac{x}{K(t)}\right)\frac{x-S(t)}{K(t)}\,,
\end{equation}
where the bounded uniformly continuous function $S$, which satisfies $S(t)+K(t)\ge 0$ for all $t\in\mathbb{R}$, determines the strength of Allee effect;
or in an additive form \cite{liu1}, adding a Holling type II functional response term, usually used to model Allee effect mechanisms related to predation,
\begin{equation}\label{eq:alleeaditivo}
x'=r(t)\, x\,\left(1-\frac{x}{K(t)}\right)-\frac{a(t)x}{x+b(t)}\,.
\end{equation}
Here, $a$ and $b$ are also bounded uniformly continuous functions positively bounded from below, which include the contribution to Allee effect of the number of attacks and the average time spent by predators on processing a food item.
The  Holling type III functional response term $-a(t)\,x^2/(x^2+b(t))$ can also be added to the cubic right-hand side of \eqref{eq:alleemultiplicativo} to describe other predation features.
 Note that the right-hand sides of \eqref{eq:alleemultiplicativo} and \eqref{eq:alleeaditivo} are coercive functions with strictly concave derivative, and that $t\mapsto K(t)$ and $t\mapsto S(t)$ are not solutions of \eqref{eq:alleemultiplicativo} unless they are constant.
We will only classify the types of Allee effect under the assumption of the existence of three hyperbolic solutions of the corresponding equations, so that it is important to highlight that this property is not guaranteed a priori for \eqref{eq:alleemultiplicativo} or \eqref{eq:alleeaditivo}.
Notice also that these three hyperbolic solutions are not necessarily positive.
We recall that Theorem~\ref{th:3three} and Remark~\ref{rm:3siempre} show that a coercive strictly d-concave equation can have
at most three uniformly separated hyperbolic solutions;
if there are three, the upper one and the lower one are attractive and delimit the set of bounded solutions, and hence they can be understood as steady population states
if they are nonnegative; and if there is just one, then it is attractive.

The Allee effect has been usually classified in two categories, {\em strong\/} or {\em weak}, depending on the existence or not of a critical population size (see e.g. \cite{boukal1}, \cite{courchamp1}, \cite{sun1}): in the autonomous case, this critical population size is given by a strictly positive repulsive hyperbolic equilibrium.
In the next paragraphs, we propose an approach to these concepts in the nonautonomous case, always under the hypotheses of existence of three hyperbolic solutions, and depending on the sign of these solutions.
This approach intends to be valid also for models which are more general than \eqref{eq:alleemultiplicativo} and \eqref{eq:alleeaditivo}, for which 0 may not solve the equation.
In this case, if the lower attractive hyperbolic solution is nonnegative and close to 0, it will not represent extinction, but a sparse (or low density) steady population.
In all the cases (either extinction, sparse steady population, or absence of both), the smallest nonnegative hyperbolic solution delimits from below the positively invariant region on which we will analyze the dynamics.

So, let us consider a population modelled by a coercive and strictly d-concave equation $x'=h(t,x)$ with three hyperbolic solutions. We shall say that it presents {\em strong Allee effect\/} if these three solutions are nonnegative, in which case the extinct or sparse population represented by the lower one is attractive. In this case, the repulsive hyperbolic solution (the middle one) plays the role of critical population size (which is time-dependent in our autonomous case): the population growth rate per individual $h(t,x)/x$ has negative average through any forward semitrajectory which is positive and below this critical solution (if the smallest steady state is small enough). We shall say that the population presents {\em weak Allee effect\/} if the equation has exactly two hyperbolic solutions which are nonnegative, in which case the lower one is repulsive. This means a positive average population growth rate per individual through the lower nonnegative hyperbolic solution (if this one is small enough), and corresponds to the autonomous idea of having weak Allee effect when the population growth rate at 0 is positive but smaller than $r$. We will not study the intermediate cases between strong and weak Allee effect in this work.

Let us analyze these concepts for the models \eqref{eq:alleemultiplicativo} and \eqref{eq:alleeaditivo}, for which 0 is always a solution, assuming the existence of three hyperbolic ones.
If \eqref{eq:alleemultiplicativo} exhibits strong Allee effect, then 0 (which is the smallest nonnegative bounded solution and hence hyperbolic) is attractive, which according to  \eqref{eq:1hurwitzinfty} is equivalent to
\begin{equation}\label{eq:5alleeffectstrengthindicator1}
\lim_{l\rightarrow\infty}\Big(\sup\Big\{\frac{1}{t-s}\int_s^t -\frac{r(\tau)S(\tau)}{K(\tau)}\, d\tau\colon\; t-s\geq l\Big\}\Big)<0\,.
\end{equation}
Conversely, \eqref{eq:5alleeffectstrengthindicator1} ensures that 0 is hyperbolic attractive, and hence that it is the smallest one since $h_{xx}(t,0)=2r(t)(K(t)+S(t))/K^2(t)\ge 0$ (see Proposition 6.5 of \cite{dno1}, and note that the minimality of the hull is not required to prove that there are no bounded solutions separated from 0 below it).
Therefore, \eqref{eq:alleemultiplicativo} exhibits strong Allee effect.
Equation \eqref{eq:alleemultiplicativo} exhibits weak Allee effect if and only if 0 is hyperbolic repulsive, which according to \eqref{eq:1hurwitz-infty} is equivalent to
\begin{equation}\label{eq:5alleeffectstrengthindicator2}
\lim_{l\rightarrow\infty}\Big(\inf\Big\{\frac{1}{t-s}\int_s^t -\frac{r(\tau)S(\tau)}{K(\tau)}\, d\tau\colon\, t-s\geq l\Big\}\Big)>0\,.
\end{equation}
In this way, once fixed $r$ and $K$, the function $S$ determines the type of Allee effect that the population exhibits, and \eqref{eq:5alleeffectstrengthindicator1} and \eqref{eq:5alleeffectstrengthindicator2} are indicators of the strength of the strong or weak  Allee effect. For instance, if $\inf_{t\in\mathbb{R}}S(t)>0$, then \eqref{eq:5alleeffectstrengthindicator1} holds; and if $\sup_{t\in\mathbb{R}}S(t)<0$, then \eqref{eq:5alleeffectstrengthindicator2} holds.

In turn, \eqref{eq:alleeaditivo} exhibits weak Allee effect if and only if  $\lim_{l\rightarrow\infty}\big(\inf\{(1/(t-s))$ $\int_s^t (r(\tau)-a(\tau)/b(\tau))\, d\tau\colon\, t-s\geq l\}\big)>0$; and strong Allee effect if and only if  $\lim_{l\rightarrow\infty}\big(\sup\{ (1/(t-s))\int_s^t(r(\tau)-a(\tau)/b(\tau))\, d\tau\colon\, t-s\geq l\}\big)<0$ and it admits a strictly positive hyperbolic solution (which is not always the case, but holds if, for instance, $h_{xx}(t,0)=2a(t)/b^2(t)-2r(t)/K(t)\ge 0$ for all $t\in\mathbb R$).

We underline that, unlike in autonomous dynamics, the difference between strong and weak Allee effects is given by averages of the involved functions.

Notice finally that we are talking of Allee effect in the case of populations modelled by dissipative and bistable equations. This bistability naturally appears in some equations $x'=h(t,x)$ with d-concave $h$, and also in some equations for which the population growth rate per individual $h(t,x)/x$ is concave (that is, $x$-concave equations; see \cite{dno1,dno2} and Section~5 of \cite{nunezobaya}). And notice also that equations for which $h$ is d-concave and with only one hyperbolic solution can represent extinction scenarios.
\subsection{Critical transitions induced by migration or harvesting}
In this subsection, we will add a migration term to the multiplicative Allee effect equation
\eqref{eq:alleemultiplicativo} with fixed quasiperiodic coefficients (to be specified below). This term is $\gamma\phi(t)$: $\gamma$ is a nonnegative parameter and $\phi$ is a quasiperiodic function positively bounded from below,
which represents the arrival of new individuals to the habitat:
\begin{equation}\label{eq:differentbifurcation1}
 x'=r(t)\, x\,\left(1-\frac{x}{K(t)}\right)\frac{x-S(t)}{K(t)}+\gamma\,\phi(t)\,.
\end{equation}
As in Section \ref{sec:mechanisms}, we will consider transition equations
\begin{equation}\label{eq:transitionmigration}
 x'=r(t)\, x\,\left(1-\frac{x}{K(t)}\right)\frac{x-S(t)}{K(t)}+\Gamma(ct)\,\phi(t)
\end{equation}
for an impulse function $\Gamma$ and a positive rate $c$. Our goal is to show that different choices of $\Gamma$ cause different types of rate-induced critical transitions: the Allee effect will be strong for the equation \eqref{eq:differentbifurcation1} corresponding to the (equal) asymptotic limits of $\Gamma$, and it will give rise to the same dynamical situation for large rates; but, as $c$ decreases, it will appear very different situations: from population close to extinction to invasion of the habitat.
\par
Departures of individuals (emigration), harvesting and hunting could be represented by adding an analogous negative parametric term. For the sake of simplicity, we will not deal with this case. \par
The function $\Gamma$ is given by $\Gamma(t)=\gamma_++(\gamma_*-\gamma_+)\exp(-t^2/10)$: it represents an impulse from its asymptotic limits at $\pm\infty$, given by $\gamma_+$, to a value $\gamma_*\neq\gamma_+$.
Clearly, a larger $c$ causes a faster transition. The remaining coefficients will be $r(t)\equiv 1$, $K(t)=30+60\sin^2(t)$, $S(t)=40+40\cos^2(t\sqrt{5}/2)$ and $\phi(t)=0.8+0.4\sin^2(t\sqrt{5}/2)$. These functions define the right-hand side  of \eqref{eq:differentbifurcation1}, which we rewrite as $f(t,x,\gamma)$. It is not hard to check that the proofs of Theorems~5.5 and 5.10 of \cite{dno1} can be repeated for \eqref{eq:differentbifurcation1} due to the condition $\inf_{t\in\mathbb R}\phi(t)>0$.
(The extension of these results require to construct the hull of $(r,K,S,\phi)$, on which $f(t,x,\gamma)$ can be represented: see Subsection \ref{subsec:2skewproduct}, and observe that all the hypotheses of both theorems are fulfilled.) These results ensure that: the set $\mathcal{R}_f$ given by \eqref{def.R} is an open interval; and that the lower bounded solution of \eqref{eq:differentbifurcation1} strictly increases as $\gamma$ increases. A simple numeric simulation shows that 0 is the unique bounded solution for $\gamma=0$.
This uniqueness, which is not possible in the autonomous formulation of \eqref{eq:differentbifurcation1}, is fundamental in what follows.
In particular, all the bounded solutions are strictly positive for all the positive values of $\gamma$.
Altogether, we can identify $\mathcal{R}_f$ with the set of values of $\gamma$ such that \eqref{eq:differentbifurcation1}$_\gamma$ exhibits strong Allee effect.

We will always take $\gamma_+\in\mathcal R_f$. (In particular, the pair $(f,\Gamma)$ satisfies \hyperlink{h1}{\textbf{h1}}-\hyperlink{h6}{\textbf{h6}}, with $\gamma_-=\gamma_+$: see Remark~\ref{rm:twovariables}). If we choose a extreme value $\gamma_*$ of $\Gamma$ also in $\mathcal R_f$, then Theorem~\ref{th:4compacity} precludes the existence of tipping for any $c\in(0,\infty)$. So, we will choose $\gamma_*\notin\mathcal R_f$.
This implies that the transition takes values outside $\mathcal R_f$ during a period of time which is determined by the rate $c$. What we will observe is that, if this period is short (i.e., if $c$ is large), then the dynamics of the transition equation is basically equal to that of the future equation (which is equal to the past equation). But if the period is long enough (given by a sufficiently small $c$), then the dynamics changes dramatically, in two possible different ways. That is, there is at least one (and, as we will see, just one) positive critical value of the rate.

\begin{figure}[h]
\definecolor{verdebosque}{rgb}{0,0.39,0}
\definecolor{rojolava}{rgb}{0.8,0.16,0.16}
\centering
\begin{subfigure}[b]{0.48\textwidth}
         \centering
         \includegraphics[width=\textwidth]{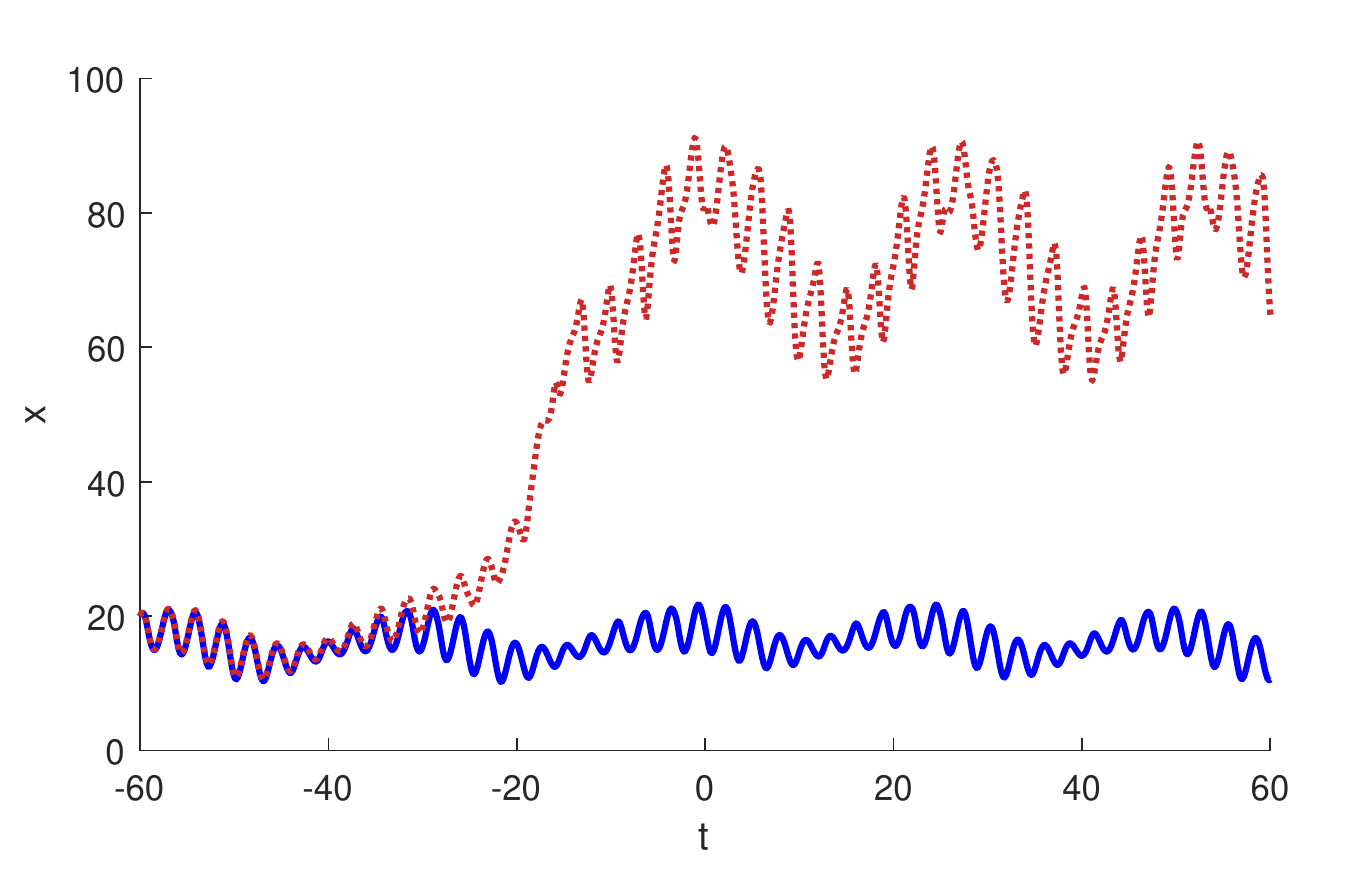}
         \caption{For $\gamma_+<\gamma_*$, small values of the transition rate can cause habitat invasion (in dotted red line). A controlled small population persists for large rate (in solid blue).}
         \label{fig:invasion}
     \end{subfigure}
\hfill
\begin{subfigure}[b]{0.48\textwidth}
         \centering
         \includegraphics[width=\textwidth]{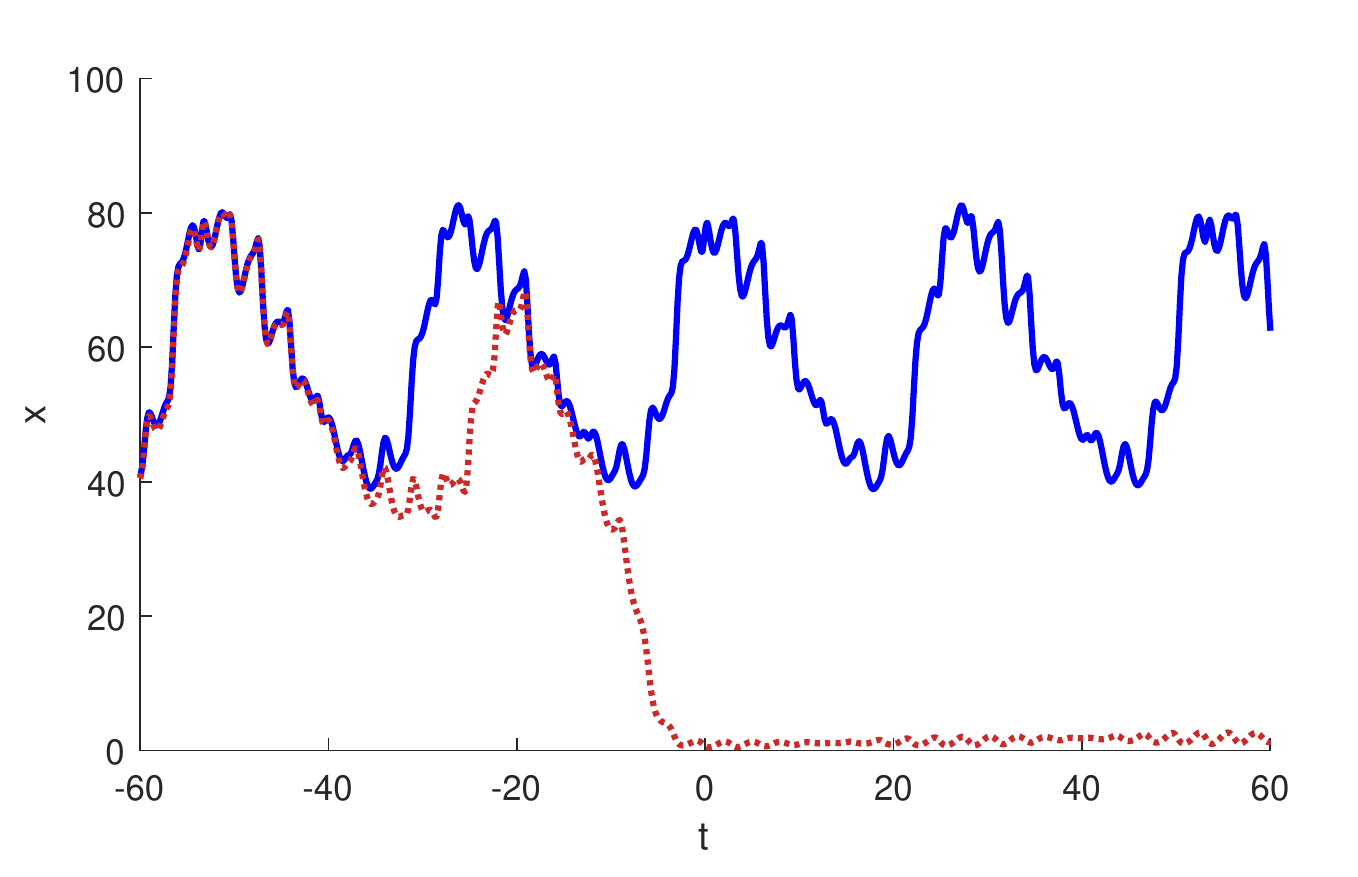}
         \caption{For $\gamma_+>\gamma_*$, small values of the transition rate can cause habitat extinction (in dotted red line). A healthy large population persists for large rate (in solid blue).}
         \label{fig:extinction}
     \end{subfigure}
\caption{Critical transitions induced by migration.}
\label{fig:invasionextinction}
\end{figure}

The robustness of hyperbolicity makes it easy to obtain numerical evidences of: $[2.0,8.5]\subset\mathcal R_f$, $1.0\notin\mathcal R_f$, and $9.0\notin\mathcal R_f$.
Figure~\ref{fig:invasion} corresponds to $\gamma_+=8.5$ and $\gamma_*=9.0$. The dotted red line represents the lowest bounded solution for $c=0.1$, and the solid blue one for $c=1.0$. In both cases (as for any $c>0$), this lowest solution may correspond to a desirable low density population which is under control for large negative values of $t$, when $\Gamma(ct)$ is practically equal to $\gamma_+$. The value of $\Gamma(ct)$ is also practically equal to $\gamma_+$ if $t$ is large enough, but it smoothly increases until $\gamma_*$ as $|t|$ approaches 0. This evolution occurs during a period of time which increases as the rate $c$ decreases. For the value $c=1.0$ (as for any large enough $c$), the previous low population remains under control for always; but, for $c=0.1$ (as for any small enough $c$), the smallest steady population undergoes an overgrowth which leads to the invasion of the habitat. So, there exists a critical rate: a threshold which must be exceeded in order to avoid an invasion. A rate above this threshold means a greater immigration for a period of time which is short enough to allow the population to keep its controlled size. But if the time is longer, invasion occurs.

Figure~\ref{fig:extinction} corresponds to $\gamma_+=2.0$ and $\gamma_*=1.0$. Now, the immigration smoothly decreases from $\gamma_+=2.0$ to $\gamma_*=1.0$ when $|t|$ approaches 0.
Here, the dotted red line represents the upper bounded solution for $c=0.1$, and the solid blue one for $c=1.0$. Let us understand this upper solution as a healthy population. As before, the behaviour does not depend on $c$ if $-t$ is large enough. As we observe in the figure, this population persists if the rate is large enough, but the population gets basically extinct
if $c$ is very small. So, again, there exists a critical rate: a threshold which must be exceeded to avoid the critical extinction. A rate above this threshold means a lower immigration for a period of time short enough to avoid extinction.

It is remarkable that these two parametric problems present an opposite behavior to other ones more usual in the literature (\cite{aspw}, \cite{lno1}), in which tracking takes place at low transition rates and tipping appears for high transition rates.
For this reason, we talk of cases of {\em rate-induced tracking} in these two scenarios.

The difference between the type of critical transitions in the two analyzed examples can be easily explained by Theorem~\ref{th:4bifurcationfunction}. Notice that we are in the cases described in this theorem both in the cases $\gamma_*>\gamma_+$ (as in Figure ~\ref{fig:invasion}) and $\gamma_*<\gamma_+$ (as in Figure \ref{fig:extinction}), but the relation between both coefficients yields a significative difference: Theorem~\ref{th:4bifurcationfunction} provides a strictly decreasing continuous bifurcation function $\varphi\colon (0,\infty)\rightarrow\mathbb{R}$, with at most a zero. If this zero $c_0$ exists and $\gamma_*>\gamma_+$ (resp. $\gamma_*<\gamma_+$), then \eqref{eq:transitionmigration}$_c$ is in
\textsc{Case \hyperlink{C1}{C1}} (resp. \textsc{\hyperlink{C2}{C2}}) if $0<c<c_0$,
\textsc{\hyperlink{B1}{B1}} (resp. \textsc{\hyperlink{B2}{B2}}) if $c=c_0)$, and
\textsc{\hyperlink{A}{A}} if $c>c_0$.

We point out once again that the occurrence of the (unique) critical transition is due to the fact that $\Gamma(ct)\notin\mathcal R_f$ for a period of time which tends to infinity as $c$ tends to cero. The radical difference appearing in the two cases described in Figures \ref{fig:invasion} and \ref{fig:extinction} depends on the relation $\inf\Gamma(ct)<\inf\mathcal R_f$ (extinction) or $\sup\Gamma(ct)>\sup\mathcal R_f$ (invasion). To this regard, it is interesting to analyze the connection with the equivalent of Theorem~5.10 of \cite{dno1} for the parametric family \eqref{eq:differentbifurcation1}: $\gamma\in\mathcal R_f$ means three hyperbolic solutions (\textsc{Case \hyperlink{A}{A}}); when $\gamma$ decreases and crosses $\inf\mathcal R_f$, just a hyperbolic solution exists, which is the hyperbolic continuation of the lower one (i.e., we arrive to \textsc{Case \hyperlink{C2}{C2}} crossing \textsc{Case \hyperlink{B2}{B2}}); and, when $\gamma$ increases and crosses $\sup\mathcal R_f$, just a hyperbolic solution exists, which is the hyperbolic continuation of the upper one (i.e., we arrive to \textsc{Case \hyperlink{C1}{C1}} crossing \textsc{Case \hyperlink{B1}{B1}}).
\subsection{Persistence or extinction due to predation depending on the type of Allee effect}
In this subsection, we will construct two examples of equations \eqref{eq:alleemultiplicativo}
with quasiperiodic right-hand side and with strong and weak multiplicative Allee effect, add a parametric term representing a Holling type III functional response to predation,
and show that a critical extinction occurs in the first case and not in the second one.
We fix $\beta>0$ and consider the $\gamma$-parametric model
\begin{equation}\label{eq:differentbifurcation}
x'=r(t)\, x\,\left(1-\frac{x}{K(t)}\right)\frac{x-S(t)}{K(t)}-\gamma\frac{x^2}{\beta+x^2}\,,
\end{equation}
which we rewrite as $x'=f(t,x,\gamma)$. It is clear that $f$ satisfies \hyperlink{h2}{\bf h2} and \hyperlink{h3}{\bf h3}.
In order to get also \hyperlink{h4}{\bf h4} and \hyperlink{h5}{\bf h5$_\gamma$}, we restrict to
\begin{equation*}
\gamma\in I_\beta=\left[0,\,\beta^{3/2}\frac{64}{5\sqrt{(5-2\sqrt{5})}\,
(7+3\sqrt{5})}\,\inf_{t\in\mathbb{R}}\frac{r(t)}{K(t)^2}\right)\,:
\end{equation*}
for these values of $\gamma$, $\sup_{(t,x)\in\mathbb R\times \mathbb R}f_{xxx}(t,x,\gamma)<0$.
As in Section \ref{sec:mechanisms}, we will substitute the parameter $\gamma$ by the transition function $\Gamma_d(t)=d\, \Gamma(t)$, with $\Gamma(t)=1/2+\arctan(t)/\pi$. We will see that this may (not necessarily) give rise to a size-induced critical transition in the dynamics of
\begin{equation}\label{eq:transitionplascubic}
 x'=r(t)\, x\,\left(1-\frac{x}{K(t)}\right)\frac{x-S(t)}{K(t)}-d\,\Gamma(t)\frac{x^2}{\beta+x^2}
\end{equation}
as the parameter $d$ increases (i.e., as predation increases). Since the left and right asymptotic limits of $\Gamma_d$ are $0$ and $d$, if we take $d$ in $I_\beta$, then the pair
$(f,\Gamma_d)$ satisfies \hyperlink{h1}{\bf h1}-\hyperlink{h5}{\bf h5}  (see again Remark \ref{rm:twovariables}\hyperlink{311}{.1}).

We will choose $r$, $K$, $S$ and $\beta$ in such a way that the past equation $x'=f(t,x,0)$ has three hyperbolic solutions. That is, 0 belongs to the set $\mathcal R_f$ defined by \eqref{def.R}. Hence, the robustness of the property of existence of hyperbolic solutions given by Proposition \ref{prop:2hiperbolicas}, combined with that of \textsc{Case \hyperlink{A}{A}} given by Proposition~\ref{prop:3CaseCrobust}, ensures that \eqref{eq:transitionplascubic}$_d$ is in \textsc{Case \hyperlink{A}{A}} for small enough $d>0$.

\begin{figure}[h]
\definecolor{verdebosque}{rgb}{0,0.39,0}
\definecolor{rojolava}{rgb}{0.8,0.16,0.16}
\definecolor{violeta}{rgb}{0.62,0.32,0.72}
\begin{subfigure}[b]{0.48\textwidth}
         \centering
         \includegraphics[width=\textwidth]{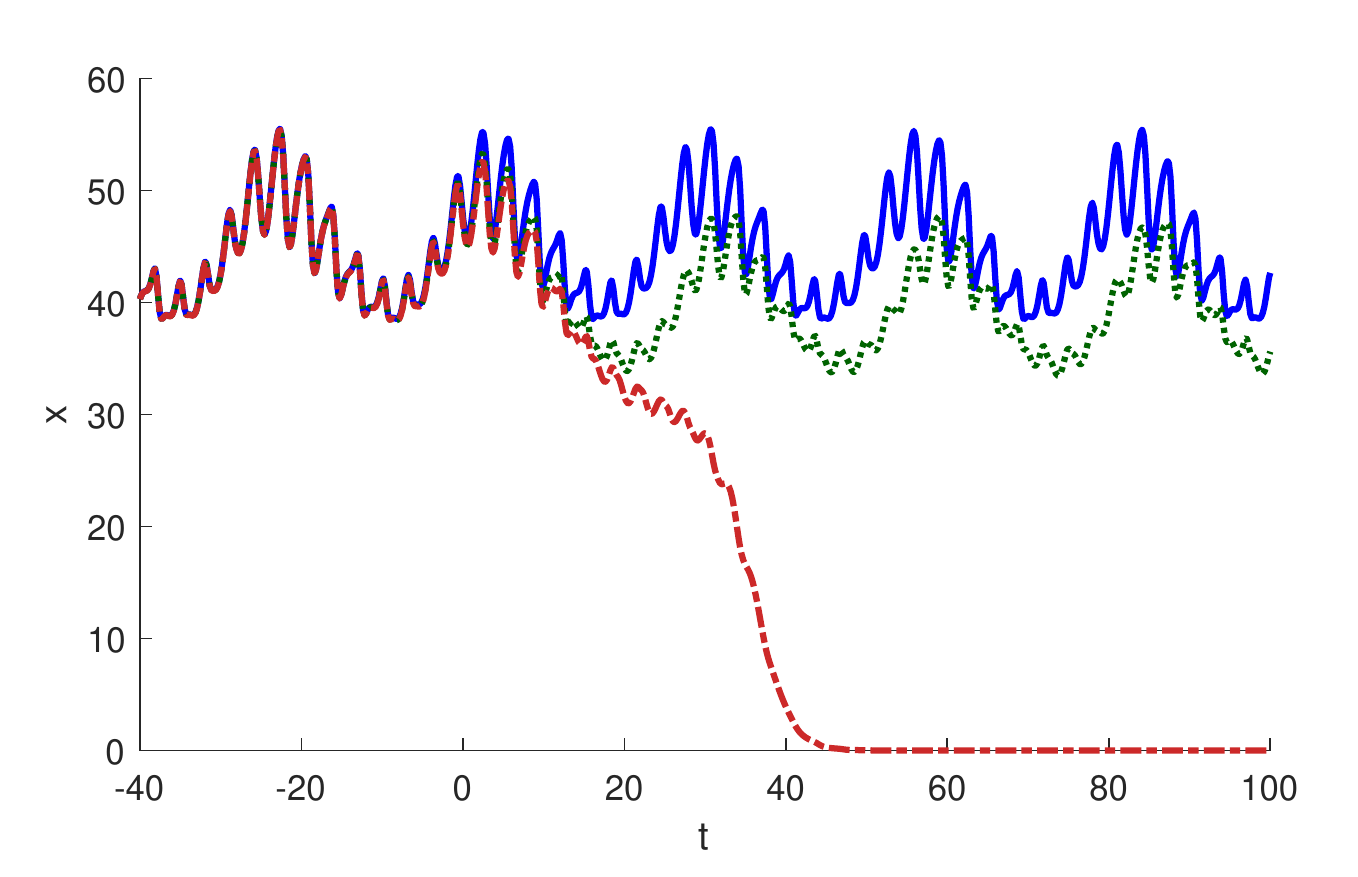}
         \caption{For a population with multiplicative strong Allee effect in the absence of predation, a critical extinction (in dashed-dotted red line) of the native species occurs as the predation increases.}
         \label{fig:strongallee2}
     \end{subfigure}
\hfill
\begin{subfigure}[b]{0.48\textwidth}
         \centering
         \includegraphics[width=\textwidth]{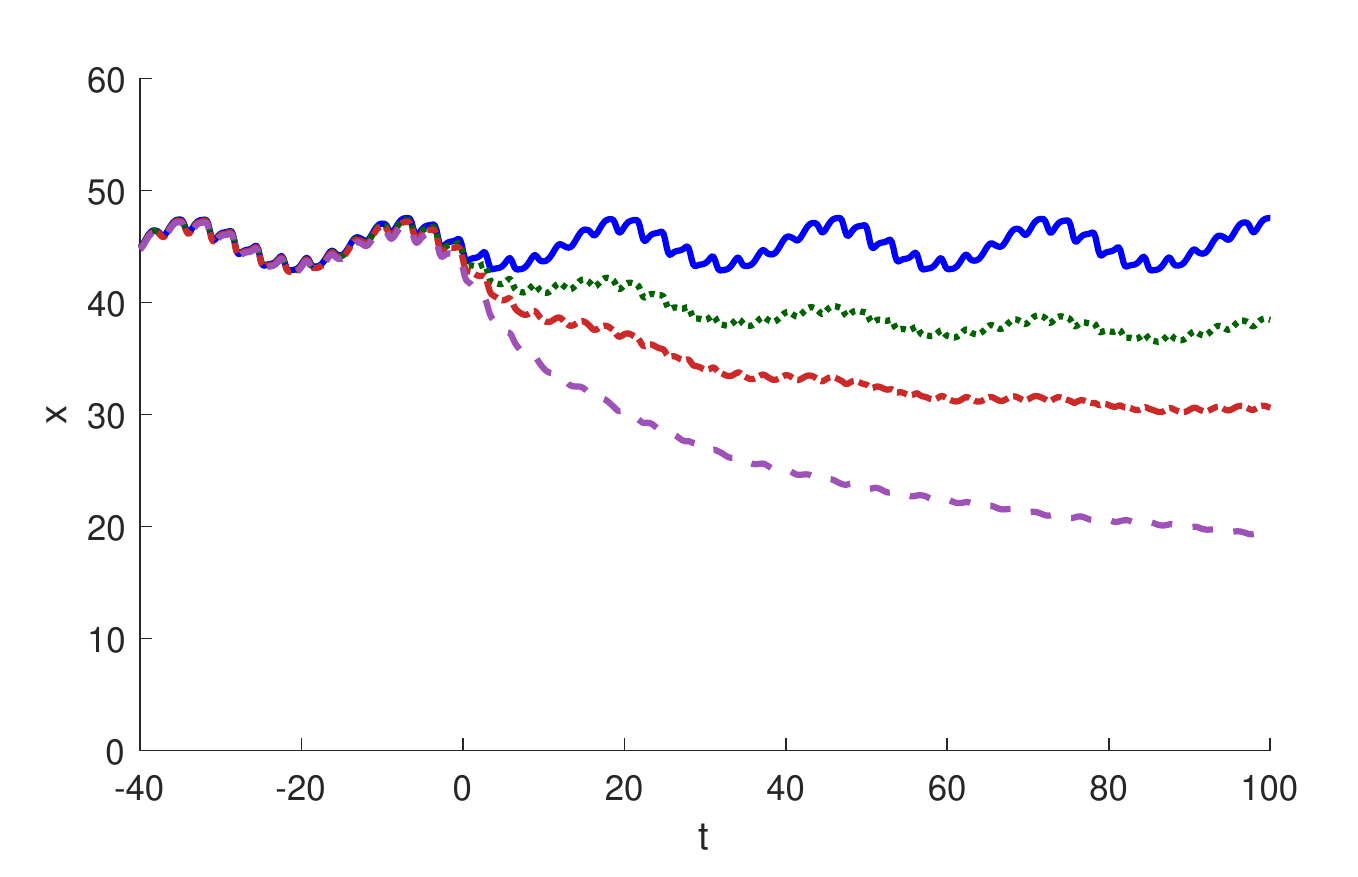}
         \caption{For a population which exhibits multiplicative weak Allee effect in the absence of predation, the native population is preserved even if the predation increases strongly.}
         \label{fig:weakallee2}
     \end{subfigure}
\caption{Extinction or persistence depending on the type of Allee effect in the absence of predation.}
\label{fig:strongweak}
\end{figure}
The two images of Figure~\ref{fig:strongweak} show the upper steady population (i.e., the upper bounded solution) of \eqref{eq:transitionplascubic}$_d$ for several values of $d$.
For Figure~\ref{fig:strongallee2}, we take $r(t)\equiv 1$, $K(t)=30+60\sin^2(t)$, $S(t)=20+20\cos^2(t\sqrt{5}/2)$, and $\beta=800$. Then, the past equation $x'=f(t,x,0)$ exhibits strong Allee effect, since the three hyperbolic solutions are nonnegative (as a simple simulation shows). In addition, $\sup I_\beta> 3.5$.
The blue solid line corresponds to $d=0$, the green dotted line to $d=1.1$, and the red dashed-dotted line to  $d=1.5$. We observe that, as already known, a small predation ensures the persistence of the upper steady population, while a greater predation causes extinction.
For Figure~\ref{fig:weakallee2}, we take $r(t)= 0.01+0.1\cos^2(t\sqrt{5}/2)$, $K(t)=30+60\sin^2(t)$, $S(t)\equiv-0.01$, and $\beta=5\cdot{10}^5$, and numerically check that $x'=f(t,x,0)$ exhibits weak Allee effect (since 0 is hyperbolic repulsive). In addition, $\sup I_\beta>560$. The blue solid line depicts the upper steady population
for $d={0}$, the green dotted line for $d=100$, the red dashed-dotted line for  $d=200$, and the violet dashed line for $d=400$.
Here, we simply observe a smooth decrease of the population as predation increases.

So, there exists at least a tipping value $d_0$ of the predation term between $1.1$ and $1.5$ in the case of strong Allee effect depicted in Figure~\ref{fig:strongallee2} (which is in fact the unique one in $I_\beta$, as we will explain below). In the case of
weak Allee effect case depicted in Figure~\ref{fig:weakallee2}, the numerical simulation shows that \textsc{Case \hyperlink{A}{A}} holds as $d$ increases in $I_\beta$.

The reason of this difference can be found in the underlying nonautonomous bifurcation diagram of minimal sets of \eqref{eq:differentbifurcation}.
As the predation perturbation which has been introduced does not change the hyperbolic character of 0 as the parameter changes, it can be
checked that the arguments leading to Proposition~6.1 and Theorem~6.3(i) and (ii) of \cite{dno2} for the perturbation $\gamma x^2$ also hold for the perturbation $-\gamma x^2/(\beta+x^2)$ of \eqref{eq:differentbifurcation} while $\gamma$ is in $I_\beta$. (Again, we need to construct the hull of $f(t,x,\gamma)$ and to check the hypotheses of these results: see Subsection \ref{subsec:2skewproduct}.) Let us consider the case of Figure~\ref{fig:strongallee2}, with strong Allee effect. The solution 0 of \eqref{eq:differentbifurcation}$_\gamma$ is hyperbolic attractive for all $\gamma\in\mathbb R$, and for $\gamma=0$ there are two more strictly positive hyperbolic solutions. We let $\gamma$ increase in order to find a strictly positive bifurcation value $\gamma_0$ for \eqref{eq:transitionplascubic}: for $\gamma\in(0,\gamma_0)$ there are three nonnegative hyperbolic solutions, while 0 is the unique one for $\gamma>\gamma_0$. (In other words, the two upper hyperbolic solutions collapse at $\gamma_0$, and there are no strictly positive bounded solutions for $\gamma>\gamma_0$: we have a saddle-node bifurcation point at $\gamma_0$.) It turns out that, for our choices, $\gamma_0\in I_\beta$. Let us check that this value of $\gamma_0$ coincides with the unique tipping value $d_0>0$ for \eqref{eq:transitionplascubic} in $I_\beta$. If $0\le d< \gamma_0$, $x'=f(t,x,\gamma)$ has three hyperbolic solutions for all $\gamma\in[0,d]$. Therefore, the pairs $(f,\Gamma_d)$ satisfy all the hypotheses of Theorem~\ref{th:4compacity}, and this result ensures that the dynamics of \eqref{eq:transitionplascubic}$_d$ is in \textsc{Case \hyperlink{A}{A}}.
Now, let us take $d\in(\gamma_0,\sup I_\beta)$: $0$ is the unique nonnegative hyperbolic solution of \eqref{eq:differentbifurcation}$_d$, and it exponentially attracts any solution as time increases. Let us check that this ensures that 0 is the asymptotic limit at $\infty$ of any positive solution \eqref{eq:transitionplascubic}$_d$, which completes the proof of our assertion concerning the uniqueness of $d_0$ in $I_\beta$. Since $0$ is a uniformly exponentially asymptotically stable solution of \eqref{eq:transitionplascubic}$_d$, a contradiction argument ensures the existence of a positive bounded solution separated from 0. This means that the corresponding skewproduct flow on the hull has a positive orbit separated from 0. The $\boldsymbol\upomega$-limit set of this orbit contains a strictly positive minimal set, and according to Lemma \ref{lemm:3union} this minimal set is also minimal for the future equation \eqref{eq:differentbifurcation}$_d$, which is impossible.

Notice that $f_x(t,0,\gamma)=-S(t)r(t)/K(t)$ independently of the value of the parameter, and therefore, in contrast to \cite{remo1}, there can exist a critical transition of the $d$-parametric family \eqref{eq:transitionplascubic} without a modification of the indicator of the strength of strong Allee effect \eqref{eq:5alleeffectstrengthindicator1}.

In the case of weak Allee effect of Figure \ref{fig:weakallee2}, 0 is a repulsive hyperbolic solution of \eqref{eq:transitionplascubic}$_d$ for all $d\in\mathbb{R}$, and therefore there exists a strictly positive attractive hyperbolic solution while $d$ remains in $I_\beta$ (see Proposition 5.3 of \cite{dno1}). Consequently,
the upper minimal set continuously decreases when the parameter increases (see Theorem~6.3(i) of \cite{dno2}).
The continuous variation of the set of bounded solutions precludes the possibility of a critical extinction (while $d$ is in $I_\beta$).

\subsection{Allee effect in the skewproduct formulation}
Let us consider a generic family of equations in the skewproduct formulation described in Subsection~\ref{subsec:2skewproduct},
\begin{equation}\label{eq:skewproductgeneral}
x'=\mathfrak{h}(\omega{\cdot}t,x)
\end{equation}
for $\omega\in\Omega$, where $\mathfrak{h},\mathfrak{h}_x,\mathfrak{h}_{xx}\colon\Omega\rightarrow\mathbb{R}$ exist and are continuous,
$\lim_{|x|\rightarrow\infty}\mathfrak{h}(\omega,x)/x=-\infty$ uniformly on $\Omega$,
$\mathfrak{h}$ is $\mathrm{(SDC)_*}$, and $\mathfrak{h}(\omega,0)=0$ for all $\omega\in\Omega$.
We assume that the compact metric set $\Omega$ is minimal, that is, every orbit is dense.
This type of equation can be obtained by taking $\Omega$ as the hull of a initial function
$h\colon\mathbb{R}\times\mathbb{R}\rightarrow\mathbb{R}$ which satisfies
\hyperlink{h2}{\textbf{h2}}-\hyperlink{h5}{\textbf{h5}}
(see Subsection~\ref{subsec:2skewproduct}) and is recurrent in time (to ensure the minimality of the hull).
In particular, if $h$ is quasiperiodic in $t$, then $\Omega$ is minimal.
Observe that $\mathcal{M}_0=\Omega\times\{0\}$ is always a minimal set for the
skewproduct flow $\tau$ induced by \eqref{eq:skewproductgeneral}
on $\Omega\times\mathbb R$.

Proposition~\ref{prop:sufconditionsalleeeffects} provides simple criteria based on results of \cite{dno1} to ensure that each one of the equations of the family \eqref{eq:skewproductgeneral} represents a model with either weak (in (i)) or strong (in (ii)) Allee effect.
Recall that a minimal set $\mathcal M\subset\Omega\times\mathbb R$
is {\em hyperbolic attractive} (resp. \emph{repulsive}) for the skewproduct flow
if it is uniformly exponentially asymptotically stable at $+\infty$ (resp. $-\infty$).
In this case, $\mathcal M$ is
the graph of a continuous map $\tilde m\colon\Omega\to\mathbb R$, and $t\mapsto \tilde m(\omega{\cdot}t)$ is an
attractive (resp.~repulsive) hyperbolic solution of \eqref{eq:skewproductgeneral}$_\omega$ for any $\omega\in\Omega$
(see e.g.~Subsection 2.4 of \cite{dno1}).

\begin{proposition}\label{prop:sufconditionsalleeeffects} The following statements hold:
\begin{enumerate}[label=\rm{(\roman*)}]
\item {\rm(Weak Allee effect).} If $\mathcal{M}_0$ is a repulsive hyperbolic minimal set, then \eqref{eq:skewproductgeneral} has three different minimal sets $\mathcal{M}_l<\mathcal{M}_0<\mathcal{M}_u$, with $\mathcal{M}_l$ and $\mathcal{M}_u$ hyperbolic attractive.
\item {\rm(Strong Allee effect).} If $\mathcal{M}_0$ is an attractive hyperbolic minimal set and there exists $\rho>0$ such that $0<\mathfrak h(\omega,\rho)$ for all $\omega\in\Omega$, then \eqref{eq:skewproductgeneral} has three different minimal sets
    $\mathcal{M}_0<\mathcal{M}_m<\mathcal{M}_u$, with $\mathcal{M}_m$ hyperbolic repulsive and $\mathcal{M}_u$ hyperbolic attractive.
\end{enumerate}
\end{proposition}
\begin{proof} (i) Proposition~5.3(i) of \cite{dno1} ensures the existence of at least three hyperbolic minimal sets, and Theorem~4.2(ii) of \cite{dno1} (applied to an invariant compact set containing all the minimal ones) ensures that there are exactly three, ordered as $\mathcal{M}_l<\mathcal{M}_0<\mathcal{M}_u$, and with $\mathcal{M}_l$ and $\mathcal{M}_u$ hyperbolic attractive.

(ii) The condition on $\rho$, Theorem 5.1(iii) of \cite{dno1} and the invariance of $0$ ensure that all the orbits $\tau(t,\omega,\rho)$ are bounded. Let us fix one of these orbits and take minimal sets $\mathcal M_u$ and $\mathcal M_m$ respectively contained in its $\boldsymbol\upomega$-limit set and $\boldsymbol\upalpha$-limit set. It is easy to check that $\mathcal M_0<\mathcal M_m<\Omega\times\{\rho\}<\mathcal M_u$ (see e.g.~the proof of Theorem 5.10 of \cite{dno1}, and have in mind that $\mathcal M_0$ is repulsive as time decreases). Hence, there exists three hyperbolic minimal sets, and Theorem~4.2(ii) of \cite{dno1} proves their type of hyperbolicity.
\end{proof}

The remark before Proposition \ref{prop:sufconditionsalleeeffects} shows that, under its hypotheses, all the equations of
the family \eqref{eq:alleecubicskewproduct} exhibit the same type of Allee effect: weak or strong.

Now, we consider the particular case of \eqref{eq:skewproductgeneral} given by the model for multiplicative Allee effect \eqref{eq:alleemultiplicativo} in the skewproduct form,
\begin{equation}\label{eq:alleecubicskewproduct}
x'=r(\omega{\cdot}t)\,x\,\left(1-\frac{x}{K(\omega{\cdot}t)}\right)\frac{x-S(\omega{\cdot}t)}{K(\omega{\cdot}t)}\,,
\quad\omega\in\Omega\,,
\end{equation}
where $r,K,S\colon\Omega\rightarrow\mathbb{R}$ are continuous, $r$ and $K$ are strictly positive, and $S$ satisfies $S(\omega)+K(\omega)\ge 0$  for all $\omega\in\Omega$. We rewrite this family of equations as $x'=\mathfrak h(\omega{\cdot}t,x)$.

Let us describe some clear and sufficient conditions ensuring the hypotheses of Proposition~\ref{prop:sufconditionsalleeeffects} for \eqref{eq:alleecubicskewproduct}. Let $\mathfrak{M}_\mathrm{erg}(\Omega,\sigma)$ be the set of ergodic measures on $\Omega$, and observe that $\mathfrak h_x(\omega,0)=-r(\omega)S(\omega)/K(\omega)$.
As explained in Subsection 2.4 of \cite{dno1}, $\mathcal M_0$ is an attractive (resp.~repulsive) hyperbolic minimal set if and only if $0>\sup\Sigma$ (resp. $0<\inf\Sigma$), where $\Sigma=\{\int_\Omega -r(\omega)S(\omega)/K(\omega)\, dm\colon\,m\in\mathfrak{M}_\mathrm{erg}(\Omega,\sigma)\}$.
In particular, $\mathcal{M}_0$ is hyperbolic repulsive (resp.~attractive) if $S$ is strictly positive (resp.~strictly negative).
The additional condition of $0<\mathfrak h(\omega,\rho)$ for all $\omega\in\Omega$ (for a certain $\rho>0$) is now $(K(\omega)-\rho)(\rho-S(\omega))> 0$ for all $\omega\in\Omega$, and is fulfilled if, for instance, $\sup_{\omega\in\Omega}S(\omega)<\inf_{\omega\in\Omega} K(\omega)$.

To end this section, the following proposition states a relationship between the functions $K$, $S$ and the minimal sets of \eqref{eq:alleecubicskewproduct} which determine the strength and range of action of Allee effect. The results are given in the situations described in Proposition~\ref{prop:sufconditionsalleeeffects}.
\begin{proposition} \label{prop:5ultima}
Let $m$ be an ergodic measure on $(\Omega,\sigma)$, and let the map $\omega\mapsto\beta(\omega)$ be the upper delimiter of the global attractor of the skewproduct flow defined by \eqref{eq:alleecubicskewproduct}.
\begin{enumerate}[label=\rm{(\roman*)}]
\item {\rm (Weak Allee effect)}. Assume that $S(\omega)<0$ for all $\omega\in\Omega$. Then,
\begin{equation}\label{eq:firstassertion-c}
\int_\Omega \frac{r(\omega)}{K(\omega)^2}\,(K(\omega)-\beta(\omega))\,(\beta(\omega)-S(\omega))\, dm=0\,.
\end{equation}
Consequently, either $K$ and $\beta$ are equal constants maps or, for every $\omega\in\Omega$, there exists $\{t_n\}_{n\in\mathbb{Z}}$ with $\lim_{n\rightarrow\pm\infty}t_n=\pm\infty$ such that $t\mapsto K(\omega{\cdot}t)-\beta(\omega{\cdot}t)$ changes sign at $t_n$ for all $n\in\mathbb{N}$.
\item {\rm (Strong Allee effect)}. Under the hypotheses of Proposition {\em \ref{prop:sufconditionsalleeeffects}(ii)}, let $\omega\mapsto\kappa(\omega)$ represent the map whose invariant graph is the middle minimal set $\mathcal{M}_m$. Assume also that
    that $S(\omega_0)<K(\omega_0)$ for a point $\w_0\in\Omega$. Then, $S(\omega)<\rho<K(\omega)$ for all $\omega\in\Omega$, and
    \eqref{eq:firstassertion-c} and
\begin{equation*}
\int_\Omega \frac{r(\omega)}{K(\omega)^2}\,(K(\omega)-\kappa(\omega))\,(\kappa(\omega)-S(\omega))\, dm=0
\end{equation*}
hold.
Consequently, either $K$ and $\beta$ (resp.~$S$ and $\kappa$) are equal constants maps or, for every $\omega\in\Omega$, there exists $\{t_n\}_{n\in\mathbb{Z}}$ with $\lim_{n\rightarrow\pm\infty}t_n=\pm\infty$ such that $t\mapsto K(\omega{\cdot}t)-\beta(\omega{\cdot}t)$
(resp.~$t\mapsto S(\omega{\cdot}t)-\kappa(\omega{\cdot}t)$)) changes sign at $t_n$ for all $n\in\mathbb{N}$.
\end{enumerate}
\end{proposition}
\begin{proof} (i) As explained before the statement of this result, the condition $S(\omega)<0$ for all $\omega\in\Omega$ ensures that the hypothesis of Proposition~\ref{prop:sufconditionsalleeeffects}(i) holds, so that there exist three minimal sets. Theorem~3.1 of \cite{dno2} ensures that $\beta$ is a continuous map. Since $\beta>0$ and $t\mapsto\beta(\omega{\cdot}t)$ solves \eqref{eq:alleecubicskewproduct}$_\omega$,
\begin{equation}\label{eq:5betaK}
\frac{\beta'(\omega)}{\beta(\omega)}=\frac{r(\omega)}{K(\omega)^2}(K(\omega)-\beta(\omega))(\beta(\omega)-S(\omega))
\end{equation}
for all $\omega\in\Omega$ and $t\in\mathbb{R}$, where $\beta'(\w)$ is the derivative of $t\mapsto\beta(\omega{\cdot}t)$ evaluated at $t=0$. Hence, $\beta'$ is also continuous. Since $\beta'/\beta=(\ln\beta)'$ and $\ln\beta$ is bounded, Birkhoff's Ergodic Theorem (see e.g.~Theorem~1 in Section 1.2 of \cite{sinai1}) yields \eqref{eq:firstassertion-c}.

Since $S<0<\beta$ and $r>0$,
$r(\beta-S)/K^2>0$. If $\beta=K$, then \eqref{eq:5betaK} evaluated on $\omega{\cdot}t$ yields $\beta'(\omega{\cdot}t)=0$ for all $\omega\in\Omega$ and $t\in\mathbb R$. So, $\beta$ is constant along any orbit in $\Omega$ and, since $\Omega$ is minimal and $\beta$ is continuous, $\beta$ (and hence $K$) is constant on $\Omega$. If this is not the case, \eqref{eq:firstassertion-c} ensures that there exist open sets of positive $m$-measure $U_+,U_-\subset\Omega$ such that $K(\omega)-\beta(\omega)>0$ for all $\omega\in U_+$ and $K(\omega)-\beta(\omega)<0$ for all $\omega\in U_-$. We fix $\omega\in\Omega$ and deduce from the minimality of $\Omega$ the existence of $\{s^\pm_n\}\uparrow\infty$ such that $\omega{\cdot}s^\pm_n\in U^\pm$ for all $n\in\mathbb{N}$. Consequently, there exists a sequence $\{t_n\}\uparrow\infty$ such that $t\mapsto K(\omega{\cdot}t)-\beta(\omega{\cdot}t)$ changes sign at $t_n$ for all $n\in\mathbb{N}$.
The same argument proves the existence of $\{\tilde t_n\}\downarrow-\infty$ with the same property.

(ii) The three hyperbolic minimal sets provided by Proposition \ref{prop:sufconditionsalleeeffects}(ii) are the graphs of the continuous maps $0<\kappa<\beta$ (see again Theorem~3.1 of \cite{dno2}; and, as seen in its proof, $\kappa<\rho<\beta$). The expression of $\mathfrak h$ and the condition $\mathfrak h(\omega,\rho)>0$ ensure that $\rho$ is between $S(\omega)$ and $K(\omega)$ for all $\omega\in\Omega$. Since these maps are continuous, either
$K<\rho<S$ or $S<\rho<K$, and the first possibility is precluded by the statement. Therefore, $K-\kappa>K-\rho>0$ and $\beta-S>\rho-S>0$. Analogous arguments to those of (i) prove the statements.
\end{proof}
The maps $\beta$ and $\kappa$ of Proposition \ref{prop:5ultima}(ii) respectively represent the carrying capacity of the environment and the critical population size above extinction. They can be naturally used to determine two more indicators of the strength of strong Allee effect:
\begin{equation*}
\inf_{m\in\mathfrak{M}_\mathrm{erg}(\Omega,\sigma)}\int_\Omega\frac{\kappa(\omega)}{\beta(\omega)}\,dm\qquad\text{and}\quad\sup_{m\in\mathfrak{M}_\mathrm{erg}(\Omega,\sigma)}\int_\Omega\frac{\kappa(\omega)}{\beta(\omega)}\,dm\,,
\end{equation*}
which measure the relative position of the upper and middle minimal sets with respect to 0. If the first quantity is close to 1, it indicates that $\kappa$ is very close to $\beta$, so that the strength of strong Allee effect is high: only a population very close to the carrying capacity of the environment can persist; and is the second quantity is close to 0, it indicates that $\kappa$ is much lower than $\beta$, and hence that the strength of the strong Allee effect is low: only very (relatively) small populations get extinct.

\end{document}